\documentclass[a4paper,english]{amsart}

\usepackage[utf8]{inputenc}%caractères accentués
\usepackage[T1]{fontenc}

\usepackage{amsfonts}	% mathbb, mathfrak, mathcal
\usepackage{amsmath}	% substack, displaystyle
\usepackage{amssymb}	% sphericalangle
\usepackage{amsthm}	% newtheorem, proof
\usepackage{bbm}	% fonction indicatrice
\usepackage{mathrsfs}
\usepackage{mathtools}	% coloneqq
\usepackage{stmaryrd}
\usepackage{bigints}
\usepackage{comment}

\usepackage{booktabs}   % provides \toprule and \bottomrule for tables

\usepackage[pdftex, pdfstartview=FitW,colorlinks=true,linkcolor=blue,%
urlcolor=blue,citecolor=blue,
pdfauthor={Timothée Bénard and Weikun He},
pdftitle={Zariski-dense random walks on homogeneous spaces of simple Lie groups},
pdfkeywords={Random walk, Homogeneous space, Slicing theorem, Equidistribution}
]{hyperref}
\usepackage{url}	% tilde dans les url
\usepackage{enumitem}	% choisir la numération dans enumerate
\setenumerate[0]{label=\alph*)}   % default enumerate tag, need enumitem package
\usepackage{tikz} %dessins
\usetikzlibrary{decorations.markings}
\usetikzlibrary{math}
\usetikzlibrary{shapes.geometric}
\newtheorem{thm}{Theorem}[section]
\newtheorem{corollary}[thm]{Corollary}
\newtheorem{lemma}[thm]{Lemma}
\newtheorem{proposition}[thm]{Proposition}

\theoremstyle{definition}
\newtheorem{definition}[thm]{Definition}
\newtheorem{example}[thm]{Example}

\newenvironment{remark}{\par\medskip \noindent \textit{Remark.} \rmfamily}{\medskip}

\usepackage{cleveref}	% Cref
\Crefname{thm}{Theorem}{Theorems}
\Crefname{lemma}{Lemma}{Lemmata}
\Crefname{proposition}{Proposition}{Propositions}
\Crefname{corollary}{Corollary}{Corollaries}
\Crefname{definition}{Definition}{Definitions}
\Crefname{appendix}{Appendix}{Appendices}
\Crefname{figure}{Figure}{Figures}
\Crefname{table}{Table}{Tables}

\newcommand{\Q}{\mathbb{Q}}

\newcommand{\N}{\mathbb{N}}
\newcommand{\Z}{\mathbb{Z}}
\newcommand{\C}{\mathbb{C}}
\newcommand{\R}{\mathbb{R}}

\newcommand{\F}{\mathbb{F}}

\newcommand{\eps}{\varepsilon}
\newcommand{\ut}{\mathbf{t}}
\newcommand{\ur}{\mathbf{r}}
\newcommand{\uk}{\mathbf{k}}
\newcommand{\us}{\mathbf{s}}
\newcommand{\ka}{\mathfrak{a}}

\newcommand{\kg}{\mathfrak{g}}
\newcommand{\ks}{\mathfrak{s}}
\newcommand{\kk}{\mathfrak{k}}
\newcommand{\kh}{\mathfrak{h}}

\newcommand{\kp}{\mathfrak{p}}
\newcommand{\kn}{\mathfrak{n}}
\newcommand{\kz}{\mathfrak{z}}

\newcommand{\so}{\mathfrak{so}}
\newcommand{\kb}{\mathfrak{b}}
\newcommand{\ku}{\mathfrak{u}}

\renewcommand{\angle}{\measuredangle}

\newcommand{\uv}{\underline v}
\newcommand{\uw}{\underline w}

\newcommand{\tnu}{\widetilde{\nu}}
\newcommand{\1}{\mathbbm{1}}

\newcommand{\uj}{{\bf j}}
\newcommand{\uV}{\mathcal{V}}
\newcommand{\uW}{\mathcal{W}}

\newcommand{\tA}{\widetilde{A}}
\newcommand{\tF}{\widetilde{F}}

\renewcommand{\phi}{\varphi}
\newcommand{\cN} {{\mathcal N}}
\newcommand{\cF} {{\mathcal F}}
\newcommand{\cG} {{\mathcal G}}
\newcommand{\cH} {{\mathcal H}}
\newcommand{\cL} {{\mathcal L}}
\newcommand{\cA} {{\mathcal A}}
\newcommand{\cB} {{\mathcal B}}
\newcommand{\cI} {{\mathcal I}}
\newcommand{\cJ} {{\mathcal J}}
\newcommand{\cV} {{\mathcal V}}

\newcommand{\cP} {{\mathcal P}}

\newcommand{\cD} {{\mathcal D}}
\newcommand{\cQ} {{\mathcal Q}}
\newcommand{\cE} {{\mathcal E}}

\newcommand{\cS} {{\mathcal S}}
\newcommand{\cT} {{\mathcal T}}
\newcommand{\cK} {{\mathcal K}}
\newcommand{\cR} {{\mathcal R}}
\newcommand{\cW} {{\mathcal W}}
\newcommand{\cO}{\mathcal{O}}

\newcommand{\sD} {{\mathscr D}}
\newcommand{\sC}{\mathscr{C}}
\newcommand{\sP}{\mathscr{P}}

\newcommand{\sS}{\mathscr{S}}

\newcommand{\bL}{\mathbf{L}}

\newcommand{\bbP}{\mathbb{P}}

\newcommand{\data}{\diamondsuit}

\DeclareMathOperator{\dimH}{\dim_{\mathrm{H}}}
\DeclareMathOperator{\Ad}{Ad}
\DeclareMathOperator{\Aut}{Aut}
\DeclareMathOperator{\End}{End}

\DeclareMathOperator{\GL}{GL}

\DeclareMathOperator{\Kill}{Kill}

\DeclareMathOperator{\Zar}{Zar}

\DeclareMathOperator{\SL}{SL}
\DeclareMathOperator{\Sym}{Sym}
\DeclareMathOperator{\SO}{SO}
\DeclareMathOperator{\OO}{O}

\DeclareMathOperator{\ad}{ad}

\DeclareMathOperator{\diag}{diag}
\DeclareMathOperator{\diam}{diam}
\DeclareMathOperator{\inj}{inj}

\DeclareMathOperator{\supp}{supp}
\DeclareMathOperator{\Span}{Span}
\DeclareMathOperator{\tr}{tr}
\renewcommand{\sl}{\mathfrak{sl}}
\DeclareMathOperator{\dist}{d}
\DeclareMathOperator{\dang}{\dist_{\measuredangle}}

\DeclareMathOperator{\SubP}{(S^-)}
\DeclareMathOperator{\SAP}{(S^+A)}
\DeclareMathOperator{\Obs}{Obs}

\newcommand{\abs}[1]{\lvert#1\rvert}    % valeur absolue
\newcommand{\abse}[1]{\left\lvert#1\right\rvert}
\newcommand{\absbig}[1]{\bigl\lvert#1\bigr\rvert}

\newcommand{\norm}[1]{\lVert#1\rVert}   % norme

\newcommand{\set}[1]{\{\, #1  \,\}}     % ensemble
\newcommand{\setbig}[1]{\bigl\{\, #1 \,\bigr\}}
\newcommand{\setBig}[1]{\Bigl\{\, #1 \,\Bigr\}}

\newcommand{\dd}{\,\mathrm{d}}
\newcommand{\defeq}{\coloneqq}
\DeclareMathOperator{\Id}{Id}
\newcommand{\acts}{\curvearrowright}

\newcommand{\transp}[1]{\prescript{t}{}{#1}}

\DeclareMathOperator{\Gr}{Gr}
\DeclareMathOperator{\Ker}{Ker}

\DeclareMathOperator{\leb}{Leb}

\DeclareMathOperator{\Vect}{Span}
\DeclareMathOperator{\too}{{}_{to}}

\renewcommand{\setminus}{\smallsetminus}
\renewcommand{\emptyset}{\varnothing}       %\varnothing est plus joli
\renewcommand{\subset}{\subseteq}

\newcommand{\increasing}{\boxslash}         %{\pmb{\boxslash}}

\newcommand{\bfparagraph}[1]{\par\medskip \noindent{\bfseries #1}}

\title[Random walks on simple homogeneous spaces]{Effective equidistribution of random walks on simple homogeneous spaces}
\author{Timoth\'ee B\'enard }
\address{CNRS – LAGA, Universit\'e Sorbonne Paris Nord, 99 avenue J.-B.
Cl\'ement, 93430 Villetaneuse}
\email{benard@math.univ-paris13.fr}

\author{Weikun He}
\address{State Key Laboratory of Mathematical Sciences, Academy of Mathematics and System Science, Chinese Academy of Sciences, Beijing 100190, China}
\email{heweikun@amss.ac.cn}
\thanks{W.H. is supported by the National Key R\&D Program of China (No. 2022YFA1007500) and the National Natural Science Foundation of China (No. 12288201).}
\date{}

\begin{document}
\large

\begin{abstract}
We consider a random walk on a homogeneous space $G/\Lambda$ where $G$ is a non-compact simple Lie group and  $\Lambda$ is a lattice. The walk is driven by a probability measure $\mu$ on $G$ whose support generates a Zariski-dense subgroup.
We show that the random walk equidistributes toward the Haar measure unless it is trapped in a finite $\mu$-invariant set.  Moreover, under arithmetic assumptions on the pair $(\Lambda, \mu)$, we show the convergence occurs at an exponential rate, tempered by the obstructions that the starting point may be high in a cusp or close to a finite orbit.

The main challenge is to show that the dimensional properties of a given probability distribution on $G/\Lambda$ improve under convolution by $\mu$.
For this, we develop a new method, which combines a dimensional interpolation result and a dimensional increase alternative.
This approach allows us to bypass inherent geometric obstructions.
To show dimensional interpolation, we establish a general subcritical projection theorem under optimal non-concentration assumptions on the projector, and a corresponding submodular inequality for irreducible representations which allows its application to random walks.
Both are of independent interest.
The dimensional increase alternative aligns with the spirit of Bourgain's projection theorem.
It is fine-tuned for random walks and has the advantage of being valid in situations lacking transversality.
\end{abstract}
\maketitle

\tableofcontents

\section{Introduction} \label{Sec-intro}

\subsection{Main results} \label{Sec-main-results}

Let $G$ be  a non-compact connected real Lie group with finite center and simple Lie algebra $\kg$, let $\Lambda\subseteq G$ be a lattice in $G$.
Let $X=G/\Lambda$ be the quotient space, and $m_{X}$ the unique $G$-invariant Borel probability measure on $X$, also called the Haar measure.

Given a Borel probability measure $\mu$ on $G$, we consider the Markov chain on $X$ with transitional probability distributions $(\mu*\delta_{x})_{x\in X}$ where $*$ denotes the convolution and $\delta_{x}$ the Dirac mass at $x\in X$. Given some initial probability distribution $\nu$ on $X$, we are interested in the asymptotic of the $n$-th step distribution of this Markov chain, in other words, $\mu^n*\nu$ where $\mu^n$ stands for the $n$-fold convolution power of $\mu$. We  show that under natural necessary constraints over $\nu$, the distribution $\mu^n*\nu$ is close to $m_{X}$ for large $n$. Our results include both qualitative and effective estimates.

%Let $G$ be  a connected real linear group with simple Lie algebra $\kg$, let $\Lambda\subseteq G$ be a lattice in $G$. Set $X=G/\Lambda$ the quotient space. Given a probability measure $\mu$ on $G$, we consider the Markov chain on $X$ with transitional probability distributions $(\mu*\delta_{x})_{x\in X}$ where $*$ denotes the convolution and $\delta_{x}$ the Dirac mass at $x\in X$. We are interested in the asymptotic of the $n$-step distribution of this Markov chain. More precisely in other terms, given $x\in X$, we  describe the behavior of $\mu^n*\delta_{x}$ for large $n$, where $\mu^n$ stands for the $n$-fold convolution power of $\mu$.

We work under the condition that $\mu$ has a \emph{finite exponential moment}, i.e. there exists $\eps>0$ such that
$$\int_{G}\|\Ad g\|^\eps\,\dd\mu(g) < +\infty $$
where $\Ad:G\rightarrow \Aut(\kg)$ stands for the adjoint representation, and $\norm{\cdot}$ is any norm on $\kg$. We denote by $\Gamma_{\mu}$ the subgroup of $G$ generated by the support of $\mu$. We assume that $\mu$ is \emph{Zariski-dense}, by which we mean that $\Ad(\Gamma_{\mu})$ is   Zariski-dense  in $\Ad(G)$.

\bfparagraph{Equidistribution in law.}
We first present   our main qualitative result: the $n$-step distribution of a Zariski-dense random walk on $X$ either converges toward the Haar measure or is trapped in a finite orbit.

\begin{thm}[Equidistribution in law] 
\label{thm:equidis}
 Let $G$ be a non-compact connected real Lie group with finite center and simple Lie algebra, let $\Lambda$ be a lattice in $G$, set $X=G/\Lambda$.
Let $\mu$ be a Zariski-dense probability measure on $G$ with a finite exponential moment.
For every $x \in X$, we have
\begin{equation}\label{eqth11}
\mu^{n}* \delta_x \rightharpoonup^* m_X
\end{equation}
unless the orbit $\Gamma_{\!\mu} x$ is finite.
\end{thm}

In the above statement, the symbol $\rightharpoonup^*$ refers to weak-$*$ convergence of measures, in other terms \eqref{eqth11} means that $\lim_{n\to +\infty}\mu^{*n}* \delta_x (f)= m_X(f)$ for every continuous bounded function $f:X\rightarrow \R$.

\begin{remark}
If $\Gamma_{\!\mu} x$ is finite, then $\mu^{n}* \delta_x$ converges to the uniform probability measure on $\Gamma_{\mu} x$ provided $\mu$ is aperiodic. By \Cref{thm:equidis},  aperiodicity is only necessary  for equidistribution within finite orbits.
\end{remark}

\Cref{thm:equidis} can be meaningfully compared with the work of Benoist-Quint \cite{BQ1, BQ3}. In \cite{BQ3}, Benoist and Quint obtain \eqref{eqth11} in \emph{Ces\`aro-average}, that is, for all $x\in X$ with $\Gamma_{\mu}x$ infinite, they prove that
\begin{equation}\label{eqth11-Ces\`aro}
\frac{1}{n}\sum_{k=0}^{n-1}\mu^{k}* \delta_x \rightharpoonup^* m_X.
\end{equation}
Their proof consists in showing  that convergence \eqref{eqth11-Ces\`aro} is equivalent to the rigidity of stationary measures. The latter is the main result of their preceding paper \cite{BQ1}, and relies on their celebrated  exponential drift argument as well as   Ratner's equidistribution theorems for unipotent flows. As part of \cite[Question 3]{BQ-survey12},  Benoist and Quint  ask whether the Ces\`aro average in \eqref{eqth11-Ces\`aro} can be removed. \Cref{thm:equidis} answers this question positively.\footnote{Note however that the question of removing the Cesàro average is still open in the broader context of $\Ad$-semisimple random walks on homogeneous spaces, see \cite[Question 3]{BQ-survey12} for more details.}
%\footnote{Note however that the question of removing the Ces\`aro average is still open in the broader context where $G$ is a general connected real Lie group, and $\mu$ has finite exponential moment and is Zariski-dense in a semisimple Lie subgroup of $G$ with no compact factor, see  \cite[Question 3]{BQ-survey12} for more details.}

Related works on that question comprise  \cite{BFLM, He2020IJM, HS2022, HS2023, HLL_Aff, HLL_Nil} in the setting of nilmanifolds,  \cite{Benard23-equidmass} for symmetric random walks, and \cite{KhalilLuethi, BHZ24, Khal-Lue-Wei25} in the context of upper triangular random walks. Our previous paper \cite{BH24} also tackles the case where $G$ is $\SO(2,1)$ or $\SO(3,1)$.

The proof of \Cref{thm:equidis} is disjoint from the work of Benoist-Quint. It does not use exponential drift nor Ratner's  theorems, and yields, in particular,  a new proof of the rigidity of stationary measures from \cite{BQ1} in the simple Lie group case. \Cref{thm:equidis} will in fact be a consequence of a quantitative equidistribution estimate which we now present.

\bigskip
\bfparagraph{Effective equidistribution.}
We  present our main effective estimate. Given an initial distribution $\nu$ on $X$ that is not too concentrated near infinity and has positive dimension, we show that $\mu^n*\nu$ converges toward the Haar measure with an exponential rate.

\bigskip
To quantify a rate of convergence, we need a class of regular functions. For that, we fix a right $G$-invariant Riemannian metric on $G$, and equip $X$ with the quotient metric. For $\beta\in (0, 1]$, we let $C^{0,\beta}(X)$ denote the space of bounded $\beta$-Hölder continuous functions on $X$, endowed with its usual norm $\norm{\,\cdot\,}_{C^{0,\beta}}$ :
\begin{align} \label{def-holder}
\forall f \in C^{0,\beta}(X), \quad \norm{f}_{C^{0,\beta}} \defeq \norm{f}_\infty + \sup_{x\neq y \in X} \frac{\abs{f(x)-f(y)}}{\dist(x,y)^\beta}.
\end{align}
The corresponding Wasserstein distance between two probability measures $\nu,\nu'$ on
X is defined as
\begin{align} \label{def-wassertstein}
\cW_\beta(\nu,\nu') \defeq \sup_{f \in C^{0,\beta}(X),\norm{f}_{C^{0,\beta}} \leq 1} \abse{\int_X f \dd \nu - \int_X f \dd \nu'}.
\end{align}
We show

\begin{thm}[Effective equidistribution I] \label{Eff-eq-thmI}
Let $G$ be a non-compact connected real Lie group with finite center and simple Lie algebra. Let $\Lambda\subseteq G$ be a lattice, $X=G/\Lambda$ equipped with a quotient right $G$-invariant Riemannian metric. Let $\mu$ be a Zariski-dense probability measure on $G$ with finite exponential moment.

Given $\beta \in {(0,1]}$ and $\kappa \in (0, 1]$, there exists $\eps = \eps(X,\mu,\beta,\kappa) > 0$  such that for small  enough $\delta>0$, the following holds.

Let $\nu$ be a probability measure on $X$ satisfying
\[\nu (B_\rho(x)) \leq \rho^\kappa \text{ for all } x \in X, \rho \in [\delta, \delta^\eps].\]
Then for all $n\geq  |\log \delta|$, one has
\[\cW_{\beta}(\mu^{n}*\nu, \,m_{X})\leq \delta^\eps+2\nu\{\inj\leq \delta^\eps\}\]
where  $m_{X}$ denotes the Haar probability measure on $X$.
\end{thm}
\bigskip

\bfparagraph{Effective equidistribution under arithmetic assumptions.}
It is natural to ask about an effective convergence rate when the initial distribution is a deterministic point, i.e., $\nu=\delta_{x}$ for some $x\in X$. We obtain such result under arithmetic assumptions, namely if $\Lambda$ is an \emph{arithmetic} lattice in $G$, and \emph{$\mu$ is algebraic with respect to $\Lambda$}. This condition on $\mu$ means that  $\Ad(\Gamma_{\!\mu})$ and $\Ad(\Lambda)$ have algebraic entries with respect to some fixed basis of $\kg$.

Note that for a deterministic starting point $x$, there are two obstructions that can delay (or prevent) equidistribution within $X$. First, $x$ may be very far in a cusp. Second, $x$ may be close to (or within) a finite $\Gamma_{\mu}$-orbit.
To quantify those, we introduce $x_{0} \defeq \Lambda/\Lambda\in X$ which we see as a basepoint for $X$, as well as
\[W_{\!\mu, R} \defeq \set{x \in X : |\Gamma_{\!\mu}x| \leq R},\]
the set of points whose $\Gamma_{\!\mu}$-orbit is finite of cardinality at most $R>0$.

\begin{thm}[Effective equidistribution II] \label{Eff-eq-thmII}
Let $G$ be a non-compact connected real Lie group with finite center and simple Lie algebra. Let
$\Lambda$ be an arithmetic lattice in $G$, set $X=G/\Lambda$ equipped with a quotient right $G$-invariant  Riemannian metric. Let $\mu$ be a Zariski-dense finitely supported probability measure on $G$ which is algebraic with respect to $\Lambda$.

Given $\beta\in (0, 1]$, there exists a constant $A=A(X, \mu, \beta)>1$  such that for all  $x\in X$, $n\in\N$ and $R \geq 2$, we have
\[\cW_{\beta}(\mu^{n}*\delta_{x}, \,m_{X}) \leq R^{-1} \]
as soon as $n \geq A\log R + A \max \{\abs{\log \dist(x,W_{\!\mu, R^A})},\, \dist(x,x_{0})\}$.
\end{thm}

%\begin{thm}[Effective equidistribution II] \label{main-thm}
%Let $G$ be a connected real linear group with simple Lie algebra. Let
%$\Lambda$ be an arithmetic lattice in $G$, set $X=G/\Lambda$ equipped with a quotient right $G$-invariant  Riemannian metric. Let $\mu$ be a Zariski-dense finitely supported probability measure on $G$ which is algebraic with respect to $\Lambda$.

%There exists a constant $A>0$  such that for all  $x\in X$,  $n\in\N$, $R \geq 2$ and  $f\in \Lip(X)$, we have
%\[|\mu^{*n}*\delta_{x}(f)-m_{X}(f)| \leq R^{-1} \|f\|_{\Lip} \]
%as soon as $n \geq A\log R + A \max \{|\log \dist(x,W_{\!\mu, R^A})|,\, \dist(x,x_{0})\}$.
%\end{thm}

\begin{remark}
In the case where $W_{\!\mu, R}=\emptyset$, we use the convention that
\[\max \{\abs{\log \dist(x,W_{\!\mu, R})},\, \dist(x,x_{0}) \}=\dist(x,x_{0}).\]
\end{remark}

Theorems \ref{Eff-eq-thmI} and \ref{Eff-eq-thmII} are connected to a vast corpus of research  dedicated to quantify equidistribution on homogeneous spaces. Most relevant to us are the works of Bourgain-Furman-Lindenstrauss-Mozes about the torus case \cite{BFLM}, and its extensions \cite{He2020IJM, HS2022, HS2023, HLL_Aff, HLL_Nil}; the works of Bourgain-Gamburd \cite{BG2008, BG2012}, Benoist-Saxcé \cite{BenoistSaxce} in the context of compact Lie groups; and the works of Kim \cite{Kim24}, Lindenstrauss, Mohammadi, Wang,  Yang \cite{LM, LMW, Yang, LMWY25},  Lin \cite{Lin25, Lin25b} for unipotent flows. Our previous paper \cite{BH24} tackles the case where $G$ is isogenous to either $\SO(2,1)$ or $\SO(3,1)$.
All these works share the feature that they crucially boil down to a dimensional bootstrap, which in turn relies on the iteration of a projection theorem or a sum product phenomenon.
Their bootstrap implementations rely on the specific structure of the ambient group— a torus, a torus fibration over a well-understood base, a compact group, or a relatively small group such as $\SL_2(\C)$ or $\SL_3(\R)$.
We will develop a bootstrap method which applies to \emph{all non-compact simple Lie groups}, regardless of dimension, rank, or other structural complexity. This general method is likely applicable in other contexts as well.

\begin{remark}
Arithmetic restrictions as in \Cref{Eff-eq-thmII}  also appear  in the aforementioned results  (e.g. algebraic entries in \cite{BFLM, BG2008, BG2012, BenoistSaxce}, artihmetic lattice in \cite{LM, LMW, Yang, LMWY25}). Getting rid of such assumptions is a well-known open question. However, it does not concern the dimensional bootstrap phase, but rather a preliminary phase where some positive initial dimension is obtained, see \S\ref{Sec-on-the-proof}. The present paper focuses on the bootstrap phase, and we leave to other works the non-arithmetic refinements of the preliminary phase.
\end{remark}

%Their implementation of the bootstrap scheme  relies crucially on how nice the ambient group is: either a torus, a torus fibration over a well-known space, a compact group, or a small group such as $\SL_{2}(\C)$ or $\SL_{3}(\R)$. The strength of our results is that they tackle all non-compact simple Lie groups, may they have large dimension,  be of higher rank, etc.

\bigskip
We record  two meaningful corollaries of \Cref{Eff-eq-thmII}. First, we identify starting points with exponential rate of convergence.
These are precisely the points which are not too well approximated by small finite $\Gamma_{\mu}$-orbits.
Given $D>1$, say $x\in X$ is \emph{$(\mu, D)$-Diophantine}  if for all $R>1$ with  $W_{\!\mu, R}\neq \emptyset$, one has
\[\dist(x, W_{\!\mu, R}) \geq \frac{1}{D} R^{-D}.\]
Observe this condition gets weaker as $D\to +\infty$. Say $x$ is \emph{$\mu$-Diophantine generic}  if it is $(\mu, D)$-Diophantine for some $D$. The set of $\mu$-Diophantine generic points $x\in X$ has full $m_X$\nobreakdash-measure. It is equal to $X$ when $\Gamma_{\!\mu}$ has no finite orbit.

\begin{corollary}[Points with exponential rate of equidistribution]
\label{cr:dioph}
In the setting of \Cref{Eff-eq-thmII}, let $\beta\in (0, 1]$, $x\in X$. The  following are equivalent:
\begin{enumerate}
\item \label{it:exprate1}The point $x$ is $\mu$-Diophantine generic.
\item \label{it:exprate2} There exists $C, \theta>0$ such that for every $n\geq 1$, $f\in C^{0, \beta}(X)$,
\begin{equation} \label{exp-eq}
|\mu^{n}*\delta_{x}(f)-m_{X}(f)| \leq \norm{f}_{C^{0,\beta}} C e^{-\theta n}.
\end{equation}

\end{enumerate}
Moreover, the constants $(C,\theta)$ can be chosen uniformly when $x$ varies in a compact subset and is $(\mu,D)$-Diophantine for a fixed $D$.
\end{corollary}

 We also derive effective equidistribution of large finite $\Gamma_{\mu}$-orbits, with polynomial rate in the cardinality of the orbit. Here we use  that    \Cref{Eff-eq-thmII} does \emph{not} require $x$ to have infinite $\Gamma_{\mu}$-orbit (contrary  to \Cref{thm:equidis}).

\begin{corollary}[Polynomial equidistribution of finite orbits]  \label{cor-orbites}
In the setting of \Cref{Eff-eq-thmII},
let $Y\subseteq X$ be a finite  $\Gamma_{\!\mu}$-orbit of cardinality $R$.
Let $m_{Y}$ denote the uniform probability  measure on $Y$.
Then for all $\beta\in (0, 1]$,  $f\in C^{0,\beta}(X)$, one has
\[|m_{Y}(f)-m_{X}(f)| \leq \|f\|_{C^{0,\beta}} C R^{-c}\]
where $C,c>0$ depend only  on $X$, $\mu$, $\beta$.
\end{corollary}

\Cref{cor-orbites} is an effective upgrade of \cite[Corollary 1.8]{BQ3}.
In the case where $\Gamma_{\mu}$ is a lattice, the result can be deduced from Maucourant-Gorodnik-Oh \cite[Corollary 3.31]{GMO08} about the effective equidistribution of Hecke points, see also \cite{Clozel-Oh-Ullmo}. In the context of unipotent flows on small groups, polynomial equidistribution of large periodic  orbits   follows from \cite{LMW,Yang, LMWY25}.

\subsection{About the proofs} \label{Sec-on-the-proof}
The proofs of the aforementioned  theorems consist of three phases: (Phase I) Starting from a point with infinite orbit, the random walk generates positive dimension above a given scale. This phase is unnecessary in the setting of \Cref{Eff-eq-thmI}, it is completed with a rate for that of \Cref{Eff-eq-thmII}, and no rate for that of \Cref{thm:equidis}. (Phase II) Starting from a measure with positive dimension above a scale,  the random walk bootstraps the dimension arbitrarily close to that of $X$.  (Phase III) Once high dimension is known, effective equidistribution follows by smoothing and a spectral gap argument.

This three-phase philosophy is shared by many works (e.g. \cite{BFLM, BG2008, BenoistSaxce, LM, LMW, Yang, LMWY25, BH24}). In our setting, Phases I and III have already been completed in \cite{BH24}.
For $G=\SO(2,1)$ or $\SO(3,1)$, Phase II was carried out in \cite{BH24} using an extension of Bourgain's projection theorem, which took the form of a multislicing estimate \cite[Corollary 2.2]{BH24}. 
%However this estimate is not sufficient to tackle a general simple Lie group.
%Applying this estimate was only made possible by the restriction on the ambient group $G$.
However, the applicability of this estimate relied crucially on the restriction imposed on the ambient group $G$.

\bfparagraph{The challenge.}
In this paper, we take up the challenge of extending Phase II of \cite{BH24}  to \emph{all non-compact simple Lie groups}.
The major obstacle is that the projection theorems \`a la Bourgain that are currently available~\cite{He2020JFG} require a strong form of non-concentration to be applicable.
To be precise, we say a subspace $V\subseteq \kg$ is \emph{transverse} if for any subspace $W\subseteq \kg$ with $\dim V+\dim W=\dim \kg$, there exists $g\in G$ such that
$$\Ad(g)V \cap W= \{0\}.$$
For random walks on $\SO(2,1)$ or $\SO(3,1)$, the  available projection theorems \`a la Bourgain require  that the highest weight subspace  of a maximal torus acting on $\kg$ via the adjoint representation be transverse. This is trivial for $\SO(2,1)$, and it has been checked for $\SO(3,1)$ in \cite{He2020IJM}.
For an arbitrary non-compact simple Lie group $G$, the corresponding transversality requirement concerns a much broader class of subspaces, namely all those of the form $E_{v,t}:=\oplus_{\alpha\,:\,\alpha(v)\geq t}\, \kg_{\alpha}$ where $v$ is an element of the Cartan subspace $\ka$, $t \in \R$ and $\kg_{\alpha}$ denotes the restricted root space\footnote{We allow $\alpha$ to be $0$, in which case $\kg_{0}$ is the centralizer of $\ka$ in $\kg$.} of root $\alpha$ relative to $\ka$.
Unfortunately, beyond small groups, the subspaces $E_{v,t}$'s are usually \emph{not} transverse.
For example, in the case $G=\SO(n, 1)$ with $n\geq 7$ odd, transversality already fails for the highest weight subspace $\kg^+$ of a maximal torus. Even worse: any four translates of $\kg^+$ under $\Ad(G)$ are never mutually in direct sum  (see \Cref{lack}), even though their dimension is much smaller than that of $\kg$. For $G=\SL_{3}(\R)$, the space $\kg^+$ is one-dimensional, whence transverse by irreducibility of $\Ad(G)\acts \kg$. However, transversality fails for some other  subspaces $E_{v,t}$, for instance the two subspaces of $\sl_{3}(\R)$ given by
$$\kb=\begin{pmatrix} *&*&* \\ 0 &*&*\\0&0&*\end{pmatrix} \quad \text{ and }  \quad  W=\left\{\begin{pmatrix} t&0&0 \\0&t&0\\ * &*&-2t\end{pmatrix}\,:\, t\in \R \right\}$$
satisfy\footnote{
Indeed, write $P$ (resp $P^-$) the upper (resp. lower) triangular subgroup of $\SL_{3}(\R)$.  Note $\kb$ and $W$ are respectively invariant under $P$ and $P^-$. As they also have nontrivial intersection, we get $\Ad(g)\kb\cap W\neq \{0\}$ for all $g\in P^-P$. But $P^-P$ is Zariski-dense in $G$ and the nontrivial  intersection condition is Zariski-closed. Hence $\Ad(g)\kb\cap W\neq \{0\}$ for all $g\in G$.
} $\dim \kb + \dim W = \dim \sl_3(\R)$ and $\Ad(g)\kb\cap W\neq \{0\}$ for all $g\in G$.

These obstructions call for the development of projection theorems and multislicing estimates with less stringent non-concentration hypotheses that would authorize application to random walks.
This is the task that we will pursue in this paper.

\bfparagraph{A subcritical projection theorem under optimal assumption.}
We  first establish a subcritical projection theorem (i.e., with small dimensional loss instead of a dimensional gain), and a corresponding subcritical multislicing estimate, both under optimal non-concentration assumptions. These provide a vast generalization of arguments in  \cite{BHZ25} and \cite{Lin25} which manage to obtain subcritical estimates despite an apparent lack of transversality. In fact \cite{BHZ25} and \cite{Lin25} rely on a combinatorial trick which exploits the specificity of their framework to reduce to the transverse case. Such a strategy seems hopeless in the context of a general simple Lie group. We use a different approach, which relies on effective upper bounds for Brascamp-Lieb constants. The output is a very general method to obtain subcritical estimates, which has applications for walks on homogeneous spaces  and certainly beyond that. 

\bfparagraph{A submodular inequality for the Grassmannian of a linear representation.}
%The weak non-concentration property required for the subcritical regime boils down, in the context of random walks, to a property of partial transversality for the random Oseledets subspaces in $\kg$, see \Cref{avoid-bad-pencil-proba}. This  in turn reduces to a beautiful submodular property regarding Borel invariant subspaces in simple complex Lie algebras. This inequality is of independent interest.
The weak non-concentration property required to apply the previous subcritical estimates in the context of random walks boils down to a property of partial transversality for the random Oseledets subspaces in $\kg$, stated as \Cref{avoid-bad-pencil-proba}. To show such transversality, we establish a beautiful submodular inequality regarding Borel invariant subspaces in simple complex Lie algebras. This inequality is new and of independent interest.
It is presented in a more general form, Theorems \ref{submod-G-F}, \ref{submod-G-F-2}, in the context of linear representations  (see also \cite{EMS97, Lin25b}).

%to an upper bound on the dimension of minimal intersections of subspaces taken from two $G$-orbits of Oseledets subspaces in $\kg$, \Cref{avoidW'}.

\bfparagraph{The supercritical regime.}
Subcritical estimates are not  enough just yet to perform the dimensional bootstrap, as they only guarantee a small dimensional loss, instead of a small gain.
If the highest weight subspace is transverse in the sense defined previously, then this dimensional gain can be obtained by means of Bourgain's projection theorem \cite{Bourgain2010} or its generalization in higher rank \cite{He2020JFG}.
Unfortunately, as discussed above, transversality may fail, even for the highest weight direction.
For that reason, we  also promote a supercritical multislicing decomposition, which is  looser than the original supercritical theorem from \cite[Theorem 2.1]{BH24}, and motivated by our application to random walks in \Cref{Sec-RWincrease}. Indeed it allows us to bypass the obstructions mentioned above by  exploiting only a weak form of transversality for the subspaces $E_{v,t}$ (namely \Cref{nc-highest-weight}) and still obtain the desired dimensional increment.
Note we do not establish an improved general supercritical projection theorem, though pursuing this direction would certainly be of interest.

\bigskip

\subsection{Conventions and notations} \label{conventions-notations}

The cardinality of a finite set $A$ is denoted by $|A|$.
The neutral element of a group is denoted by $\Id$. We write $\R^+, \N, \N^*$ for  the sets of non-negative real numbers, non-negative integers, and positive integers.

\bfparagraph{Metric spaces.}
Given a metric space $X$, and $\rho > 0$, we denote by $B^X_\rho(x)$ the open ball of radius $\rho$ and center $x$. 
If the metric space in which $x$ is taken is clear from context, we may simply write $B_{\rho}(x)$.
If the space has a distinguished point (say the zero vector $0$ of a vector space, or the neutral element $\Id$ of a group), then $B^X_\rho$ refers to the ball centered at the distinguished point.
For example, taking $X=G/\Lambda$, with $G$ equipped with a right $G$-invariant metric, and $X$ with the quotient metric, we have $B^X_{\rho}(x)=B^G_{\rho}x$. In this context, we also set
$$\inj(x):=\sup\{\, \rho>0 \,:\,\text{the map $B^G_{\rho}\rightarrow X, g\mapsto gx$ is injective} \,\} $$
to be the injectivity radius of $X$ at the point $x$.
We write $\{ \inj \geq \rho \} = \set{ x \in X : \inj(x) \geq \rho }$, and  $\{ \inj < \rho \}$ for its complement.

\bfparagraph{Grassmannian.}
Given $d\geq 2$, we equip $\R^d$ with its standard Euclidean structure. It extends  naturally to the exterior algebra $\bigwedge^*\R^d$. Namely, if $e_{1}, \dots, e_{d}$ is an orthonormal basis of $\R^d$, then $\{e_{i_{1}}\wedge\dots \wedge e_{i_{k}} \,:\, 1\leq i_{1}<\dots <i_{k}\leq d\}$ is an orthonormal basis of $\Lambda^* \R^d$.
We let $\Gr(\R^d, k)$ denote the collection of $k$-planes in $\R^d$, and set $\Gr(\R^d)=\bigcup_{k=0}^d \Gr(\R^d, k)$ as usual.

Given  $V,W\in \Gr(\R^d)$, we define the \emph{angle functional}
$$\dist_{\angle}(V, W):=\frac{\|\uv\wedge \uw\|}{\| \uv\|\| \uw\|}$$
where $\uv,\uw$ are any non-zero vectors $\uv\in \Lambda^{\dim V}V$, $\uw \in \Lambda^{\dim W}W$.  It makes sense provided  $V,W\neq \{0\}$, and in this case   $\dist_{\angle}(V, W)=0$ amounts to $V\cap W\neq \{0\}$.  We extend the notation to $\Gr(\R^d)^2$ by setting $\dist_{\angle}(V, \{0\})=\dist_{\angle}(\{0\},W)=0$ for any $V,W\in \Gr(\R^d)$. Note  $\dist_{\angle}$ is $\OO(d)$-invariant.

We define the \emph{distance from $V$ to $W$}  by
$$\dist(V \too W):=\sup_{\R v \subseteq V}\inf_{\R w\subseteq W} \dist_{\angle}(\R v, \R w).$$
It is well defined for $V,W\neq \{0\}$. To allow potentially zero subspaces, we use the convention $\dist(V \too \,\{0\})=\1_{V=\{0\}}$, $\dist(\{0\} \too W)=0$ for any $V,W\in \Gr(\R^d)$.
In particular, we find that $\dist(V \too W)=0$ if and only if $V\subseteq W$. We also have the triangle inequality $\dist(V \too W)\leq \dist(V \too S)+\dist(S \too W)$.

We define a \emph{distance on   $\Gr(\R^d)$}  (in the standard sense)  by
\begin{align*}
\dist(V, W) &:= \max\left\{\dist(V \too  W), \,\dist(W \too V) \right\}\\
&= \left\{
    \begin{array}{ll}
        \dist(V\too W) & \mbox{if } \dim V=\dim W \\
        1 & \mbox{else.}
    \end{array}
\right.
\end{align*}

We record $\dist(\cdot, \cdot)$ is  $\OO(d)$-invariant, and equivalent to any distance induced by a Riemannian metric on $\Gr(\R^d)$.
For $r>0$, we let $B_{r}(V)$ denote the open ball in $\Gr(\R^d)$  of center $V$ and radius $r$ for this distance. Note that for $r\leq 1$ we have $B_{r}(V)\subseteq \Gr(\R^d, \dim V)$. It can be checked that for every $V,W\in \Gr(\R^d)$, we have $\dist(V \too W)=\dist(W^\perp \too V^\perp)$, in particular we get:
\begin{equation}\label{grasmannian-distance}
\dist(V, W)=\dist(V^\perp, W^\perp).
\end{equation}

\bfparagraph{Asymptotic notations.}
We use the Landau notation $O(\,\cdot\,)$ and the Vinogradov symbol $\ll$. Given $a,b>0$, we  write $a\simeq b$ for $a\ll b\ll a$. We also say that a statement involving $a,b$  is valid under the condition $a\lll b$ if it holds provided  $a\leq \eps b$ where $\eps>0$ is a small enough constant. When the implicit constants involved in the asymptotic notations $O(\,\cdot\,)$, $\ll$, $\simeq$, $\lll$   depend on some parameters, those are indicated as subscripts.   For instance, $a \lll_p b$ means that the constant $\eps$ above can be taken as a function of the parameter $p$ and nothing else. The absence of subscript indicates absolute constants.

\medskip
\begin{center}
\rule{2cm}{0.4pt}
\end{center}

\bfparagraph{Acknowledgement.}
We are indebted to Yves Benoist for pointing out to us the obstruction concerning $\SL_{3}(\R)$ discussed in \S\ref{Sec-on-the-proof}, and for sharing his insights related to \Cref{avoidW'}.
We are grateful to Nicolas de Saxcé for precious discussions regarding non-concentration assumptions in projection theorems.
We thank Ruotao Yang for showing us an elegant alternative proof of a special case of \Cref{submod-G-F} over the field $\C$.
We thank Yuxiang Jiao and Chengyang Wu for pointing out various typos in an earlier version of this paper.
\medskip

\section{Reduction of the main results and overview}\label{Sec-roadmap}

As we already mentioned, we follow the strategy of \cite{BH24} and both phase I and phase III are already taken care of in that paper, leaving only phase II.
So the main results of the paper all stem from the iteration of a dimensional increment property concerning measures on a homogeneous space transformed under the action of a random walk (\Cref{pr:bootstrap}). This increment property, in turn, is obtained from the conjunction of two phenomena, whose study underpins the entire paper.
 The first concerns dimensional interpolation under the random walk (\Cref{pr:dim-quasi-pres}), the second is about a dimensional increase modulo decomposition (\Cref{pr:splitnu}).
In this section, we present these two key results and derive the main statements from them. We also explain how the remainder of the paper is organized around the proofs of these results.

\subsection{Dimensional interpolation and supercritical decomposition} \label{reductions}

%In this section, we present these two key results and derive the main statements from them. We also outline how the remainder of the paper is organized around the proofs of these results, making this section serve as both a conceptual overview and a roadmap for the paper.

%In this Section, we reduce the proof of the main result, \Cref{main thm}, to two fundamental results. They concern respectively the question of dimensional interpolation and dimensional increase for a measure on a homogeneous space transformed under the action of a random walk. We also explain how the rest of text is structured around the proof of those two results.

\bigskip

%Let $G$ be a non-compact connected real Lie group with finite center and simple Lie algebra $\kg$. Let $\norm{\cdot}$ be a Euclidean norm on $\kg$ and endow $G$ with the associated right-invariant Riemannian metric. Let $\Lambda \subset G$ be a lattice and set $X = G / \Lambda$ equipped with the quotient metric. Below we will denote by $\triangle=(G, \norm{\cdot}, \Lambda)$ these data.

Let $G$ be a non-compact connected real Lie group with finite center and simple Lie algebra. Fix a Euclidean norm $\norm{\cdot}$ on the Lie algebra of $G$. Let $\Lambda \subset G$ be a lattice. Below we will denote by $\triangle=(G, \norm{\cdot}, \Lambda)$ these data. Endow $G$ with the right $G$-invariant Riemannian metric associated to $\norm{\cdot}$, and $X$ with the quotient metric. Recall $m_X$ denotes the Haar probability measure on $X$.
%Writing $d:=\dim G$, we have in particular $m_{X}(B^G_\rho x)\simeq _{\triangle}\rho^d$ for every $\rho>0$, $x\in \{\inj \geq \rho\}$.

Let $\mu$ be a Zariski-dense probability measure on $G$ with finite exponential moment.
%Let $\check\mu$ denote the image measure of $\mu$ under the map $g \mapsto g^{-1}$.
Let $\lambda_1 > \lambda_2 > \dotsb > \lambda_{m+1}$ be the \emph{Lyapunov exponents} of $\Ad(\mu)$, enumerated without repetition and by decreasing order. Let $(j_i)_{1 \leq i \leq m+1}\in (\N^*)^{m+1}$ denote their respective multiplicities. Formally, this means the following.
Consider any choice of maximal compact subgroup $K\subseteq G$ and compatible\footnote{This means $\ka$ is orthogonal to the Lie algebra of $K$ for the Killing form.} Cartan subspace $\ka$ with an open Weyl chamber $\ka^{++}\subseteq \ka$. Let $\kappa_{\mu}\in \ka^{++}$ be the Lyapunov vector of $\mu$ \cite{BQ_book}. Then the pairs $(\lambda_{i}, j_{i})_{i=1, \dots, m+1}$ are given by the eigenvalues and multiplicities of $\ad(\kappa_{\mu})$. By \cite[Theorem 10.9]{BQ_book}, there is also a more concrete characterization: for every $\eps>0$, for large enough $n$, for most $g\sim \mu^n$, the singular values of $\Ad(g)$ are of the form $e^{n\kappa_{1}(g) }\geq \dots \geq e^{n\kappa_{d}(g) }$ where $d=\dim G$, $\kappa_{i}(g)\in \R$, and the vectors $\bigl(\lambda_{1}^{\otimes j_{1}}, \dots, \lambda_{m+1}^{\otimes j_{m+1}}\bigr)$ and $(\kappa_{1}(g), \dots, \kappa_{d}(g))$ are $\eps$-close.

%Given any choice of maximal compact subgroup $K\subseteq G$, the couples $(\lambda_{i}, j_{i})_{i=1, \dots, m+1}$ is characterized by the property that, given any choice of Cartan decomposition\footnote{This means a choice of a maximal compact subgroup $K\subseteq G$, and a Cartan subspace $\ka$ orthogonal to the Lie algebra of $K$ for the Killing form.} on $G$, for every $\eps>0$, for large enough $n$, for most $g\sim \mu^n$, the singular values of $\Ad(g)$ are of the form $e^{n\kappa_{1}(g) }\geq \dots \geq e^{n\kappa_{d}(g) }$ where $d=\dim G$, $\kappa_{i}(g)\in \R$, and the vectors $(\lambda_{1}^{\otimes j_{1}}, \dots, \lambda_{m+1}^{\otimes j_{m+1}})$ and $(\kappa_{1}(g), \dots, \kappa_{d}(g))$ are $\eps$-close. See \cite[Theorem 13.17]{BQ_book} for details.

\begin{definition}
Let $\alpha, \tau\in \R^+$ be parameters.
Let $\nu$ be a Borel measure on $X$, let $\cB$ be a collection of measurable subsets in $X$.
We say $\nu$ is \emph{$(\alpha, \cB, \tau)$-robust}\footnote{
It should be noted that this definition deviates from our previous work~\cite{BH24} where $\nu_1$ is additionally required to be supported away from the cusps.
The definition here is adapted to the argument employed in the present paper.}, if we can decompose $\nu$ as a sum of two Borel measures $\nu = \nu_1 + \nu_2$ such that
\begin{enumerate}
\item $\forall B \in \cB$, $\nu_1(B) \leq m_X(B)^\alpha$
\item $\nu_2(X) \leq \tau$.
\end{enumerate}
\end{definition}
In practice, $\alpha$ and $\tau$ will be smaller than $1$, and $\cB$ will be a collection of balls.
For $\rho > 0$, we let $\cB_\rho$ denote the collection of all balls of radius $\rho$ in $X$.
For $I \subset \R$, we set $\cB_I = \bigcup_{\rho \in I} \cB_\rho$.

%\begin{proposition}\label{pr:bootstrap}
%Given $\varkappa > 0$, there is $\eps > 0$ such that the following holds for all $\alpha \in [\varkappa, 1 - \varkappa]$, and all $\delta > 0$ sufficiently small.
%Let $\nu$ be a $(\alpha, \cB_{[\delta, \delta^\eps]}, \tau)$-robust measure with
%$\nu(\{\inj < \delta^\eps\}) \leq \tau$.
%Then for some $n \asymp \frac{1}{\lambda_1} \abs{\log \delta}$,\marginpar{Notation $\asymp$?}
%$\mu^{*n} * \nu$ is $(\alpha + \eps, \cB_{\delta^{1/2}}, 2 \tau + \delta^\eps)$-robust and satisfies $\mu^{*n} * \nu (\{ \inj < \delta^{1/2} \}) \leq 2 \tau +\delta^\eps$.
%\end{proposition}

\bigskip

Our first key result is the following dimensional interpolation property concerning the action of the $\mu$-walk on a given initial distribution on $X$.

\begin{proposition}[Dimensional interpolation]
\label[proposition]{pr:dim-quasi-pres} 
Let $X$, $\mu$, $(\lambda_i)$, $(j_i)$ be as above.
Let $s \in (0, \frac{1}{4\lambda_1}]$ and $\eps_0, \eps, \delta > 0$.

Let $\nu$ be a Borel measure on $X$ of mass at most $1$, supported on $\{ \inj \geq \delta^{2/3} \}$, and which is $(\alpha_i, \cB_{\delta^{1 - s \lambda_i}}, 0)$-robust for some parameter $\alpha_i>0$, for all $1 \leq i \leq m+1$.
Let $\beta \in \R$ be such that
\[
(\dim G) \beta = \sum_{i = 1}^{m+1} (1 - s \lambda_i) j_i \alpha_i.
\]

If $\eps, \delta \lll_{\triangle, \mu, s, \eps_{0}}1$, then setting $n = \lfloor s \abs{\log \delta} \rfloor$, we have that
\[
\text{$\mu^{n} * \nu$ is $(\beta - \eps_0, \cB_\delta, \delta^\eps)$-robust.}
\]
\end{proposition}

\begin{remark}
In \Cref{pr:dim-quasi-pres}, if all the $\alpha_{i}$'s are equal to some $\alpha$, then $\beta=\alpha$ as well. This is because $\sum_{i = 1}^{m+1} (1 - s \lambda_i) j_i=\dim G- s\sum_{i = 1}^{m+1} \lambda_i j_i=\dim G$ where the last inequality uses that $\Ad(G)\subseteq \SL(\kg)$ (so the sum of Lyapunov exponents, counted with their multiplicity,  is zero).
%the observation that the Lyapunov spectrum of $\Ad(\mu)$ is symmetric with opposite eigenvalues $\lambda_{i}, -\lambda_{i}$ having the same multiplicity.
Therefore, \Cref{pr:dim-quasi-pres} expresses in particular that \emph{the random walk does not decrease the dimension of a prescribed initial distribution}.
We will also use it in a context where the $\alpha_{i}$'s are not all equal, via \Cref{pr:preservenu} below.
\end{remark}

\begin{proof}
\Cref{pr:dim-quasi-pres} is a direct consequence of \Cref{RW-preserve-dim}, whose proof will be established in \Cref{Sec-dim-pres}, relying on Sections \ref{Sec-multislicing-thms}, \ref{Sec-optimal-subcritical}, \ref{Sec-submod}.
\end{proof}

Applying \Cref{pr:dim-quasi-pres} at scale $\delta^{1/2}$ and with $s = \frac{1}{8\lambda_1}$, we deduce easily
\begin{corollary}
\label[corollary]{pr:preservenu}
Let $X$, $\mu$ and $(\lambda_i)$ be as above. Let $\eps_{0}, \eps, \delta>0$.

Let $\nu$ be a Borel measure  on $X$ of mass at most $1$, supported on $\{ \inj \geq \delta^{1/3} \}$, and which is $(\alpha, \cB_{\delta^{1/2-\lambda_{i}/(16\lambda_{1})}}, 0)$-robust for some $\alpha\in \R^+$ and every $i=1, \dots, m+1$. Assume also that $\nu$ is
$$\text{either $(\alpha + \eps_0, \cB_{\delta^{1/2}}, 0)$-robust\,\,\,\, \,\,or\,\,\,\,\,\, $(\alpha + \eps_0, \cB_{\delta^{7/16}}, 0)$-robust}.$$

If $\eps, \delta \lll_{\triangle, \mu, \eps_{0}}1$, then for $n = \lfloor \frac{1}{16\lambda_1}\abs{\log \delta} \rfloor$ and $d=\dim G$,
$$ \text{$\mu^{n}*\nu$ is $\left(\alpha + {\textstyle\frac{1}{4d}}\eps_{0}, \cB_{\delta^{1/2}}, \delta^\eps\right)$-robust.}$$
\end{corollary}

Our second key result claims that the $\mu$-walk on $X$ in fact improves the dimensional properties of a given initial distribution, but for that we need to partition the new distribution into two submeasures, and look at different scales for each piece.

\begin{proposition}[Supercritical decomposition]
\label[proposition]{pr:splitnu}
Let $X$, $\mu$ and $\lambda_1, \lambda_2$ be as above.
Let $\varkappa, \eps, \delta > 0$ and $\alpha \in [\varkappa, 1 - \varkappa]$.
Let $\nu$ be a Borel measure  on $X$, supported on $\{ \inj \geq \delta^{\eps} \}$, and which is $(\alpha, \cB_{[\delta, \delta^{\eps}]}, 0)$-robust.
Set $n = \bigl\lfloor \frac{1}{16(\lambda_1 + \lambda_2)} \abs{\log \delta} \bigr\rfloor$.

If $\eps, \delta\lll_{\triangle, \mu, \varkappa}1$, then $\mu^{n}*\nu$ is the sum of a $(\alpha + \eps, \cB_{\delta^{1/2}}, \delta^\eps)$-robust measure and a $(\alpha + \eps, \cB_{\delta^{7/16}}, \delta^\eps)$-robust measure.
\end{proposition}

We note that $0$ must be among the Lyapunov exponents of $\Ad(\mu)$, hence $\lambda_{1}>\lambda_{2}\geq 0$, so the above denominator $\lambda_{1}+\lambda_{2}$ is indeed positive.

\begin{proof} This is  a direct consequence of \Cref{sup-mult-X} (applied with $t_{1}=1/2$ and $t_{2}=7/16$, and noting the assumptions of support and non-concentration  on $\nu$ imply $\nu(X)\ll_{\triangle}\delta^{-\eps \dim X}$).  The proof of \Cref{sup-mult-X} will be carried out in \Cref{Sec-RWincrease}, relying on Sections \ref{Sec-multislicing-thms}, \ref{Sec-optimal-subcritical}, \ref{Sec-submod}, \ref{Sec-dim-pres}.
%Let us justify that we may indeed suppose $\nu(X)\leq 1$. Note first we can further assume $\delta$ small enough  in terms of $\eps$ as well. Indeed, if the conclusion holds for a pair $(\eps, \delta)$ then it holds for all $(\eps', \delta)$ with $\eps'\in (0, \eps)$.
 %Now, the robustness assumption applied at scale $\delta^{\eps}$ and the hypothesis on the support of $\nu$ together imply $\nu(X)\ll_{\triangle} \delta^{- \eps \dim X }$. Up to taking $\eps\lll_{\triangle, \mu, \varkappa}1$ and renormalizing, \Cref{pr:splitnu} thus reduces to the case $\nu(X)\leq 1$. 
\end{proof}

\bigskip
Admitting \Cref{pr:dim-quasi-pres} and \Cref{pr:splitnu} for now, we conclude the proof of the main results. We need the next quantitative recurrence estimate, which follows from \cite[Lemma 4.8]{BH24}.
\begin{lemma}[Non-escape of mass \cite{BH24}]
\label[lemma]{lm:reccurence} 
There exists a constant $c =c(\triangle, \mu)> 0$ such that for every Borel measure $\nu$ on $X$ of mass at most $1$, every $n \geq 0$, and every $\rho, r \in (0, 1)$, we have
\[
(\mu^{ n} * \nu)(\{ \inj < r \}) \ll_{\triangle, \mu} r^{c} (e^{- c n} \rho^{-1} + 1) + \nu( \{ \inj < \rho \}).
\]
\end{lemma}

Combining all the previous results, we deduce the announced dimensional increment property. For the purpose of iteration, it also comes with a control of the injectivity radius.

\begin{proposition}[Dimensional increment]
\label[proposition]{pr:bootstrap}
Let $X$, $\mu$, $\lambda_{1}$, be as above.
Let $\varkappa, \eps, \delta \in (0, 1)$, $\tau\in \R^+$, and $\alpha \in [\varkappa, 1 - \varkappa]$.
Let $\nu$ be a $(\alpha, \cB_{[\delta, \delta^\eps]}, \tau)$-robust measure on $X$ satisfying  and
$\nu(\{\inj < \delta^\eps\}) \leq \tau$.

If $\eps, \delta\lll_{\triangle,\mu, \varkappa}1$,
then for some $n \simeq \frac{1}{\lambda_1} \abs{\log \delta}$, the measure
$\mu^{n} * \nu$ is $(\alpha + \eps, \cB_{\delta^{1/2}}, 2 \tau + \delta^\eps)$-robust and satisfies $(\mu^{n} * \nu) (\{ \inj < \delta^{1/2} \}) \leq 2 \tau +\delta^\eps$.
\end{proposition}

\begin{proof}[Proof of \Cref{pr:bootstrap}] 
By definition of robustness, we can write $\nu = \nu_0 + (\nu - \nu_0)$ where $\nu_0$ is a $(\alpha, \cB_{[\delta, \delta^\eps]}, 0)$-robust measure supported on $\{\inj \geq \delta^\eps\}$ and $\nu - \nu_0$ is a  Borel measure of total mass at most $ 2\tau$. It is enough to show the lemma for $\nu_{0}$, in other terms we may assume $\tau=0$.
Noting $\nu(X)\ll_{\triangle} \delta^{-\eps \dim X}$ and renormalizing if necessary, we may assume $\nu$ has mass at most $1$. 
Moreover, throughout the proof, we may assume $\delta$ small enough depending on $\eps$ as well. Indeed, if the conclusion holds for a pair $(\eps, \delta)$ then it holds for all $(\eps', \delta)$ with $\eps'\in (0, \eps)$.

Provided $\eps, \delta\lll_{\triangle, \mu, \varkappa}1$, 
we may apply \Cref{pr:splitnu} to $\nu$. Writing $n_1 = \lfloor \frac{1}{16(\lambda_1 + \lambda_2)} \abs{\log \delta} \rfloor$, we obtain a constant $\eps_0=\eps_{0}(\triangle, \mu, \varkappa) \in (0, \varkappa)$ and a decomposition
\[
\mu^{ n_1} * \nu = \nu_1 + \nu_2
\]
where $\nu_1$ is a $(\alpha + \eps_0, \cB_{\delta^{1/2}}, \delta^{\eps_0})$-robust measure, while $\nu_2$ is a $(\alpha + \eps_0, \cB_{\delta^{7/16}}, \delta^{\eps_0})$-robust measure. This is not enough just yet, because the scales $\delta^{1/2}, \delta^{7/16}$ where the gain $\eps_{0}$ occurs are different. For the rest of the proof, we aim to apply more convolutions by $\mu$ in order to reconcile the scales (via  \Cref{pr:preservenu}). 

Note the measure $\mu^{ n_1} * \nu$  enjoys robustness properties at other scales. Indeed \Cref{pr:dim-quasi-pres} (and its remark) apply to $\nu$ at any scale in the range $\bigl[\delta^{9/16}, \delta^{7/16}\bigr]$, with dimensional loss $\eps_{0}/(8d)$. 
More precisely, \Cref{pr:dim-quasi-pres}  (applied several times) yields some constant some $\eps_{1}=\eps_{1}(\triangle, \mu, \varkappa)>0$ such that for any finite subset $I \subset [\delta^{9/16}, \delta^{7/16}]$,  and provided $ \delta\lll_{\triangle, \mu, \varkappa, |I|}1$,   the measure
$\mu^{ n_1} * \nu$ is $(\alpha - \frac{1}{8d}\eps_0, \cB_{I}, \delta^{\eps_1})$-robust where $d=\dim G$.
To prepare for the use of \Cref{pr:preservenu}, we choose $I=\setbig{\delta^{1/2-\lambda_{i}/(16\lambda_{1})} \,:\, i=1, \dotsc, m+1}$. We also note that the robustness of $\mu^{ n_1} * \nu$ automatically transfers to $\nu_1$, $\nu_2$.

Observe also that  the measure  $\mu^{ n_1} * \nu$ is not too concentrated near the cusps.
Indeed, using \Cref{lm:reccurence} with $r = \delta^{1/3}$ and $\rho = \delta^\eps$, we have
\[
(\mu^{ n_1} * \nu) \{ \inj < \delta^{1/3} \} \ll_{\triangle, \mu} \delta^{\frac{c}{3}}\left(\delta^{\frac{ c}{16(\lambda_1 + \lambda_2)} - \eps} + 1\right).
\]
Hence $(\mu^{ n_1} * \nu) \{ \inj < \delta^{1/3} \} \leq \delta^{c/4}$ as soon as $\eps < \frac{ c}{16(\lambda_1 + \lambda_2)}$ and $\delta\lll_{\triangle, \mu}1$. This automatically transfers to $\nu_1$, $\nu_2$

Combining the three previous paragraphs,  we can write $\nu_1 = \nu_3 + (\nu_1 - \nu_3)$ and $\nu_2 = \nu_4 - (\nu_2 - \nu_4)$ where $\nu_3$, $\nu_{4}$ are Borel measures that are supported on $\{ \inj \geq \delta^{1/3} \}$, as well as
$(\alpha - \frac{1}{8d}\eps_0, \cB_{I}, 0)$-robust, and respectively $(\alpha + \eps_0, \cB_{\delta^{1/2}}, 0)$-robust, $(\alpha + \eps_0, \cB_{\delta^{7/16}}, 0)$-robust; while $\nu_1 - \nu_3$ and $\nu_2 - \nu_4$ are Borel measures of total mass at most $\delta^{\eps_0} +\delta^{\eps_1} + \delta^{c/4}$.

We are now in a position to apply \Cref{pr:preservenu} to the measures $\nu_3$, $\nu_4$, and with $\alpha - \frac{1}{8d}\eps_0$ in the place of $\alpha$.
Write $n_2 = \lfloor \frac{1}{16\lambda_1}\abs{\log \delta} \rfloor$,
We obtain that for $\eps_{2}, \delta\lll_{\triangle, \mu, \varkappa}1$, we have
$\mu^{n_2}*\nu_3$ and $\mu^{n_2}*\nu_4$ both $(\alpha + \frac{1}{8d}\eps_{0}, \cB_{\delta^{1/2}}, \delta^{\eps_{2}})$-robust.

Set $n = n_1 + n_2$, from the above, $\mu^{ n} * \nu - \mu^{ n_{2}} *(\nu_3 + \nu_4)$ has total mass at most $ 2(\delta^{\eps_0} +\delta^{\eps_1} + \delta^{c/4})$.
It follows that $\mu^{ n} * \nu$ is $(\alpha + \frac{1}{8d}\eps_{0}, \cB_{\delta^{1/2}}, \tau')$-robust, where $\tau' := 2(\delta^{\eps_0} +\delta^{\eps_1} + \delta^{c/4} +\delta^{\eps_{2}})$.
Provided $\eps<\min(\eps_{0}, \eps_{1}, \eps_{2}, c/4)$, and $\delta \lll_{\eps}1$, we have $\tau' <  \delta^{\eps}$.

Finally, by \Cref{lm:reccurence} applied with $\rho = \delta^\eps$ and $r = \delta^{1/2}$,
we have
\[
(\mu^{n} * \nu) \{ \inj < \delta^{1/2} \} \ll_{\triangle, \mu} \delta^{c/2} \left(\delta^{(\frac{1}{16(\lambda_1 + \lambda_2)} + \frac{1}{16 \lambda_1}) c - \eps} + 1\right),
\]
leading to the desired bound on $(\mu^{n} * \nu) \{ \inj < \delta^{1/2} \}$ provided $\delta \lll_{\triangle, \mu, \eps}1$.
\end{proof}

We can now derive from \Cref{pr:bootstrap} the main results announced in \Cref{Sec-intro}.

\begin{proof}[Proof of \Cref{Eff-eq-thmI}]
%For $G$ isogenous to $\SO(2,1)$ or $\SO(3,1)$, \Cref{Eff-eq-thmI} follows from our first paper \cite[Proof of Theorem 4.1]{BH24}. This restriction on $G$ is only used to obtain the dimensional increment property \cite[Proposition 4.9]{BH24}, generalized here to arbitrary non-compact simple Lie groups with finite center via \Cref{pr:bootstrap}. 
Once we know there is a dimensional increment, effective equidistribution can be deduced verbatim from \cite{BH24}.
Namely, arguing as in \cite[Section 4.3.4]{BH24}, we may apply \Cref{pr:bootstrap} iteratively in order to bootstrap the dimension of $\nu$ arbitrarily close to the ambient dimension, $\dim X$. The argument can be performed exactly as in \cite{BH24}, using our \Cref{pr:bootstrap} instead of \cite[Proposition 4.9]{BH24}, and noting the notion of robustness used in \cite{BH24} already takes into account the injectivity radius. Then, we go from high dimension to equidistribution using \cite[Proposition 4.14]{BH24}, concluding the proof of \Cref{Eff-eq-thmI}.
\end{proof}

\begin{proof}[Proof of \Cref{Eff-eq-thmII}]
Invoking the extra arithmeticity assumptions, \cite[Theorem 3.3]{BH24} guarantees that $\mu^n*\delta_{x}$ acquires positive dimension above scale $R^{-1}$ for $n \geq A\log R + A \max \{|\log \dist(x,W_{\!\mu, R^A})|,\, \dist(x,x_{0})\}$. We then apply \Cref{Eff-eq-thmI} and \Cref{lm:reccurence} to conclude. See \cite[Section 5, Proof of Theorem 1.3]{BH24} for details.
\end{proof}

\begin{proof}[Proof of \Cref{thm:equidis}] It is identical to the proof of \Cref{Eff-eq-thmII}, but using \cite[Proposition 5.1]{BH24} instead of \cite[Theorem 3.3]{BH24}. This frees us from artihmeticity assumptions, but we loose the rate of equidistribution.
\end{proof}

\begin{proof}[Proof of \Cref{cr:dioph} and \Cref{cor-orbites}]
It is identical to that of \cite[Corollaries 1.4, 1.5]{BH24}, using \Cref{Eff-eq-thmII} instead of \cite[Theorem 3.3]{BH24}.
\end{proof}

 As we have just seen, all our main results reduce to \Cref{pr:dim-quasi-pres} and \Cref{pr:splitnu}. The remainder of the paper is dedicated to the proof of these two propositions.

\bigskip
\subsection{Overview of the paper} The rest of the article is organized as follows.
\bigskip

\noindent In \Cref{Sec-multislicing-thms}, we present \emph{ \bf multislicing estimates}. We consider a random box in $\R^d$ with side lengths of the form $(\delta^{r_{1}}, \dots, \delta^{r_{d}})$ for some small $\delta>0$, and parameters $r_{i}\in (0, 1)$ not all equal. It determines a random partial flag. Given a measure $\nu$ on $\R^d$ with dimension $\alpha$ above scale $\delta$, we establish an upper bound on the mass granted by $\nu$ to these random boxes. More precisely, we assume that subspaces from the random partial flag satisfy a subcritical projection theorem, and derive that $\nu$ has dimension at least $\alpha-\eps$ with respect to all translates of a typical random box (subcritical estimate). Moreover, under extra supercritical assumptions, we also prove a supercritical estimate, i.e., with $\eps$ gain instead of $\eps$ loss. It takes the form of a supercritical decomposition, better suited for our application to random walks. The proofs are similar to \cite[Section 2]{BH24} and postponed to \Cref{ss:Appendix-Multislicing}.

%\noindent In \Cref{Sec-multislicing-thms}, we present \emph{ \bf multislicing estimates}. We consider a random box in $\R^d$ with side lengths of the form $(\delta^{r_{1}}, \dots, \delta^{r_{d}})$ for some small $\delta>0$, and parameters $r_{i}\in (0, 1)$ not all equal. We assume some non-concentration for the associated random partial flag. Given a measure $\nu$ on $\R^d$ with dimension $\alpha$ above scale $\delta$, we establish an upper bound on the mass granted by $\nu$ to all translates of a typical random box. Under weak non-concentration assumptions, this upper bound conveys that $\nu$ has dimension at least $\alpha-\eps$ from the point of view of those translates (subcritical estimate). Under stronger non-concentration assumptions, we also prove a supercritical estimate, i.e. with $\eps$ gain instead of $\eps$ loss. It takes the form of a supercritical decomposition, better suited for our application to random walks.

\bigskip

\noindent In \Cref{Sec-optimal-subcritical}, we establish a \emph{\bf subcritical projection theorem under optimal non-concentration assumptions}. We consider a random orthogonal projector on $\R^d$ whose kernel has a high probability not to intersect too much any given proper subspace of $\R^d$. We conclude that for every set $A\subseteq \R^d$ with (discretized) dimension at least $\alpha$, for an event with high probability, the image of $A$ under the projector has (discretized) dimension at least $\alpha/d-\eps$ where $\eps$ is arbitrarily small. The proof makes use of a quantitative upper bound for Brascamp-Lieb constants. The latter is deferred to a separate paper \cite{BH25-BrascampLieb} and builds upon the work of Bennett-Carbery-Christ-Tao \cite{BCCT08}.

\bigskip

\noindent In \Cref{Sec-submod}, we establish a \emph{ \bf submodular inequality for irreducible representations over arbitrary fields}. This inequality is of significance on its own and the section can be read independently of the rest of the paper. It will be applied later in the context of random walks in order to check the non-concentration assumptions relevant to subcritical estimates. The proof of the submodular inequality makes use of supermodular functions on the Grassmannian and resonates with the theory of Harder-Narasimhan filtrations.

%The proof of the submodular inequality relies on a case by case approach. It uses the classification of simple complex Lie algebras, and combinatorial arguments to exhibit common convexity properties.

\bigskip
\noindent In \Cref{Sec-dim-pres}, we prove \Cref{pr:dim-quasi-pres}, i.e the \emph{\bf dimensional interpolation property under the action of random walks}. The proof combines Sections \ref{Sec-multislicing-thms}, \ref{Sec-optimal-subcritical}, \ref{Sec-submod}. We also put forward a linearization technique which allows for linearization at microscopic scales. This technique is inspired by Shmerkin \cite{Shmerkin}. It improves upon the linearization procedure used in \cite{BH24}, which was taking place at macroscopic scales, and failed for higher rank groups such as $\SL_{d}(\R)$ where $d\geq3$ (see the remark following the proof of Lemma 4.10 in \cite{BH24}).

\bigskip
\noindent In \Cref{Sec-RWincrease}, we prove \Cref{pr:splitnu}, i.e. the \emph{\bf supercritical decomposition under the action of a random walk}. The proof makes use of Sections \ref{Sec-multislicing-thms}, \ref{Sec-optimal-subcritical}, \ref{Sec-submod}, \ref{Sec-dim-pres}. It boils down to a supercritical alternative property regarding projections onto maximally expanded and maximally contracted directions for $\Ad(g)$ where $g\sim \mu^n$.
 As discussed in \S\ref{Sec-on-the-proof}, we may only rely on a weak form of non-concentration for those subspaces. It is incarnated by \Cref{nc-highest-weight}. Note that if the adjoint representation of $G$ is proximal, then the section can be simplified a lot: there is no need to discuss a supercritical alternative because the maximally contracted direction of $\Ad(g)$ is known to satisfy a supercritical projection theorem. The point of the section is to deal with simple Lie groups which are \emph{not} $\Ad$-proximal, and for which the maximally contracted direction of $(\Ad(g))_{g\sim \mu^n}$ fails to satisfy the transversality property required in projection theorems \`a la Bourgain (e.g. $G=\SO(7,1)$).

\bigskip
\noindent In \Cref{ss:Appendix-Multislicing}, we detail the proof of the multislicing estimates from \Cref{Sec-multislicing-thms}. In \Cref{lack}, we highlight a drastic form of non-transversality for highest weight subspaces of $\SO(n,1)\acts \so(n,1)$.

\section{Multislicing machinery} \label{Sec-multislicing-thms}

In this section, we explain how a collection of projection theorems can be combined into a multislicing theorem. More precisely, we consider a probability measure $\nu$ on the unit cube of a Euclidean space and we suppose $\nu$ satisfies  certain dimensional estimates with respect to balls. We partition the unit cube into smaller cubes, and cover each one of them with translates of an asymmetric box, which is chosen randomly according to some measure. For each small cube, the associated random box determines  a random partial flag of $\R^d$. Assuming each random subspace involved in the flag satisfies a subcritical projection theorem,  we show the dimension estimates for $\nu$ with respect to such boxes are almost as good as those assumed for balls (subcritical regime). If moreover, one random subspace enjoys a supercritical projection theorem, we  obtain a dimensional gain when estimating the $\nu$-mass of the boxes (supercritical regime). More generally, under a weaker condition which we call the supercritical alternative property, we prove that $\nu$ can be partitioned into two Borel submeasures that each enjoy dimensional gain, although at different scales. This extension will be crucial for our application to random walks.

%In this section, we present various multislicing theorems. We consider a probability measure $\nu$ on the unit cube of a Euclidean space and assume  certain dimensional properties with respect to balls. We partition the unit cube into smaller cubes, and cover each one of them with translates of an asymetric box whose axes are randomized so as to satisfy certain projection theorems. Under weak assumptions on non-concentration, we show the dimension estimates for $\nu$ with respect to such boxes are almost as good as those assumed for balls (subcritical case). We also explore under which circumstances we may obtain a dimensional gain (supercritical case).

\bigskip
We place ourselves in $\R^d$ where $d\geq2$, endowed with its standard Euclidean structure.

\bfparagraph{Pixelization.}
Given $\eta>0$, we write $\cD_{\eta}$ the partition of $\R^d$ generated by the cell
$$ Q_{\eta}:=[0,2^{-k}[^d$$
where $2^{-k}$ is the dyadic upper-approximation of $\eta$, i.e., $k\in \Z$ and satisfies
$2^{-k-1}< \eta \leq  2^{-k}$.

\bfparagraph{Boxes.}
Let $m\in \{0, \dots, d-1\}$. We set
$$\cP_{m}(d):=\{(j_{i})_{i=1}^{m+1}\in \N^{m+1}_{\geq 1} :\,  d=j_{1}+\dots+j_{m+1}\},$$
$$\increasing_{m}:=\{(r_{i})_{i=1}^{m+1} :\,  0\leq r_{1} < \dots  < r_{m+1}\leq 1\}.$$
Every $\uj\in \cP_{m}(d)$ determines a collection of partial flags $\cF_{\uj}$, consisting of all the tuples  $(V_{i})_{i=1}^{m+1} \in \Gr(\R^d)^{m+1}$ such that
\[
\{0\}\subsetneq  V_{1}\subsetneq \dots  \subsetneq V_{m+1}=\R^d \qquad\text{with}\quad\dim V_{i}=j_{1}+\dots +j_{i},\, \forall i.
\]
%$$\cF_{\uj}:=\{(V_{i})_{i=1}^{m+1} \in \Gr(\R^d)^{m+1}\,:\, \{0\}\subsetneq  V_{1}\subsetneq \dots  \subsetneq V_{m+1}=\R^d\},$$
For $\uV=(V_{i})_{i=1}^{m+1} \in \cF_{\uj}$, $\ur=(r_{i})_{i=1}^{m+1}\in \increasing_{m}$, and $\delta\in (0, 1)$,  we introduce the box
$$B_{\delta^\ur}^{\uV} :=\sum_{i=1}^{m+1}B^{V_{i}}_{\delta^{r_{i}}}.$$
We call $\uV$ the partial flag (or the filtration) carrying the box $B_{\delta^\ur}^{\uV}$.

\bfparagraph{Dimension.}
For heuristics, it will be convenient to talk about the dimension of a measure with respect to certain shapes in $\R^d$. We say a measure $\nu$ on $\R^d$ has \emph{normalized dimension} at least $\alpha \in [0, 1]$ with respect to a collection $\cS$ of measurable subsets of $\R^d$ if every  $S\in\cS$ satisfies $\nu(S)\leq (\leb S)^\alpha$. When $\cS$ is the collection of balls of  given radius $r>0$, we just talk about normalized dimension at scale $r$.

\bigskip

In order to state our subcritical multislicing theorem, we formalize what it means for a measure on the Grassmannian to satisfy a subcritical projection theorem. Given $A\subseteq \R^d$, $\delta>0$, we write $\cN_{\delta}(A)$ the smallest number of  $\delta$-balls needed to cover $A$.

\begin{definition}
\label[definition]{Sub-crit-P} 
 Let $\sigma$ be a probability measure on $\Gr(\R^d)$, let $ \delta, \eps, \tau>0$.  We say $\sigma$ has the \emph{subcritical property} $\SubP$ with  parameters $(\delta, \eps,\tau)$ if  for every set $A\subseteq B^{\R^d}_{1}$, the exceptional set
\begin{equation*}
\begin{split}
\cE := \Bigl\{\, V   : \exists A' \subseteq A \,\,&\text{ with }\,\,\cN_{\delta}(A')\geq \delta^{\eps} \cN_{\delta}(A)\\
& \text{ and }\, \cN_{\delta}(\pi_{|| V}A')< \delta^\tau \cN_{\delta}(A)^{\frac{\dim V^\perp}{d}}\Bigr\}
\end{split}
\end{equation*}
has measure $\sigma(\cE)\leq \delta^\eps$.
\end{definition}

We now present our subcritical multislicing theorem. The unit ball of $\R^d$ is subdivided into cubes $Q\in \cD_{\eta}$ for some fixed $\eta>0$. Within each $Q$, we consider a  box $B^{\uV_{Q, \theta}}_{\delta^\ur}$ where $\uV_{Q,\theta}=(V_{Q, \theta, i})_{i=1}^{m+1}$ is a partial flag, randomized through a common parameter $\theta$. We assume  that for each $i=1, \dots, m$, the random subspace $V_{Q,\theta, i}$ satisfies a subcritical projection theorem at a scale $\delta^{r_{i+1}}$, uniformly in $Q\in \cD_{\eta}$.
The main output is that if  a measure $\nu$ has normalized dimension \emph{at least} $\alpha$ at scales $(\delta^{r_{k}})_{k=1, \dots, m+1}$, then $\nu$  must have normalized dimension \emph{almost} $\alpha$ with respect to translates  $(B^{\uV_{Q,\theta}}_{\delta^\ur}+v)_{v\in \R^d}$ in each block $Q$, up to choosing $\theta$ outside of an event of small $\sigma$-mass and putting aside a small part of the measure $\nu$ (that may depend on $\theta$).

\begin{thm}[Subcritical multislicing]
\label{subc-mult}
Let $d>m\geq 1$, $\uj\in \cP_{m}(d)$, $\ur \in \increasing_{m}$,  $\delta\in (0, 1)$.  Let $\eta\in [\delta^{r_{1}}, 1]$  and $\tau, \eps, \eps'>0$.

Let $(\Theta, \sigma)$ be a probability space. For each $Q\in \cD_{\eta}$, consider a measurable map $\Theta \rightarrow  \cF_{\uj}, \theta \mapsto \uV_{Q,\theta}=(V_{Q,\theta, i})_{i}$.  Assume that  for every  $i \in \{1, \dots, m\}$, the distribution of $(V_{Q, \theta, i})_{\theta\sim \sigma}$ satisfies $\SubP$ with parameters $(\delta^{r_{i+1}}, \eps, \tau)$.

Let $\nu$ be a Borel measure on $B^{\R^d}_{1}$ of mass at most $1$, and for $i=1, \dots, m+1$, let $t_{i}>0$ such that for all $v\in \R^d$,
$$\nu(B^{\R^d}_{\delta^{r_{i}}}+v)\leq t_{i}.$$

If $\eps'\lll\eps$ and $\delta^{r_{2}} \lll_{d, \eps} 1$, then there exists $\cE \subseteq \Theta$ such that $\sigma(\cE)\leq \delta^{r_{2}\eps'}$ and for all $\theta\in\Theta \setminus \cE$, there is a set  $F_{\theta} \subseteq \R^d$ with  $\nu(F_{\theta}) \leq \delta^{r_{2}\eps'}$ and such that for all $Q\in \cD_{\eta}$, $v\in \R^d$,
\[ \nu_{|Q\smallsetminus F_\theta}\left(B^{{ \uV_{Q,\theta}}}_{\delta^\ur} + v \right) \leq  \delta^{-(\tau+\eps) \sum_{i=2}^{m+1}r_{i}}\prod_{i}t_{i}^{j_{i}/d}.\]
\end{thm}

\begin{remark}
The implicit constant in the upper bound $\delta^{r_{2}}\lll_{d, \eps} 1$ only depends on $d$ and  a positive lower bound on $\eps$.
\end{remark}

\bigskip
The term $\delta^{-(\tau+\eps) \sum_{i=2}^{m+1}r_{i}}$ in the conclusion represents a dimensional \emph{loss}.
We now explain that we obtain a dimensional \emph{gain} under the extra assumption that for at least one $i$, the distributions of $(V_{Q, \theta, i})_{\theta\sim \sigma}$ where $Q\in \cD_{\eta}$ satisfy a a supercritical projection theorem. Motivated by our application to random walks on simple homogeneous spaces, we present in fact a more general statement, which only assumes a supercritical alternative. In order to present this notion,  given $\alpha, \tau, \delta>0$, $A\subseteq \R^d$, we set
\begin{equation} \label{notation-cEaed}
\begin{split}
\cE^{(\alpha, \tau)}_{\delta}(A) \defeq \bigl\{\, V\in \Gr(\R^d) : \exists A' \subseteq A \,\,&\text{ with }\,\,\cN_{\delta}(A')\geq \delta^\tau \cN_{\delta}(A)\\
 \text{and }& \cN_{\delta}(\pi_{||V}A') < \delta^{- \alpha \dim V^\perp  -\tau} \,\bigr\}.
\end{split}
\end{equation}

\begin{definition}
\label[definition]{SAP}
 Let $\sigma_{1},\sigma_{2}$ be  probability measures on $\Gr(\R^d)$, let $\delta, \varkappa, \tau>0$.  We say $(\sigma_{1},\sigma_{2})$ has the \emph{supercritical alternative property} $\SAP$ with parameters $(\delta, \varkappa,\tau)$ if the following holds.

Let $A\subseteq B^{\R^d}_{1}$ be any non-empty subset satisfying for some
 $\alpha\in[\varkappa, 1-\varkappa]$, for  $\rho \geq \delta$,
\begin{equation}\label{nc-dim-dalpha}
\sup_{v \in \R^d} \cN_\delta\left(A \cap B^{\R^d}_\rho(v)\right) \leq \delta^{-\tau} \rho^{d \alpha } \cN_\delta(A).
\end{equation}
Then there  exists $A'\subseteq A$ such that
$$ \min_{p=1,2}\sigma_{p}\left(\cE^{(\alpha, \tau)}_{\delta}(A')\right) \leq \delta^\tau.$$
\end{definition}

Roughly speaking, the above property considers an arbitrary $\delta$-separated set $A$ on which the uniform probability measure has  normalized dimension almost  $\alpha$ at scales above $\delta$. It requires the existence of a subset $A'$ of $A$ and $p\in \{1, 2\}$ such that for most projections $(\pi_{||V})_{V\sim \sigma_{p}}$, all rather large subsets of $A'$ have a big image under $\pi_{||V}$ (say normalized box dimension at least $\alpha+\tau/d$). The  term  $\varkappa$ constrains $\alpha$ to be away from $0$ and $1$,  while  $\tau$ controls the dimensional increase, the size of $A'$, and tempers the non-concentration of $A$.

\bigskip
With this notion at hand, we can formulate a supercritical multislicing decomposition theorem for measures. We keep the partition of $\R^d$ into $\cD_{\eta}$-cubes  $Q$ for some fixed $\eta>0$. We consider two types of boxes, whose geometries are locally given by partial flags $\cV_{Q,\theta}, \cW_{Q,\theta}$ randomized through $\theta\sim \sigma$, and  fixed exponents $\ur, \us$. We keep the non-concentration assumption from \Cref{subc-mult}. We consider  exponents $\ur, \us$ that coincide on a pair of consecutive entries, say $r_{i_{1}}=s_{i_{2}}$ and $r_{i_{1}+1}=s_{i_{2}+1}$,  and we assume that  the corresponding random projectors $(\pi_{||V_{Q,\theta, i_{1}}})_{\theta\sim \sigma}$ and $(\pi_{|| W_{Q,\theta, i_{2}}})_{\theta\sim \sigma}$ satisfy the aforementioned supercritical alternative at an appropriate scale. We conclude that any measure $\nu$ with normalized dimension at least $\alpha$ at scales within $\{\delta^{r_{i}}\}_{i=1}^{m+1}\cup \{\delta^{s_{i}}\}_{i=1}^{n+1}\cup [\delta^{r_{i_{1}+1}},\delta^{r_{i_{1}}}]$ can be partitioned into two submeasures which respectively have improved dimensional properties for translates of $B^{\uV_{Q,\theta}}_{\delta^\ur}$ and $B^{\uW_{Q,\theta}}_{\delta^\us}$ in each $\cD_{\eta}$-block $Q$.

%Note that checking the supercritical alternative assumption in the context of random walks will require some work. This will be the essence of \Cref{Sec-RWincrease}

\begin{thm}[Supercritical multislicing decomposition] \label{mult-sup-dec}
Let $d>m, n\geq1$, fix $(\uj, \ur)\in \cP_{m}(d)\times \increasing_{m}$ and $(\uk, \us)\in \cP_{n}(d)\times \increasing_{n}$. Let $\delta, \eps, \eps', \varkappa, c, \tau, \tau'>0$, let $\eta\in [\max(\delta^{r_{1}}, \delta^{s_{1}}), 1]$.

Let $(\Theta, \sigma)$ be a probability space. For each $Q\in \cD_{\eta}$, consider  measurable families $(\uV_{Q,\theta})_{\theta\in \Theta} \in \cF_{\uj}^\Theta$ and $(\uW_{Q,\theta})_{\theta\in \Theta} \in \cF_{\uk}^\Theta$.

For every $Q\in \cD_{\eta}$,  $i=1, \dots, m$, assume the distribution of $(V_{Q,\theta, i})_{\theta\sim \sigma}$  satisfies $\SubP$ with parameter $(\delta^{r_{i+1}}, \eps, \tau)$. Make the corresponding  assumption for the collection $(W_{Q,\theta, i})_{\theta\sim \sigma}$ at scale $\delta^{s_{i+1}}$ for  $i=1, \dots, n$.

Assume that for some subscripts $i_{1}, i_{2}$ we have $r_{i_{1}}=s_{i_{2}}$ and $r_{i_{1}+1}=s_{i_{2}+1}$, and that for every  $Q\in \cD_{\eta}$,   the distributions of $(V_{Q,\theta, i_{1}})_{\theta\sim \sigma}$ and $(W_{Q,\theta, i_{2}})_{\theta\sim \sigma}$ together satisfy $\SAP$ with parameters $(\delta^{r_{i_{1}+1}-r_{i_{1}}}, \varkappa,\tau')$.

Let $\nu$ be a Borel measure on $B^{\R^d}_{1}$ of mass at most $\delta^{-c}$, and such that for some $\alpha\in [\varkappa, 1-\varkappa]$, for all $v\in \R^d$,  all $\rho \in \{\delta^{r_{i}}\}_{i=1}^{m+1}\cup \{\delta^{s_{i}}\}_{i=1}^{n+1}\cup [\delta^{r_{i_{1}+1}},\delta^{r_{i_{1}}}]$, we have
$$\nu(B^{\R^d}_{\rho}+v)\leq\delta^{-c} \rho^{d\alpha}.$$

%$\eps\lll_{d,r_{i_{1}}, r_{i_{1}+1}, C,c,\tau} 1$ and $\delta\lll_{d, \ur, \us, c,  \tau, \eps} 1$,

Let $t_{2}>0$ be the second minimum of $\{r_{i}\}_{i=1}^{m+1}\cup \{s_{i}\}_{i=1}^{n+1}$, and $u:=r_{i_{1}+1} -r_{i_{1}}$. %\marginpar{Contrairement au cas sous-critique, \c ca ne fait pas sens ici d'écrire la condition sur $\delta$ sous la forme $\delta^{t_{2}}\lll_{d, u, \eps,  \tau'} 1$ car $\eps$ est déjà très petit en fonction de $t_{2}$, donc mettre $t_{2}$ en exposant   n'apporte pas d'information}

If $\eps'\lll \eps$; and $\eps, c, \tau  \lll_{d, t_{2}, u, \tau'}1$; and  $\delta \lll_{d, t_{2}, u, \tau', \eps}1$,  then there exists a decomposition
$$\nu=\nu_{1}+\nu_{2}$$
into mutually singular Borel measures, and an event $\cE \subseteq \Theta$ such that $\sigma(\cE)\leq \delta^{t_{2}\eps'}$ and for $p \in \{1, 2\}$,  $\theta\in\Theta \setminus \cE$, there is a set  $F_{p, \theta} \subseteq \R^d$ with  $\nu_{p}(F_{p,\theta}) \leq \delta^{t_{2}\eps'}$ and such that for every $Q\in \cD_{\eta}$, $v \in \R^d$,
\[ \nu_{1 |Q\smallsetminus  F_{1,\theta}} \left(B^{\uV_{Q,\theta}}_{\delta^\ur}+v\right) \leq \delta^{u\tau'/(100d)}   \leb \left(B^{\uV_{Q,\theta}}_{\delta^\ur}\right)^{\alpha},\]
 while $ \nu_{2 |Q\smallsetminus  F_{2,\theta}}$ satisfies the analogous bound with $(\us,\uW_{Q,\theta})$ in the place of $(\ur, \uV_{Q,\theta})$.
\end{thm}

\begin{remark}
The implicit constant in the upper bound $\delta\lll_{d, t_{2}, u, \tau', \eps} 1$ only depends on $d$ and a positive lower bound on $t_{2}$, $u$,  $\tau'$, $\eps$.
\end{remark}

\bigskip
The proofs of \Cref{subc-mult} and \Cref{mult-sup-dec}  are similar to those in \cite[Section 2]{BH24}. We postpone them to \Cref{ss:Appendix-Multislicing}.

\section{Optimal subcritical projection theorem} \label{Sec-sub-proj} \label{Sec-optimal-subcritical}

This section can be read independently of the rest of this paper.
We consider a probability measure $\sigma$ on $\Gr(\R^d,k)$ for  fixed integers $d>k\geq1$. It defines  a random orthogonal projector $(\pi_{L})_{L \sim \sigma}$. 
We wish to find a criterion on $\sigma$ to guarantee that for any  set $A \subset B^{\R^d}_{1}$ of dimension at least $s \in [0, d]$, for most realizations of $L\sim \sigma$,  the dimension of $\pi_L A$ is at least $\frac{k}{d}s$, up to an arbitrary small loss. 

%We consider a random subspace $L\subseteq \R^d$ of fixed rank $k\geq 1$. 
%We write $\pi_L : \R^d \to \R^d$  the orthogonal projector of image $L$. 
%We study whether for any bounded subset $A \subset \R^d$ of dimension $s \in [0, d]$ at some
%scale $\delta > 0$, most of its projections $\pi_L A$ have dimension roughly at least $\frac{k}{d}s$ at scale $\delta$.

There are obvious linear obstructions.
Indeed, consider a subspace $W\subseteq \R^d$.
If $\sigma$ is supported on the constraining pencil\footnote{Using the terminology of \cite{ABRS}.}
\begin{equation}
\label{eq:constrPW}
\cP^W=\setBig{ L \in \Gr(\R^d,k) \,:\, \dim( \pi_L W) < \frac{k}{d} \dim W },
\end{equation}
then taking $A$ to be the unit ball in $W$, every projection $\pi_L A$ is of dimension less than the desired threshold.

The main result of this section, recorded below as \Cref{subcritical-projection}, states in a quantitative way that these linear obstructions are the only obstructions. It is presented in  a discretized form, i.e. in terms of covering numbers at a fixed small scale. A limiting version in terms of Hausdorff dimension is recorded in \Cref{cor-haus-dim}.
In the rest of the paper,  \Cref{subcritical-projection}
will be crucial to check the subcritical assumptions in the multislicing theorems from  \Cref{Sec-multislicing-thms}.  
%A limit version, in terms of Hausdorff dimension, is derived in \Cref{cor-haus-dim}. 

%In this section, we consider a random orthogonal projector $(\pi_{L})_{L\sim \sigma}$, a set $A\subseteq \R^d$, and we give a subcritical lower bound for the covering number of the image $\pi_{L}(A)$ by balls of radius $\delta$, provided a natural non-concentration assumption on $\sigma$. This result gives a concrete way of checking the subcritical assumptions in the multislicing theorems from .

%We work under the assumption that given any $W\in \Gr(\R^d)$, there is a high probability that the intersection $L^\perp\cap W$ has smaller or equal proportion in $W$ that $L^\perp$ in $\R^d$. This assumption is optimal. It improves upon a previous version of the subcritical projection theorem \cite{BH24, He2020JFG} which required the stronger condition that $L^\perp$ is typically in direct sum with any subspace $W$ of complementary dimension.
\bigskip

Recall from \S\ref{conventions-notations} that we have fixed a distance on $\Gr(\R^d)$.
For $W\in \Gr(\R^d)$, $\rho>0$,  the notation $B_\rho(W)$ stands for the open ball of radius $\rho$ and center $W$. By convention, every subspace $W'\in B_{1}(W)$ satisfies $\dim W'=\dim W$. We introduce a thickening of the constraining pencil $\cP^W$. It is defined  for  $\rho\in (0, 1)$ by 
\[
\cP_\rho^W = \setBig{ L \in \Gr(\R^d,k) \,:\, \exists W' \in B_\rho(W),\; \dim( \pi_L W') < \frac{k}{d} \dim W },
\]
or equivalently 
\[
\cP_\rho^W = \setBig{L \in \Gr(\R^d,k) \,:\, \exists W' \in B_\rho(W),\; \dim (L^\perp \cap W') > \frac{d-k}{d} \dim W}.
\]
We show
\begin{thm}[Subcritical projection theorem] \label{subcritical-projection}
Let $d>k \geq 1$ be integers. Let $D>1$, let $\kappa, \eps, \rho, \delta \in (0, 1)$ with $\rho\geq \delta$.
Let $\sigma$ be a probability measure on $\Gr(\R^d, k)$ satisfying
%Assume that for every subspace $W\in \Gr(\R^d)$, we have
\begin{equation} \label{avoiding-pencils}
\forall W \in \Gr(\R^d),\quad
\sigma\left(\cP_{\rho}^W\right)  \leq \rho^{\kappa}.
\end{equation}

If $D \ggg_{d}1 +\frac{\eps}{\kappa} \abse{\frac{\log \delta}{\log \rho}}^2   $; $\delta\lll_{\eps}1$; and  $\rho \leq \delta^{4d^2\eps/\kappa}$, then  for every set $A\subseteq B^{\R^d}_{1}$, the exceptional set 
\begin{equation}
\label{eq:Eset}
\begin{split}
\cE \defeq \{\, L\, :\, \exists A' \subseteq A \,\,&\text{ with }\,\,\cN_{\delta}(A')\geq \delta^{\eps} \cN_{\delta}(A)\\
& \text{ and } \,\, \cN_{\delta}(\pi_{L}A') < \rho^{D} \cN_{\delta}(A)^{\frac{k}{d}} \}
\end{split}
\end{equation}
satisfies $\sigma(\cE)\leq \delta^{\eps}$.
\end{thm}

\Cref{subcritical-projection} improves upon a previous version of the subcritical projection theorem \cite[Proposition 29]{He2020JFG} (see also \cite[Proposition A.2]{BH24}) which required the stronger condition that $L$ is typically in direct sum with any subspace $W$ of complementary dimension, or in other words, that $L$ avoids all pencils, not only the constraining ones.
In this regard, \Cref{subcritical-projection} is \emph{optimal}, since constraining pencils are indeed obstructions.

\begin{remark}
In the particular case where $\rho=\delta^{\sqrt{\eps}}$, the lower bound on the exponent $D$ only depends on $d, \kappa$, namely one can take 
$D=O_{d}(\kappa^{-1})$.  With the terminology of \Cref{Sub-crit-P}, the conclusion then means that the distribution of $L^\perp$ as $L\sim \sigma$ satisfies the subcritical property $\SubP$ with parameters $(\delta, \eps, D\sqrt{\eps})$.
\end{remark}

%\begin{remark} An obstruction to subcritical estimates is the following: if $\dim L^\perp \cap V > \frac{d-k}{d}\dim V$, then taking $A=B^V_{1}$, we have a macroscopic loss of dimension:
%$$\cN_{\delta}(\pi_{L}A)= \cN_{\delta}(A)^{\frac{\dim V-\dim L^\perp\cap V}{\dim V}}\leq \cN_{\delta}(A)^{\frac{k-1}{d}}.$$
%\Cref{subcritical-projection} assumes (via \eqref{avoiding-pencils}) this obstruction does not occur for most $V\sim \sigma$, and guarantees subcritical estimates. In this sense, our theorem is \emph{optimal}. 
%\end{remark}

%\begin{remark} The role of the non-concentration assumption \eqref{avoiding-pencils} is to prevent the obstruction that if $\dim L^\perp \cap V > \frac{d-k}{d}\dim V$, then taking $A=B^V_{1}$, we have a macroscopic loss of dimension:
%$$\cN_{\delta}(\pi_{L}A)= \cN_{\delta}(A)^{\frac{\dim V-\dim L^\perp\cap V}{\dim V}}\leq \cN_{\delta}(A)^{\frac{k-1}{d}}.$$
%Our theorem is optimal in the sense.
%\end{remark}

\begin{remark}
Assumption \eqref{avoiding-pencils} is invariant by replacing $k$ by $d - k$ and $\sigma$ by its image under $L\mapsto L^\perp$. Indeed, this follows from the fact that the distance on the Grassmannian is invariant under taking the orthogonal (see \Cref{grasmannian-distance}), combined with \Cref{orthog-red-sub} below.
\end{remark}

\begin{lemma}
\label[lemma]{orthog-red-sub}
Let $E, F$ be subspaces of a given real Euclidean vector space $T$. Then the relation
$$\dim E \dim F \geq \dim T \dim E\cap F$$
is equivalent to its orthogonal counterpart
$$\dim E^\perp \dim F^\perp \geq \dim T \dim E^\perp \cap F^\perp.$$
\end{lemma}

\begin{proof}
Set $e=\dim E$, $f=\dim F$, $t=\dim T$, $s=\dim E+F$. The first relation can be written $ef\geq t(e+f-s)$. Note that
$\dim E^\perp \cap F^\perp= t- \dim (E^\perp \cap F^\perp)^\perp=t-s$. Hence the second relation can be written $(t-e)(t-f)\geq t(t-s)$. Both relations are then clearly equivalent.
\end{proof}

%Taking $\rho=\delta^{\sqrt{\eps}}$ with $\eps$ arbitrarily small, and letting $\delta$ go to zero, 
\Cref{subcritical-projection} implies a corresponding statement for the Hausdorff dimension of analytic sets. %but it is of independent interest and  also enlightening for understanding \Cref{subcritical-projection}.
%\comW{To be checked.}
\begin{corollary}
\label[corollary]{cor-haus-dim}
Let $d>k \geq 1$ be integers. Let $A \subset \R^d$ be an analytic set. %of Hausdorff dimension $\dimH A = s$.
The set of exceptional directions 
\[
\cE_{\mathrm{H}}(A) = \set{ L \in \Gr(\R^d, k) \,:\, \dimH (\pi_L A) < \frac{k}{d} \dimH A}
\]
does not support any non-zero Borel measure $\sigma$ satisfying
\[
\exists \kappa > 0,\; \forall \rho > 0,\; \forall W \in \Gr(\R^d), \quad \sigma( \cP_\rho^W) \leq \rho^\kappa.
\]
\end{corollary}
Although it will not be used in the rest of the paper, \Cref{cor-haus-dim} is interesting in its own right. %, and in view of the current literature.
It implies\footnote{Together with Frostman's Lemma.} for example an estimate on the Hausdorff dimension of the exceptional set:
\[
\dimH \cE_{\mathrm{H}}(A) \leq \dim \Gr(\R^d,k) - \min\{k, d-k\},
\]
which is precisely \cite[Theorem 1]{Gan24}.

\bigskip
The proof of \Cref{subcritical-projection} relies on effective Brascamp-Lieb inequalities, which take the form of a visual inequality presented below.
Those inequalities are established in our companion paper \cite{BH25-BrascampLieb}.
The strategy to use  Brascamp-Lieb inequalities in order to derive a lower bound on the dimension of a  projected set is inspired by  \cite{Gan24}.

\subsection{Visual inequality}
We start by stating the precise input we need from \cite{BH25-BrascampLieb}.

Let $J\in \N^*$, and consider a collection
$$\sD=((\pi_{L_{j}})_{1\leq j\leq J}, (q_{j})_{1\leq j\leq J})$$ where
$L_{j}\in \Gr(\R^d)\smallsetminus \{0\}$, $\pi_{L_{j}}$ is the orthogonal projector of image $L_{j}$, $q_{j}>0$. Assume they together satisfy 
$$\sum_{j=1}^Jq_{j}\dim L_{j}=d.$$

\begin{definition}[Perceptivity]
\label[definition]{def-BL}
Given $\alpha\in (0, 1], \beta\in \R^+$,
we say the datum $\sD$ is \emph{$(\alpha, \beta)$-perceptive}\footnote{This terminology diverges slightly from that in \cite[Equation (9)]{BH25-BrascampLieb}.
However, for $\alpha\lll_{d}1$, \cite[Lemma 2.5]{BH25-BrascampLieb} implies that $(\alpha,\beta)$-perceptivity in our sense implies $(O_{d}(\alpha),\beta)$-perceptivity in the sense of \cite{BH25-BrascampLieb} and conversely. }
 if for all non-zero $W \in \Gr(\R^d)$,
\begin{equation}\label{eq-percep-proportion}
\sum_{j=1}^J q_{j} \max_{W'\in B_{\alpha}(W)}\frac{\dim L_{j}^{\perp} \cap W'}{\dim W} \,\,\leq\,\,\frac{\beta}{\dim W} \,+\, \sum_{j=1}^J q_{j} \frac{\dim L^\perp_{j}}{d}.
 \end{equation}
\end{definition}

For $\beta=0$, perceptivity expresses that, in average, the orthogonal subspaces $L_{j}^{\perp}$ fill up (proportionally) less $W$ than the whole space $\R^d$.
It actually allows for some perturbations of $W$, which is a way to say that in average the $L_{j}^{\perp}$'s have a large subspace making a large angle with $W$.

\bigskip
The following is a special case of \cite[Theorem 1.6]{BH25-BrascampLieb}. 
\begin{proposition}[Visual inequality]
\label[proposition]{thm-visual-ineq}
Let $\sD=((\pi_{L_{j}})_{1\leq j\leq J}, (q_{j})_{1\leq j\leq J})$ be as above. Assume $\sD$ is $(\alpha, \beta)$-perceptive for some $\alpha\in (0, 1]$, $\beta\in \R^+$. Then for every $\delta\in (0,1)$, every subset $A\subseteq B^{\R^d}_{1}$, we have
\begin{equation} \label{eq-thm-visual-ineq}
 \cN_{\delta}(A)\leq C  \delta^{-\beta} \alpha^{-d} \prod_{j=1}^J \cN_{\delta}(\pi_{L_{j}}A)^{q_{j}}
 \end{equation}
where $0< C \leq e^{O_d(1+\sum_{j}q_{j})} (1+\sum_{j}q_{j})^{\frac{\beta}{2}} \prod_{j}q_{j}^{-q_{j}\dim L_{j}/2} $.
\end{proposition}

This inequality can be seen as a generalization of the trivial inequality that for any finite set $A\subseteq \R^d$, any basis $(v_{1}, \dots, v_{d})$ of $\R^d$, one has
$$|A|\leq \prod_{i=1}^d|\pi_{\R v_{i}}A|.$$

%\begin{proof}
%The proof is established in \cite{BH25-BrascampLieb}. It arises from a quantitative bound for Brascamp-Lieb constants, applied with a suitable choice of test functions.
%\end{proof}

\subsection{Proof of the subcritical projection theorem}

Let $d>k\geq 1$ be integers. Given a finite collection  $\bL = (L_j)_{j}\in \Gr(\R^d, k)^J$ of $k$-planes in $\R^d$,  consider the datum
\begin{equation}\label{sDLj}
\sD_\bL:=\left((\pi_{L_{j}})_{1\leq j\leq J}, \left(\frac{d}{kJ}\right)^{\otimes J} \right).
\end{equation}
In the next lemma, we assume the $L_{j}$'s are chosen randomly and independently via a probability measure $\sigma$ on $\Gr(\R^d, k)$ which is not concentrated near constraining pencils. 
%We select randomly and independently with law $\sigma$ a finite collection of subspaces $L_{j}\in \Gr(\R^d, k)$ ($1\leq j\leq J$), and assign to each $L_{j}$ the weight  $q_{j}:=\frac{d}{kJ}$.
 We then obtain a lower bound on the probability that the associated datum  be perceptive.
%We place ourselves in the context of \Cref{subcritical-projection}. Aiming to apply \Cref{thm-visual-ineq}, we first show that by choosing a finite collection of subspaces $(L_{j})_{1\leq j\leq J}$ with law $\sigma^{\otimes J}$, and assigning weight $q_{j}:=\frac{d}{kJ}$ to each of them, we obtain with high probability a $(\delta^{C\sqrt{\eps}}, \beta)$-perceptive datum.

\begin{lemma}
\label[lemma]{generating-perceptive}
Let $\alpha, \gamma \in (0,1]$. Let $\sigma$ be a probability measure on $\Gr(\R^d, k)$ satisfying 
\begin{equation}\label{avoid-pencils2}
\forall W \in \Gr(\R^d), \quad \sigma\bigl(\cP_{2\alpha}^W\bigr) \leq \gamma.
\end{equation}
Then for every $J \geq 1$, $\beta > 0$,
\[
\sigma^{\otimes J}\setbig{\bL \,:\, \text{$\sD_{\bL}$ is not $(\alpha, \beta)$-perceptive}} \ll_d O(1)^J \alpha^{- \dim \Gr(\R^d)} \gamma^{J \beta/d}.
\]
\end{lemma}

The proof combines the non-concentration assumption \eqref{avoid-pencils2} with Chernoff's additive tail bound for sum of i.i.d. Bernoulli variables. We recall the latter.

\begin{lemma}[Chernoff's bound]
\label[lemma]{Chernoff}
 Let $J\geq 1$, let $Z_{1}, \dots, Z_{J}$ be i.i.d. Bernoulli random variables.
Then for any $t\in \R^+$,
\[
\bbP\left[\frac{1}{J}\sum_{j=1}^JZ_{j}\geq t \right]^{1/J} \ll \bbP[Z_{1} = 1]^{t}.
\]
\end{lemma}
\begin{proof}
We record a short proof from \cite{Chvatal79}.
Write $p = \bbP[Z_1 = 1]$. 
Note one may assume $t\in (p, 1)$.
Set $k:=\lceil tJ \rceil$. Given $s>1$,
$\bbP\left[\sum_{j=1}^JZ_{j}\geq tJ \right] =\sum_{i=k}^J \binom{J}{i}p^i(1-p)^{J-i} \leq \sum_{i=0}^J \binom{J}{i}p^i(1-p)^{J-i}s^{i-k}=s^{-k}(sp+(1-p))^J$.
Plugging $s=\frac{t(1-p)}{p(1-t)}$ and using $s^{-k}\leq s^{-tJ}$, $\sup_{r\in (0, 1)} |r\log r|<\infty$, the bound follows.
\end{proof}

%We give a short proof extracted from \cite{chvatal 1979 tail of hypergeometric distribution end of the article}
%$$\mathbb{P}\left(\frac{1}{J}\sum_{j=1}^JZ_{j}\geq t \right)\leq \exp\left(-J \DD(t \,\|\, p)\right) $$
%where $$\DD(t \,\|\, p)= t\log \frac{t}{p} + (1-t) \log \frac{1-t}{1-p}.$$
%In particular, using that $\sup_{s\in (0, 1)} |s\log s|<\infty$, we deduce
%$$\mathbb{P}\left(\frac{1}{J}\sum_{j=1}^JZ_{j}\geq t \right)^{1/J}\ll p^{t}.$$

\begin{proof}[Proof of \Cref{generating-perceptive} ]
Let $L_1, \dotsc, L_J$ be i.i.d. random variables taking value in $\Gr(\R^d,k)$ and following the law $\sigma$.
Writing $\bL = (L_1, \dotsc, L_J)$, we bound from above the probability of the event (denoted by $\Obs$) that $\sD_{\bL}$ is not $(\alpha, \beta)$-perceptive.
By definition, we have 
\[
\Obs = \bigcup_{W \in \Gr(\R^d)\smallsetminus \{0\}} \Obs_\alpha^W,
\]
where $\Obs_\alpha^W$ is the event that there exists $(W_j)_j \in B_\alpha(W)^J$ such that
\[
\frac{1}{J} \sum_{j = 1}^J \frac{d}{k} \left(\frac{ \dim ( L_j^\perp \cap W_j )}{ \dim W} - \frac{d-k}{d} \right) > \frac{ \beta}{\dim W}.
\]
Clearly $\Obs_\alpha^{W'} \subset \Obs_{2\alpha}^W$ for any $W' \in B_\alpha(W)$.
Covering $\Gr(\R^d)\smallsetminus \{0\}$ by $O_d(\alpha^{- \dim \Gr(\R^d)})$ balls of radius $\alpha$, we obtain
\begin{equation} \label{eqd}
\bbP[\Obs] \ll_d \alpha^{- \dim \Gr(\R^d)} \sup_{W \in \Gr(\R^d)\smallsetminus \{0\}} \bbP\left[\Obs_{2 \alpha}^W\right].
\end{equation}

We now bound the probability of $\Obs_{2 \alpha}^W$ for a given non-zero $W\in  \Gr(\R^d)$.
First, observing the relation $\frac{d}{k}\bigl(1-\frac{d-k}{d}\bigr)=1$ and recalling the definition of $\cP_{2 \alpha}^W$, we have for each $j \in \{1,\dotsc,J\}$,
\[
\frac{d}{k} \max_{W' \in B_{2\alpha}(W)} \left(\frac{ \dim ( L_j^\perp \cap W' )}{ \dim W} - \frac{d-k}{d} \right) \leq Z_j,
\]
where $Z_j = \1_{\cP_{2 \alpha}^W}(L_j)$.
Therefore,
\[
\bbP\left[\Obs_{2 \alpha}^W \right] \leq \bbP\left[\frac{1}{J}\sum\nolimits_{j = 1}^J  Z_j > \frac{\beta}{\dim W} \right] \leq 
\bbP\left[\frac{1}{J}\sum\nolimits_{j = 1}^J Z_j > \frac{\beta}{d} \right].
\]
Note that the $(Z_j)_j$ are i.i.d. Bernoulli random variables with $\bbP[Z_j=1] = \sigma(\cP_{2 \alpha}^W) \leq \gamma$, therefore we can use \Cref{Chernoff} to obtain
\[
\sup_{W \in \Gr(\R^d)\smallsetminus \{0\}} \bbP\left[\Obs_{2 \alpha}^W \right] \leq O(1)^J\gamma^{ \beta J/d }.
\]
Together with \eqref{eqd}, this gives the desired estimate.
\end{proof}

We shall also make use of the following lemma. It guarantees that i.i.d. random events have a reasonable chance to occur simultaneously.

\begin{lemma}
\label[lemma]{lem:intersection}
Let $(\Omega,\bbP)$, $(A, \lambda)$ be two probability spaces, let $(A_{\omega})_{\omega\in \Omega}$ be a measurable\footnote{Measurability means that the map $\Omega\times A\rightarrow \R, (\omega, x)\mapsto \1_{A_{\omega}}(x)$ is measurable.}
 collection of subsets of $A$. Assume $\inf_{\omega\in \Omega} \lambda(A_{\omega})\geq t$ where $t\in (0, 1)$. Then for every integer $J\geq 1$,
$$\bbP^{\otimes J}\setbig{(\omega_{j})_{1\leq j\leq J} \,:\, \lambda(\cap_{j} A_{\omega_{j}}) \geq t^{J}/2 } \geq t^{J}/2.$$
\end{lemma}

%\begin{lemma}\label{lem:intersection}
%Let $(\Omega,\bbP)$, $(A, \lambda)$ be two probability space, let $A\subseteq \R^d$ be a finite set, and let $(A_{\omega})_{\omega\in \Omega}$ be a measurable collection of subsets of $A$. Assume for some $r\in (0, 1/2)$ that $\inf_{\omega\in \Omega}|A_{\omega}|\geq r|A|$. Then for every integer $J\geq 1$,
%$$\bbP^{\otimes J}\{(\omega_{j})_{1\leq j\leq J} \,:\,|\cap_{j} A_{\omega_{j}}|\geq r^{2J} |A| \}\geq r^{2J}.$$
%\end{lemma}

\begin{proof}
It follows by applying Markov's inequality, then Fubini's theorem, and Jensen's inequality. See \cite[Lemma 19]{He2020JFG} for details.
\end{proof}

We are now able to conclude the proof of the subcritical projection theorem.

\begin{proof}[Proof of \Cref{subcritical-projection}]
We may suppose that $A$ is $2\delta$-separated, hence finite.
Let $\sP(A)$ denote the collection of subsets of $A$, endowed with the discrete $\sigma$-algebra.

Assume for a contradiction that $\sigma(\cE)>\delta^\eps$.
%Replacing $\sigma$ by its restriction $\sigma_{\mid \cE}$ on $\cE$, and replacing the condition \eqref{avoiding-pencils} by
%\begin{equation} \label{avoiding-pencils-weak}
%\forall W \in \Gr(\R^d),\quad
%\sigma\left(\cP_{\rho}^W\right)  \leq \delta^{-\eps} \rho^{\kappa},
%\end{equation}
For every $L\in \cE$, there is a subset $A_L \subseteq A$ such that 
\[
\abs{A_{L}}\geq \delta^{\eps} \abs{A}\quad \text{and} \quad \cN_{\delta}(\pi_L A_L) < \rho^{D}\abs{A}^{\frac{k}{d}}.
\]
Note that the same set $A_L$  can serve as $A_{L'}$ for every $L'$ sufficiently close enough to $L$.
Hence we may choose the map $\cE\rightarrow \sP(A),  \,L\mapsto A_{L}$ to be measurable on $\cE$. We then extend it arbitrarily into a measurable map on $\Gr(\R^d,k) \rightarrow \sP(A), \,L\mapsto A_{L}$.

We consider parameters $J\in \N^*$ and $\beta > 0$ to specify below.
Let $L_1, \dotsc, L_J$ be i.i.d. random variables following the law $\sigma$.
Write $\bL = (L_1,\dotsc, L_J)$ and set $A_\bL = \cap_{j} A_{L_j}$.
By \Cref{lem:intersection} applied to the probability measure $\sigma(\cE)^{-1} \sigma_{\mid \cE}$ and the uniform probability measure on $A$, we know that the event 
\begin{equation}
\label{eq:Badset}
\abs{A_\bL} \geq 2^{-1}\delta^{J \eps} \abs{A} \quad \text{and} \quad \forall 1 \leq j \leq J,\, \cN_\delta(\pi_{L_j} A_\bL) < \rho^{D}\abs{A}^{\frac{k}{d}}
\end{equation}
happens with probability at least $\delta^{2 J \eps}/2$.

On the other hand, let $\sD_\bL$ be as in \eqref{sDLj}.
Then \Cref{generating-perceptive} implies that the event that
\begin{equation}
\label{eq:Percepset}
\sD_\bL \text{ is $(\rho/2, \beta)$-perceptive}
\end{equation}
happens with probability at least $1 - \rho^{- d^3 + \kappa J \beta / d^2}$, provided $\rho\lll_{d}1$ and $\rho^{\kappa  \beta / d}\lll 1$.
%happens with probability at least $1 - O_d\bigl(\rho^{-O_d(1) + \kappa J \beta / d} \bigr)$.

Now, choose 
\(
\beta = \frac{4 d^2 \eps}{\kappa} \frac{\log \delta}{\log \rho}
\)
so that 
\(
\rho^{\kappa J \beta / d^2} = \delta^{4 J \eps}
\)
and then choose $J = \left\lceil 3 \kappa^{-1} \beta^{-1} d^5\right\rceil$ so that $\kappa J \beta / d^2 \geq 3 d^3$. Assume $\delta \lll_{\eps}1$.
Then the lower bounds on the probability of the events \eqref{eq:Badset} and \eqref{eq:Percepset} imply that they  happen simultaneously with non-zero probability. We may thus consider a realization of $\bL$ satisfying both  \eqref{eq:Badset} and \eqref{eq:Percepset}. Assume $\rho\leq \delta^{\frac{4 d^2 \eps}{\kappa}}$ so that $\beta\leq 1$.
Invoking the visual inequality from \Cref{thm-visual-ineq}, we obtain
\begin{align*}
2^{-1}\delta^{J \eps} \abs{A} \leq \abs{A_\bL}& \ll_d J^{d/2} \delta^{-\beta} \rho^{-d} \prod_{j =1}^J \cN_{\delta}(\pi_{L_j} A_{L_j})^{\frac{d}{kJ}}\\
& \ll_d J^{d/2} \delta^{-\beta} \rho^{d D / k - d} \abs{A}.
\end{align*}
Provided $\rho\leq \eps$, and noting $J\eps\leq \frac{\log \rho}{\log \delta} d^3=O_{d}(1)$, this implies by direct computation
\[
D \leq \frac{\log \delta}{\log \rho}  \beta+ O_{d}(1).
\]
Hence, we obtain a contradiction if 
\[
D - \frac{\log \delta}{\log \rho}  \beta \ggg_{d}1.\qedhere
\]
\end{proof}

\section{Submodular inequality for irreducible representations} \label{Sec-submod}

This section can be read independently of the rest of this paper.
Its goal  is to establish a general submodular inequality for strongly irreducible group representations over arbitrary fields, \Cref{submod-G-F}.
 We also study the equality case, see  \Cref{submod-G-F-2}.

In the context of random walks, \Cref{submod-G-F} will be used to justify \Cref{mu-Sub}, which checks that  a translate $gB_{\rho}x$ where $g\sim \mu^n$, $\rho>0$, $x\in X$ looks like a random box, whose associated partial flag satisfies the non-concentration estimates relevant to the subcritical projection theorem.
As such, \Cref{submod-G-F} is a crucial ingredient for proving the properties of dimensional interpolation and supercritical decomposition for the action of random walks on homogeneous spaces (namely \Cref{RW-preserve-dim} and \Cref{sup-mult-X}, or their simplified versions Propositions \ref{pr:dim-quasi-pres}, \ref{pr:splitnu} from \Cref{Sec-roadmap}).

\bigskip
\begin{thm}[Submodular inequality]
\label{submod-G-F}
Let $\cH$ be a finite-dimensional vector space over some field $\F$.
Let $G \subset \GL(\cH)$ be a  Zariski-connected  subgroup  whose action on $\cH$ is irreducible.
Then for all subspaces $V, W\subseteq \cH$, one has
\begin{equation}\label{submod-ineq-GF}
\min_{g \in G} \dim(g V \cap W)  \leq \frac{\dim V \dim W}{\dim \cH}.
\end{equation}
%Moreover, the set $\cL_{W}\subseteq \Gr(\R^d)$ of subspaces $V$ such that equality holds in \eqref{submod-ineq-GF} is a $G$-invariant sublattice and every $V\in \cL_{W}$ satisfies
%\begin{equation}\label{submod-eq-GF}
%\dim(V \cap W)  = \frac{\dim V \dim W}{\dim \cH}.
%\end{equation}
\end{thm}

%We recall the \emph{strong irreducibility} assumption of $G\acts \cH$ means that every $G$-invariant finite union of non-zero subspaces of $\cH$ must be equal to $\cH$.Equivalently, the identity component of the Zariski-closure of $G$ within $\GL(\cH)$ acts irreducibly on $\cH$.

We emphasize that $\F$ may have positive characteristic, and is potentially not algebraically closed. The case of interest for our application to random walks is that of $\F=\R$ and $G\acts \cH$ is the adjoint representation of a simple connected real\footnote{In fact, the complex case is also sufficient.} Lie group. Note also that $G$ is potentially not Zariski-closed, although it is easy to reduce to that case.

\begin{remark}
It is an immediate consequence that if $\cH$ is a strongly irreducible linear representation of some group $G$, then the same inequality holds.
\end{remark}

\begin{remark}
Equation \eqref{submod-ineq-GF} can be seen as a \emph{transversality principle}, stating that  $gV$ and $W$ cannot intersect too much for $g$ generic. More precisely, the orbit $G \cdot V$ cannot be contained in any constraining pencil (see \eqref{eq:constrPW} and \Cref{orthog-red-sub}).
Recall that such an  avoidance condition has already played a role as an assumption for the subcritical projection theorem, \Cref{subcritical-projection}.
\end{remark}

\begin{remark}
Multiplying both sides of \eqref{submod-ineq-GF} by $\dim \cH$, Equation \eqref{submod-ineq-GF} takes approximately the form of a \emph{multiplicative submodular inequality}, with $\cH$ playing the role of the ``join'' of  $gV$ and $W$.
In the proof, some submodularity in the sense of lattice theory will be a key ingredient.
%Finally, dividing \eqref{submod-ineq-GF} by $\dim W$, it can be interpreted as a \emph{scarcity principle}, saying that a $gV$ becomes scarcer in restriction to  $W$. Scarcity under restriction has already played a role in \Cref{Sec-sub-proj}, through the notion of perceptiveness, and as an assumption for the subcritical projection theorem, \Cref{subcritical-projection}.
\end{remark}

\begin{remark} The theorem fails if we replace the irreducibility assumption by the weaker condition that the $G$-orbits $G\cdot V$ and $G\cdot W$ span $\cH$. Indeed, consider for instance $G\acts \cH$ given by $\GL(\R^k)\times \GL(\R^l) \acts \R^k\oplus \R^l$ ($k,l\in \N^*$), and $V=W=\R^k\oplus \R e_{1}$. Then $\min_{g}\dim(gV\cap W)=k$ while $\dim V\dim W/\dim \cH=(k+1)^2/(k+l)$ can be much smaller than $k$ for large $l$.
\end{remark}

We  illustrate \Cref{submod-G-F}  with a few examples.
\begin{example}
For $F=\R$ and $G$ isomorphic to $\SL_{2}(\R)$, \Cref{submod-G-F} is trivial. In fact, representation theory of $\SL_{2}(\R)$ then yields the stronger result: 
\begin{equation}
\label{eq:dimVdimW-sl2case}
\min_{g\in G}\dim(gV \cap W) = \max\{0, \dim V + \dim W - \dim \cH \}.
\end{equation}
In other words, $gV$ (for generic $g$) and $W$ do not intersect unless they have to because of the Grassmann formula. 
\end{example}

\begin{example}
Consider the case where $G\acts \cH$ is $\Ad(\SL_{d}(\R))\acts \sl_{d}(\R)$. Set $V=\Vect\set{E_{i,d}\,:\, 1 \leq i \leq {d-1} } \subseteq \sl_{d}(\R)$ given by the last column, and  $W=\Vect\set{E_{i,j}\,:\, d-1 \leq i \leq d,\, 1 \leq j \leq d } \cap \sl_{d}(\R)$ given by the last two rows.
\Cref{submod-G-F}  predicts $\min_{g} \dim(gV \cap W) \leq (d-1)(2d-1)/(d^2-1)$, while the reality is $\min_{g} \dim(gV \cap W) = 1$.%, while $\dim V=d-1$, and $\dim W=2d-1$.
\end{example}

\begin{example}[Equality case]
\label[example]{Ex-eg}
If either $V$ or $W$ is the zero subspace or $\cH$, \eqref{submod-ineq-GF} holds trivially as an equality.
As observed in \cite[Remark 1.8]{Lin25b}, nontrivial equality may arise in \eqref{submod-ineq-GF}.
Indeed, assume $\F$ algebraically closed, and for $i=1,2$ consider a $\F$-vector space $\cH_{i}$, and a subgroup $G_{i} \subseteq \GL(\cH_{i})$  which is Zariski-connected and acts irreducibly on $\cH$.
Then $G:=G_{1}\times G_{2}$ is also Zariski-connected and acts irreducibly on $\cH:=\cH_{1}\otimes  \cH_{2}$. 
Moreover, given any $E_{1}\in  \Gr(\cH_1)$, $E_{2}\in \Gr(\cH_2)$, setting $V:=\cH_{1}\otimes E_{2}$, $W:=E_{1}\otimes \cH_{2}$ gives $\dim(gV\cap W) =\dim V \dim W/\dim \cH$ for all $g\in G$. 
\end{example}

\bigskip

We obtain additional rigidity properties when equality occurs in \eqref{submod-ineq-GF}. To present those properties, we need the notion of Grassmannian lattice.
We call \emph{sublattice} of $\Gr(\cH)$  any  collection $\cL$ of subspaces of $\cH$ which is stable by  sum and  intersection, and contains\footnote{This last requirement differs from standard terminology.} both $\{0\}$ and $\cH$.
We say $\cL$ is \emph{$G$-invariant} if for any $g \in G$, $V \in \cL$, we have $g V \in \cL$.
We say $\cL$ is \emph{complemented} if for every $V\in \cL$, there exists $V'\in \cL$ satisfying $\cH=V \oplus V'$.
Below, all sublattices will be $G$-invariant, whence  complemented (see \Cref{split-lattice}). 

\begin{thm}[Submodular inequality - equality case] \label{submod-G-F-2}
Let $\F,\cH,G, V,W$ be as in \Cref{submod-G-F}.
Denote by $\cL_{V},\cL_{W}$ the $G$-invariant sublattices of $\Gr(\cH)$ generated by $V$ and $W$ respectively.
If equality holds in $\eqref{submod-ineq-GF}$, then for every $V'\in \cL_{V}$, $W'\in \cL_{W}$, one has
$$\dim(V'\cap W') = \frac{\dim V' \dim W'}{\dim \cH}.$$
Moreover, $\cL_{W}$ is complemented, and given any splitting $\cH=W'\oplus W''$ where $W',W''\in \cL_{W}$, one has for all $V'\in \cL_{V}$, 
$$V'=(V'\cap W') \oplus (V'\cap W''). $$
\end{thm}

Before passing to the proof, we explain how our results above relate to other works. 
\bigskip

\noindent\emph{Related works}. 
For \Cref{submod-G-F}, the case where $\F=\Q$ is due to Eskin–Mozes–Shah \cite[Corollary 1.4]{EMS97},  who deduced submodularity from non-divergence estimates for translates of algebraic measures in homogeneous spaces. 
In a former arXiv version of the present paper, \Cref{submod-G-F} was presented in the case where $\F=\C$ and $G\acts \cH$ is the adjoint representation of a simple complex Lie group. This case is sufficient to derive all the main results of the paper. The proof relied on a  case by case analysis, using the classification of simple complex Lie algebras, and exploiting certain convexity properties.
Afterwards, an independent proof of \Cref{submod-G-F} appeared in Lin~\cite[Theorem 1.7]{Lin25b}; introducing a notion of stratified set, and arguing inductively on taking sums and intersections, see \cite[Sections 10,11]{Lin25b}.
We have been informed by Ruotao Yang of yet another proof of the case $\F = \C$, which uses  Weyl's unitary trick and an averaging argument.
Neither \cite{EMS97, Lin25b} give any description of the equality case.
We establish \Cref{submod-G-F}  via a  different approach, which has the double advantage of being more direct than the above and to yield meaningful information regarding the equality case, see  \Cref{submod-G-F-2}.
The   readers may notice a resonance with Petridis' proof~\cite{Petridis}\footnote{See also Terence Tao's comment in this blog post~\url{https://gowers.wordpress.com/2011/02/10/a-new-way-of-proving-sumset-estimates} by Timothy Gowers.} of the Plünnecke-Ruzsa inequality,  and the theory of Harder-Narasimhan filtrations (see e.g.~\cite[Lemma 2]{BS}).

%\noindent\emph{Related works}. 
%For \Cref{submod-G-F}, the case where $\F=\Q$ is due to Eskin–Mozes–Shah \cite[Corollary 1.4]{EMS97},  who deduced submodularity from non-divergence estimates for translates of algebraic measures in homogeneous spaces. An independent proof of \Cref{submod-G-F} has recently appeared in Lin~\cite[Theorem 1.7]{Lin25b}; introducing  a notion of stratified set, and arguing inductively on taking sums and intersections, see \cite[Sections 10 \& 11]{Lin25b}.
%We have been informed by Ruotao Yang of yet another proof of the case $\F = \C$, which uses  Weyl's unitary trick and an averaging argument. Neither \cite{EMS97, Lin25b} give any description of the equality case.
%We obtain \Cref{submod-G-F} via a  different approach, which has the double advantage of being more direct and to yield meaningful information regarding the equality case, see  \Cref{submod-G-F-2}.
%The readers may notice a resonance with Petridis' proof~\cite{Petridis}\footnote{See also Terence Tao's comment in this blog post~\url{https://gowers.wordpress.com/2011/02/10/a-new-way-of-proving-sumset-estimates} by Timothy Gowers.} of the Plünnecke-Ruzsa inequality,  and the theory of Harder-Narasimhan filtrations (see e.g.~\cite[Lemma 2]{BS}).
\bigskip

We now engage in the proof of Theorems \ref{submod-G-F}, \ref{submod-G-F-2}.
We may assume without loss of generality that $G$ is closed for the Zariski-topology on $\GL(\cH)$.
Indeed, let $\overline{G}^{\Zar}$ denote the Zariski-closure of $G$ in $\GL(\cH)$. 
Given $W\in \Gr(\cH)$ and $l\in \N$, the set $\{g\in \GL_{N}(\F)\,:\, \dim(gV\cap W)\geq l\}$ is Zariski-closed, therefore $\min_{g\in G}\dim(gV\cap W)=\min_{g\in \overline{G}^{\Zar}}\dim(gV\cap W)$, which allows the reduction.

From now on, we suppose $G=\overline{G}^{\Zar}$. We set
$$\deg_{W}(V) :=\min_{g\in G}\dim(gV\cap W).$$
Observe the map $\deg_W$ is $G$-invariant. The next lemma lists a few additional properties.
\begin{lemma}
\label[lemma]{basics-degW}
Let  $V, V_{1},V_{2}, W\in \Gr(\cH)$.
\begin{itemize}
\item[i)]  There exists a Zariski-dense open subset $O\subseteq G$ such that for all $g\in O$, 
$$\deg_{W}(V) =\dim(gV\cap W).$$
\item[ii)]  The function $\deg_{W}$ is supermodular:
$$\deg_{W}(V_{1})+\deg_{W}(V_{2})\leq \deg_{W}(V_{1}\cap V_{2})+\deg_{W}(V_{1}+V_{2}).$$
\end{itemize}
\end{lemma}

\begin{proof}
For item i), observe that the set $\set{g\in G\,:\, \dim(gV\cap W) \geq \deg_{W}(g) + 1}$ is Zariski-closed.
Its complement, being $\set{g\in G\,:\, \dim(gV\cap W) = \deg_{W}(g)}$, is nonempty and a Zariski-open subset of $G$.
The latter is therefore Zariski-dense because connected algebraic groups are irreducible \cite[Chapter I, \S1.2]{Borel91}.

We deal with item ii). For generic $g$, we have, using item i), then the Grassmann formula, then item i) again
\begin{align*}
\deg_{W}(V_{1})+\deg_{W}(V_{2})
&=\dim(gV_{1}\cap W)+\dim(gV_{2}\cap W)\\
&=\dim(gV_{1}\cap gV_{2} \cap W)+\dim(gV_{1}\cap W+gV_{2}\cap W)\\
&\leq \dim(g(V_{1}\cap V_{2})\cap W)+\dim(g(V_{1} +V_{2})\cap W)\\
&= \deg_{W}(V_{1}\cap V_{2})+\deg_{W}(V_{1}+V_{2}). \qedhere
\end{align*}
\end{proof}

For the next lemma, we write
 $$r_{W}:=\max \left\{\, \frac{\deg_{W}(V)}{\dim V} \,:\, V\in \Gr(\cH)\smallsetminus \{0\} \,\right\}.$$
 
\begin{lemma}  \label[lemma]{eq-case-subl}
The set 
$$\cL=\left\{\, V\in \Gr(\cH)\,:\, \deg_{W}(V) = r_{W} \dim V \, \right\} $$
is a $G$-invariant sublattice. In particular, $$r_{W}= \frac{\dim W}{\dim \cH}.$$
\end{lemma} 

\begin{proof}
Note $\cL$ is $G$-invariant because $\deg_{W}$ is $G$-invariant.

We show that $\cL$ is stable under taking intersection and sum.
Let $V_{1},V_{2}\in \cL$. 
Using the the definition of $r_W$, then \Cref{basics-degW} ii), then the definition of $r_{W}$ again, we have
\begin{align*}
r_{W}\dim(V_{1})+r_{W}\dim(V_{2})
&= \deg_{W}(V_{1})+\deg_{W}(V_{2})\\
&\leq \deg_{W}(V_{1}\cap V_{2})+\deg_{W}(V_{1}+V_{2})\\
&\leq  r_{W}\dim(V_{1}\cap V_{2})+ r_{W}\dim_{W}(V_{1}+V_{2}).
\end{align*}
By the Grassmann formula, we know there is equality between the left-hand side and the final sum on the right-hand side, whence all inequalities involved are in fact equalities. In particular, $\deg_{W}(V_{1}\cap V_{2})=r_{W}\dim(V_{1}\cap V_{2})$ and similarly for $V_{1}+V_{2}$. This justifies that $V_{1}\cap V_{2}, V_{1}+ V_{2}\in \cL$.

We also need to check $\cL$ contains $\{0\}$ and $\cH$. The first is obvious. The second comes from the stability of $\cL$ by sum, the $G$-invariance of $\cL$, and the irreducibility of $G\acts \cH$.

Finally, the formula for $r_{W}$ follows from the relation $\cH\in \cL$ established at the moment.
\end{proof}

\begin{remark}
Note that the irreducibility assumption is only used at the end of this proof.
Without it, we can still conclude that there is a unique subspace $V_W \in \cH \setminus\{0\}$ satisfying $r_W = \frac{\deg_W(V_W)}{\dim V_W}$ and containing all other such subspaces.
It is moreover $G$-invariant.
When we do have irreducibility, $V_W$ must be $\cH$. In this case, in the language of Harder-Narasimhan theory, we say $\cH$ is semistable.
\end{remark}

We are now able to establish the submodular inequality announced in \Cref{submod-G-F}.

\begin{proof}[Proof of \Cref{submod-G-F}]
This is a direct consequence of the formula for $r_{W}$ in \Cref{eq-case-subl}.
\end{proof}

To derive \Cref{submod-G-F-2} about the equality case, we need a preliminary observation.
\begin{lemma}    \label[lemma]{split-lattice}
Every $G$-invariant sublattice of $\Gr(\cH)$ is complemented.
\end{lemma}

\begin{proof}
Call $\cL$ such a lattice. Let $L\in \cL$ with non-zero minimal dimension. Let $V\in \cL$. For every $g\in G$, the intersection $gL\cap V$ belongs to $\cL$, therefore coincides with $gL$ or is $\{0\}$ by minimality of $\dim L$. If $V\neq \cH$, by irreducibility of $G\acts \cH$, we may therefore find $g$ such that $V$ and $gL$ are in direct sum. Repeating the argument with $V\oplus gL$ in the place of $V$, we obtain the desired complement.
 %$g_{1}, \dots, g_{k}\in G$ ($k\geq0$) such that $\cH=V\oplus \bigoplus_{i=1}^k g_{i}L$. 
\end{proof}

We may now conclude regarding the equality case.
\begin{proof}[Proof of \Cref{submod-G-F-2}]

We see from \Cref{eq-case-subl} that for all $V'\in \cL_{V}$, all $g\in G$, we have 
$$\dim(gV' \cap W)\geq \deg_W(V') = r_W \dim V' = \frac{\dim V' \dim W}{\dim \cH}.$$
Noting the left-hand side is also $\dim(V' \cap g^{-1}W)$, then applying  \Cref{eq-case-subl} to $W$ for every fixed $V'$, we deduce that for all $W'\in \cL_{W}$, 
\begin{equation}\label{V'W'ineq}
\dim(V' \cap W')\geq  \frac{\dim V' \dim W'}{\dim \cH}.
\end{equation}
By \Cref{split-lattice},  every $W'\in \cL_{W}$ admits a complement  $W''\in \cL_{W}$. Applying Equation \eqref{V'W'ineq}  to  both $W'$ and $W''$, we obtain
\begin{align*}
\dim V' &\geq \dim(V' \cap W')+\dim(V' \cap W'')\\
&\geq \frac{\dim V' \dim W'}{\dim \cH}+\frac{\dim V' \dim W''}{\dim \cH}=\dim V'.
\end{align*}
Inequalities in this chain must therefore be equalities. This forces equality in \eqref{V'W'ineq} and the splitting
\[
V'=(V' \cap W') \oplus (V' \cap W'').\qedhere
\]
\end{proof}

\section{Dimension interpolation under random walks}\label{Sec-dim-pres}
In this section, we establish dimensional interpolation properties for the action of a Zariski-dense random walk on a simple homogeneous space. The main result is \Cref{RW-preserve-dim}. It implies \Cref{pr:dim-quasi-pres}, and thus validates the first of the two key steps toward the main results of the paper (see \Cref{reductions}).

%In this section, we establish \Cref{pr:dim-quasi-pres}, essentially stating that a random walk on a homogeneous space arising from a simple Lie group almost preserves the dimensional properties of its initial distribution. \Cref{pr:dim-quasi-pres} is a consequence of \Cref{RW-preserve-dim} below, which allows for the action of a single random element rather than the

\bigskip

Let $G$ be a non-compact connected real Lie group with finite center and simple Lie algebra $\kg$. Fix a maximal compact subgroup  $K\subseteq G$, write $\kk\subseteq \kg$  its Lie algebra, and $\ks$ the orthogonal of $\kk$ in $\kg$ for the Killing form. Fix a Cartan subspace  $\ka\subseteq \ks$. Write $\Phi\subseteq \ka^*\smallsetminus\{0\}$ the associated restricted root system, fix a choice of positive roots $\Phi^+\subseteq \Phi$. We write $\ka^+$ the corresponding Weyl chamber, $\ka^{++}$ its interior, and $d=\dim G$. We endow $\kg$ with the scalar product $-\Kill\bigl( \cdot,\vartheta(\cdot)\bigr)$ where $\Kill$ is the Killing form, and $\vartheta$ is the Cartan involution associated to $K$, namely $\vartheta=\Id_{\kk}\oplus -\Id_{\ks}$. We write $\norm{\cdot}$ the associated Euclidean norm on $\kg$. Note that $\Ad(K)$ preserves $\norm{\cdot}$, and $\ad(\ka)$ consists of self-adjoint endomorphisms. We endow $G$ with the induced right $G$-invariant Riemannian metric.

%Let $G$ be a non-compact connected real Lie group with finite center and simple Lie algebra. Suppose $G$ endowed with a maximal compact subgroup $K$ and a compatible Cartan subspace $\ka\subseteq \kg$. Write $\Phi\subseteq \ka^*\smallsetminus\{0\}$ the associated restricted root system, fix a choice of positive roots $\Phi^+\subseteq \Phi$. We write $\ka^+$ the corresponding Weyl chamber, $\ka^{++}$ its interior, and $d=\dim G$. We 
%\cite[Lemma 6.33]{BQ_book}, we may equip $\kg$ with a Euclidean norm $\norm{\cdot}$ such that $\Ad(K)\subseteq O(\kg)$, $\ad(\ka)\subseteq \Sym(\kg)$. We consider the induced right-invariant Riemannian metric on $G$.

Let $\Lambda\subseteq G$ be a lattice.
Equip $X=G/\Lambda$ with the quotient metric.

\bigskip Below, the geometric data $G$, $K$, $\ka$, $\Phi^+$, $\Lambda$ will be considered as fixed, and we will occasionally use the notation $\data$ to refer to this setting.

%in the notations $\ll, \lll, O(.), \simeq$ will use the subscript $\clubsuit$ to indicate that the implicit constant depends on the geometric $G$, $K$, $\ka$, $\Phi^+$, $\norm{\cdot}$, $\Lambda$. Additional dependences will be indicated as subscript.
%refer implicitely to constants possibly depending on the above geometric data $G$, $K$, $\ka$, $\Phi^+$, $\norm{\cdot}$, $\Lambda$. Additional dependences will be indicated as subscript.
\bigskip

Let $\mu$ be a Zariski-dense probability measure on $G$ with finite exponential moment.
We write $\kappa_{\mu}\in \ka^{++}$ its Lyapunov vector \cite[Section 10.4]{BQ_book}, and set $\lambda_{1}> \dots > \lambda_{m+1}$ the collection of the eigenvalues of $\ad(\kappa_{\mu}) \in \End(\kg)$ ordered by decreasing order.
Let $j_{i}\geq1$ denote the multiplicity of $\lambda_{i}$.

 The next theorem considers a measure $\nu$ on $X$ and a small scale $\delta>0$. As a first approximation, it guarantees that if $\nu$ has normalized dimension at least $\alpha$ at scales above $\delta^2$, then for $n\simeq \frac{1}{4\lambda_{1}}\abs{\log \delta}$, and most $g\in G$ selected by $\mu^n$, it has dimension at least $\alpha-\eps$ with respect to the sets $(gB^G_{\delta}x)_{x\in X}$. It is in fact a bit more general as the only scales that matter for the dimensional assumption are those occuring as side lengths of $gB^G_{\delta}x$. Also $\alpha$ does not need to be uniform among those scales, and the output dimensional estimate interpolates the dimensional inputs. As we saw in \Cref{reductions} (via the use of \Cref{pr:preservenu}), this flexibility is crucial for performing the bootstrap.

\begin{thm} \label{RW-preserve-dim}
Let $s, \eps_{1}, \eps_{2}, \delta \in (0, 1)$. Let $\nu$ be a Borel measure on $X$ of mass at most $1$ and which is supported on $\{\inj \geq \delta^{2/3}\}$. For $i=1,\dots, m+1$, let $t_{i}>0$ such that
$$\sup_{x\in X}\nu(B^G_{\delta^{1-s\lambda_{i}}}x)\leq t_{i}.$$

Assume $s\leq \frac{1}{4\lambda_{1}}$ and $\eps_{2}, \delta\lll_{\data, \mu,s, \eps_{1}} 1$. 
Set $n = \lfloor s \abs{\log \delta} \rfloor$. Then there exists $E\subseteq G$ with $\mu^n(E)\leq \delta^{\eps_{2}}$ such that for every $g\in G\smallsetminus E$, for some $F_{g}\subseteq X$ satisfying $\nu(F_{g})\leq \delta^{\eps_{2}}$, we have
$$\sup_{x\in X} \nu_{|X\smallsetminus F_{g}} (gB^G_{\delta}x)\leq \delta^{-\eps_{1}}\prod_{i}t^{j_{i}/d}_{i}.$$
\end{thm}

We will deduce \Cref{RW-preserve-dim} from our subcritical multislicing estimate \Cref{subc-mult}. For this estimate to apply, we need suitable linearizing charts in which the translates of balls by an element $g$ look like boxes carried by a partial flag. 
Those charts are constructed in \S\ref{Sec-covering charts}, and the boxes are described in \S\ref{Sec-random-boxes}.
We also need the subspaces involved in the partial flag to satisfy a subcritical projection property as $g$ varies according to $\mu^n$. Non-concentration properties for this random flag are studied in \S\ref{Sec-nc-random-box}. The analysis is based on our submodular inequality from \Cref{Sec-submod}. Combined with \Cref{subcritical-projection}, we obtain the relevant subcritical projection property. The proof \Cref{RW-preserve-dim} is then concluded in \S\ref{Sec-proof-dimstab-conclusion}.

%We will deduce \Cref{RW-preserve-dim} from our subcritical multislicing estimate \Cref{subc-mult}. For this estimate to apply, we need suitable linearization charts in which the translates of balls by an element $g$ look like boxes carried by a partial flag. \marginpar{These charts are a new technique? Compared to LMWY?}
%We also need the subspaces involved in the partial flag to satisfy a subcritical projection property as $g$ varies according to $\mu^n$. We start by introducing those charts. We then investigate this subcritical property. Using \Cref{subcritical-projection}, we show the latter reduces to a submodular inequality for Borel invariant subspaces in simple complex Lie algebras. This inequality will be established in \Cref{Sec-submod}, thus completing the proof of \Cref{RW-preserve-dim}.

\subsection{A covering of linearizing charts} \label{Sec-covering charts}

We cover $X$ by local exponential charts at a small scale $r>0$. We show that those charts linearize into Euclidean boxes the translates of balls that are not too distorted, namely the subsets $(gB^G_{\rho}y)_{g\in G, \rho>0, y\in X}$ for which $B^\kg_{r^2}\subseteq \Ad(g)B^\kg_{\rho} \subseteq B^\kg_{r}$.

Given $x\in X$,
we recall the injectivity radius of $X$ at $x$ is given by
$$\inj(x):=\sup\{r>0\,:\, B^G_{r}\rightarrow X, g\mapsto gx \text{ is injective }\}.$$
As $B^G_{r}$ denotes an \emph{open} ball, the above supremum is in fact a maximum. We also let $c_{0}=c_{0}(G, \norm{\cdot})>0$ be the largest\footnote{This maximality condition on $c_{0}$ will not be used, it is merely a way to define $c_{0}$ canonically in terms of $G, \norm{\cdot}$.} constant such that $\exp: B^\kg_{c_{0}}\rightarrow G$ is injective and we set
$$\exp_{x}:B^\kg_{\inj(x)\wedge c_{0}} \rightarrow X, \,v\mapsto \exp(v)x$$
where for any $a,b>0$, we use the notation $a\wedge b=\min(a,b)$.
Noting that for any $r>0$, we have $\exp(B^\kg_{r})\subseteq B^G_{r}$, we see the map $\exp_{x}$ is injective.

\begin{lemma}
\label[lemma]{local-linear}
Let $x\in X$, let $0<r\lll_{G}\inj(x) \wedge 1$.
Let $g\in G$, $\rho>0$, $y\in X$, such that $gB^G_{\rho}y\cap B^G_{r}x\neq \emptyset$ and $B^\kg_{r^2}\subseteq \Ad(g)B^\kg_{\rho} \subseteq B^\kg_{r}$.
Then $\exp_{x}^{-1}\bigl(gB^G_{\rho}y \bigr)$ is covered by $O_{G}(1)$ many translates of $\Ad(g)B^\kg_{\rho}$.
\end{lemma}

\begin{remark}
 The exponential map does not linearize translates of balls which are too asymetric, this is why we require the condition $B^\kg_{r^2}\subseteq \Ad(g)B^\kg_{\rho} \subseteq B^\kg_{r}$.
\end{remark}

\begin{remark} There is no dependence on the norm $\norm{\cdot}$ in \Cref{local-linear}. This is because any other norm $\norm{\cdot}'$ on $\kg$ that arises from a maximal compact subgroup $K'$ of $G$ satisfies $C^{-1}\norm{\cdot}\leq \norm{\cdot}'\leq C\norm{\cdot}$ for some $C=C(G)$ independent\footnote{Indeed, the set of pairs $(\kk, \ks)\in \Gr(\kg)^2$ where the Killing form $\Kill$ is  negative definite on $\kk$, positive definite on $\ks$, and $\kk$ is the orthogonal of $\ks$ for $\Kill$, is compact. } of $\norm{\cdot}'$.
\end{remark}

\begin{proof}
Note the assumption $\Ad(g)B^\kg_{\rho} \subseteq B^\kg_{r}$ implies $\rho\leq r$. Combined with $r\lll_{G}1$, we have both $B^G_r \subset \exp(B^\kg_{2r})$ and $B^G_{2\rho} \subset \exp(B^\kg_{4\rho})$.
Since $g B^G_\rho y \cap B^G_{r}x \neq \emptyset$, there is $w_0 \in B^\kg_{2 r}$ such that
\(
\exp(w_0) x \in  g B^G_\rho y,
\)
or equivalently,
\(
y \in B^G_\rho g^{-1} \exp(w_0) x.
\)

Let $v\in B^\kg_{\inj(x)\wedge c_{0}}$ such that $\exp_{x}(v) \in g B^G_\rho y$.
We have
\begin{align*}
\exp_x(v) &\in g B^G_{2\rho} g^{-1} \exp(w_0) x \\
& \subset g \exp( B^\kg_{4\rho} ) g^{-1} \exp(w_0) x \\
& = \exp\bigl( \Ad(g) B^\kg_{4\rho} \bigr) \exp(w_0) x.
\end{align*}
In other words, there is a vector 
\(
w \in \Ad(g) B^\kg_{4\rho}
\)
such that 
\[
\exp_x(v) = \exp(w) \exp(w_0) x.
\]
By the assumption $\Ad(g)B^\kg_{\rho} \subseteq B^\kg_{r}$, we have $\norm{w} \leq 4r$, we derive from the Baker-Campbell-Hausdorff formula that 
$$\exp(w)\exp(w_0) = \exp( w+w_0 + O_{G}(r^2)).$$
Since $r \lll_{G} \inj(x) \wedge 1$, the vector on the right-hand side is in $B^\kg_{\inj(x)\wedge c_0}$.
The injectivity of $\exp_x$ then implies
\[
v = w + w_0 + O_{G}(r^2).
\]
This justifies
$$\exp_{x}^{-1} (g B^G_\rho y) \subseteq \Ad(g)B^\kg_{4\rho} + w_0 +B^\kg_{O_{G}(r^2)}.$$
Using the assumption $B^\kg_{r^2}\subseteq \Ad(g)B^\kg_{\rho}$, we see that the set on the right-hand side is covered by $O_{G}(1)$-many translates of $\Ad(g)B^\kg_{\rho}$. This finishes the proof.
\end{proof}

Patching together charts from the previous lemma, we deduce the following. It allows to convert \Cref{RW-preserve-dim} into a linear statement.

%The next lemma allows to see a measure $\nu$ on $X$ through a net of local charts in order to straighten simultaneously images of balls by elements of $g$ into Euclidean boxes of $\kg$.

\begin{lemma}
\label[lemma]{net-charts}
Let $0<r \lll_{G}1$. There exists a measurable map $\varphi : \{\inj \geq r\}\rightarrow B^\kg_{1}$ satisfying the following.
\begin{itemize}
\item[1)] For every $\rho\in (0, r)$, $v\in \kg$, the preimage $\varphi^{-1}(B^\kg_{\rho}+ v)$ is covered by $O_{\data}(1)$ many balls of the form $(B^G_{\rho}x)_{x\in X}$
\item[2)] For every $\rho\in (0, r)$, $g\in G$ such that $B^\kg_{r^2}\subseteq \Ad(g)B^\kg_{\rho}\subseteq B^\kg_{r}$, and $x\in X$, the translate $gB^G_{\rho}x \cap \{\inj \geq r\}$ is covered by $O_{\data}(1)$ many preimages of boxes of the form $(\varphi^{-1}(\Ad(g)B^\kg_{\rho}+ v))_{v\in \kg}$.
\end{itemize}
\end{lemma}

In particular, given a measure $\nu$ on $X$ supported on $\{\inj \geq r\}$, we see that the $\varphi_{\star}\nu$-measure of balls on $\kg$ is controlled by the $\nu$-measure of balls on $X$ (up to radius $r$), while the $\nu$-measure of translates $gB^G_{\rho}x$ is controlled by the $\varphi_{\star}\nu$-measure of boxes $\Ad(g)B^\kg_{\rho}+ v$ provided the size of the box belongs to a certain window prescribed by $r$.

\begin{remark} The map $\varphi$ depends on $r$.  In practice, the parameter $r$ will be a  power of $\delta$, with exponent  macroscopic and smaller than $1$, e.g. $\delta^{2/3}$ in the proof of \Cref{RW-preserve-dim}. We note that the  radius $\rho$ appearing in item 1)  is required to be smaller than $r$. It would be possible to refine the construction of $\varphi$ in order to allow $\rho$ bigger than $r$ in item 1), say $\rho \in [r, \eta]$ where $\eta>r$  satisfies $\supp \nu\subseteq \{\inj \geq \eta\}$.  We stick to the above version for simplicity.
Finally, we note that in item 2), the condition on $g$ forces $\rho \in [r^2, r]$, in particular $\rho$ is not arbitrarily small in item 2). 
\end{remark}

\begin{proof}

We let $C>1$ be a parameter to be specified later depending on $G$.
Let $\{x_j\}_{j \in \cJ}$ be a maximal $r/C$-separated set of points in $\{\inj\geq r\}$.
Then since the balls of radius $r/(2C)$ centered in $\{x_j\}_{j \in \cJ}$ are disjoint, we have $\abs{\cJ} \ll_\data (r/C)^{-d}$.
For each $j \in \cJ$, let $U_j = B^G_{r/C} x_j$.
By maximality, we have $\{\inj\geq r\} \subseteq \cup_{j}U_{j}$.
By the triangle inequality, we have $\cup_{j}U_{j} \subset \{\inj \geq r/2 \}$ provided $C\geq 2$.
Taking  $C \ggg_{G}1$ large enough,  we may assume that the map $\exp_{x_{j}}^{-1}$ defines a $2$-bi-Lipschtiz diffeomorphism from $U_{j}$ to an open subset $V_{j}'\subseteq B^\kg_{2r/C}$.
One may compose by similarities to make those $V_{j}'$'s disjoint in $B^{\kg}_{1}$.
More precisely, one may choose  $s=s(\data)>0$ (small), some vectors $v_{j}\in \kg$, such that writing $\tau_{j}=s\Id_{\kg}+v_{j}$ and $V_{j}=\tau_{j}(V'_{j})$, the sets $(V_{j})_{j\in \cJ}$ are included in  mutually disjoint balls of radius $2rs/C$ in $B^\kg_{1}$.
Let $\varphi_j = \tau_{j}\circ {\exp^{-1}_{x_{j}}}_{|U_{j}}: U_j \to V_j$ denote the resulting diffeomorphisms.
Then define $\varphi : \{\inj \geq r\}\rightarrow B^\kg_{1}$ to be a measurable map coinciding with one of the $(\varphi_{j})_{j}$ at every point, i.e., such that for all $x\in \{\inj \geq r\}$ we have $\varphi(x)\in \{\varphi_{j}(x)\,:\,j\in \cJ\}$.

We check that $\varphi$ satisfies  item 1). Let $\rho\in (0, r)$,  $v\in \kg$. The separation condition on the $(V_{j})_{j\in \cJ}$ implies that   $\cJ':=\{j\,:\, V_{j}\cap (B^\kg_{\rho}+ v)\neq \emptyset\}$ has cardinality $|\cJ'|=O_{\data}(1)$. Moreover, for $j\in \cJ'$,  the preimage $\varphi^{-1}_{j}((B^\kg_{\rho}+ v) \cap V_{j})$ has diameter $O_{\data}(\rho)$, so it is covered by $O_{\data}(1)$ $\rho$-balls in $G$. Hence item 1). 

For item 2), note it is sufficient to establish the claim with $\rho_{1}=\rho/C$ instead of $\rho$.
The assumption $\Ad(g)B^\kg_{\rho}\subseteq B^\kg_{r}$ implies that $gB^G_{\rho_{1}}x$ has diameter $O(r/C)$.
It follows from the separation condition on the $(U_{j})_{j\in \cJ}$ that $\cJ'':=\{j\,:\, U_{j}\cap gB^G_{\rho_{1}}x \neq \emptyset\}$ has cardinality $|\cJ''|=O_{\data}(1)$.
Assuming $C$ large enough (depending on $G$ again),  we can apply \Cref{local-linear} to guarantee that for each $j\in \cJ''$, the set $\varphi_{j}( U_{j} \cap gB_{\rho_{1}}x)$ is included in $O_{\data}(1)$ translates of $\Ad(g)B^\kg_{\rho_{1}}$. Hence item 2). 
%We now check item 2). Consider $\rho,g, x$ as in item 2). The assumption $\Ad(g)B^\kg_{\rho}\subseteq B^\kg_{r}$ implies that $g B^G_{\rho}x$ is included in a ball of radius $O_{d}(r)$, in particular $g B^G_{\rho}x\cap \{\inj \geq r\}$ is covered by $O_{d}(1)$ sets $U_{j}$. Hence there exists $j$ such that  $\nu(gB_{\rho}x)\ll_{d}  \nu(gB_{\rho}x \cap U_{j})=\tnu(\varphi_{j}(gB_{\rho}x \cap U_{j}))$. The assumption $B^\kg_{r^2}\subseteq \Ad(g)B^\kg_{\rho}\subseteq B^\kg_{r}$ allows to apply \Cref{local-linear} to guarantee that $\varphi_{j}(gB_{\rho}x \cap U_{j})$ is included in $O_{d}(1)$ translates of $\Ad(g)B_{\rho}$. Item 2) follows.
\end{proof}

\subsection{The random boxes in the Lie algebra} \label{Sec-random-boxes}

Given a random parameter $g\sim \mu^n$, we describe the box $\Ad(g)B^\kg_{1}$.

%We also see that the random filtration of subspaces of $\kg$ induced by $\Ad(g)B^\kg_{1}$ as $g\sim \mu^n$ satisfies some angle control with respect to any given subpace of $\kg$ unless a trivial obstruction occurs.
\bigskip

Every $g\in G$ admits a Cartan decomposition
\begin{equation} \label{def-Cartan-decomp}
g=\theta_{g} a_{g} \theta'_{g}
\end{equation}
where $\theta_{g}, \theta'_{g}\in K$ and $a_{g}=\exp(\kappa(g))$ with $\kappa(g) \in \ka^+$. The element $\kappa(g)$ is uniquely determined by $g$ and called the Cartan projection of $g$. The components $(\theta_{g}, \theta'_{g})$ are not uniquely defined, we choose them to depend measurably on $g$.

Set $\kappa_{\mu}=\lim_{n \to +\infty} n^{-1}\int_G \kappa(g)  \dd\mu^n(g)$ to be the Lyapunov vector of $\mu$.
It is known that $\kappa_{\mu}$ is well defined and belongs to $\ka^{++}$, see \cite[Theorem 10.9]{BQ_book}.
For $\alpha\in \ka^*$, set $\kg_{\alpha}:=\{v\in \kg\,:\, \forall w\in \ka, \,[w,v]=\alpha(w)v\}$, so that $\kg=\oplus_{\alpha\in \Phi\cup \{0\}}\kg_{\alpha}$ is the restricted root space decomposition associated to our choice of Cartan subspace $\ka$.
Enumerate $\{\alpha(\kappa_{\mu}) \,:\, \alpha\in \Phi\cup \{0\}\}=\{\lambda_{1}>\lambda_{2}>\dots>\lambda_{m+1}\}$, set for $i=1, \dots, m+1$,
\begin{equation}\label{Cartan-filtration}
V_{i}:=\bigoplus \set{\kg_{\alpha}\,:\,\alpha(\kappa_{\mu})\geq \lambda_{i}}.
\end{equation}
In particular, $V_{m + 1} = \kg$.

The next lemma states that for $g\sim \mu^n$, the set $\Ad(g)(B^{\kg}_{1})$ is essentially a Euclidean box with associated partial flag $(\Ad(\theta_{g})V_{i})_{i}$ and size parameters $(e^{n\lambda_{i}})_{i}$. 
\begin{lemma}
\label[lemma]{descrip-AdgB1}
Given $\eps>0$, there exists $\eta=\eta(\mu, \eps)>0$ such that for $n\ggg_{\mu,\eps}1$, for $g\in G$ oustide of a set of $\mu^n$-measure at most $e^{-\eta n}$, we have
\begin{equation*}
\Ad(g) B^{\kg}_1 \subseteq \Ad(\theta_{g})\bigl(B^{V_{1}}_{e^{n(\lambda_{1}+\eps)}}+\dots+ B^{V_{m+1}}_{e^{n(\lambda_{m+1}+\eps) }}\bigr),
\end{equation*}
while the converse inclusion holds provided $\eps$ is replaced by $-\eps$.
\end{lemma}

\begin{proof}
For every $g\in G$, we have
$$\Ad(g) B^{\kg}_1 = \Ad(\theta_{g})( \Ad(a_{g})B^{\kg}_{1}) \subset \Ad(\theta_{g}) \Bigl(\sum_{\alpha\in \Phi\cup \{0\}} B^{\kg_{\alpha}}_{e^{\alpha(\kappa(g))}} \Bigr).$$
Given $\eps>0$, the large deviation principle for the Cartan projection \cite[Theorem 13.17]{BQ_book} yields some $\eta=\eta(\mu,\eps)>0$ such that for $n\ggg_{\mu,\eps}1$, 
$$\mu^n\setbig{g\,:\,\max_{\alpha\in \Phi\cup \{0\}} \abs{\alpha(\kappa(g)-n\kappa_{\mu})} \leq 2^{-1} \eps n } \geq 1-e^{-\eta n}.$$
For $g$ in the above set, we deduce
\begin{align*}
\Ad(g)(B^{\kg}_{1}) & \subseteq \Ad(\theta_{g}) \Bigl(\sum_{\alpha\in \Phi\cup \{0\}} B^{\kg_{\alpha}}_{e^{n(\alpha(\kappa_{\mu})+\eps/2) }} \Bigr) \\ %\label{Adg=box2}\\
&\subseteq \Ad(\theta_{g})\Bigl(B^{V_{1}}_{e^{n(\lambda_{1}+\eps)}}+\dots+ B^{V_{m+1}}_{e^{n(\lambda_{m+1}+\eps)}}\Bigr).
\end{align*}
The converse inclusion is similar.
\end{proof}

%And we need to check that the largest slice $(\Ad(\theta_{g})V_{1})_{g\sim \mu^{*n}}$ satisfies the strong non-concentration property $\sP_{0}$ with scale parameter $e^{-n}$, while for $i\geq 2$, $(\Ad(g)V_{i})_{g\sim \mu^{*n}}$ satisfies $\sP_{1}$, $\sP_{2}$, or $\sP_{3}$ with scale parameter $e^{-n}$. For this, it is useful to replace $\Ad(\theta_{g})$ by $\Ad(g)$.

\subsection{Non-concentration for the random boxes}\label{Sec-nc-random-box}

We establish two non-concentration properties for the partial flag $(\Ad(\theta_{g})V_{i})_{i=1}^{m+1}$ associated to the random box $\Ad(g)B^\kg_{1}$ where $g\sim \mu^n$.

\bigskip
The first property, \Cref{non-conc-Ad-thetagVi}, states that the random subspace $(\Ad(\theta_{g})V_{i})_{g\sim \mu^n}$ is typically transverse to any prescribed subspace $W\subseteq \kg$, unless $W$ intersects every subspace of the orbit $\Ad(G)V_{i}$, in which case the statement clearly fails. \Cref{non-conc-Ad-thetagVi} also allows for a H\"older-regular control of the angle. This result will be used on many occasions in the rest of the paper.

\begin{proposition}[Angle control]
\label[proposition]{non-conc-Ad-thetagVi}
Let $i\in \{1, \dots m\}$, $W\in \Gr(\kg)$, and $\eps\in (0,1)$, such that $\sup_{g\in G}\dist_{\angle}(\Ad(g)V_{i}, \, W)>\eps$. There exist $C=C(\data, \mu)>1$ and $c=c(\mu)>0$ such that for $n\geq 1$, $\rho\geq e^{-n}$,
$$\mu^n\{ g\,:\,\dist_{\angle}(\Ad(\theta_{g})V_{i}, \, W)\leq \rho\} \leq C\eps^{-c}\rho^c.$$
\end{proposition}

\begin{proof}
The adjoint action $G \acts \kg$ induces an action $G\acts \bigwedge^{\dim V_{i}}\kg$. By definition of $V_{i}$, the endomorphism $\exp(\kappa_{\mu}) \acts \bigwedge^{\dim V_{i}}\kg$ has a unique dominant eigenvalue, which is simple, with corresponding eigenspace $\bigwedge^{\dim V_{i}}V_{i}:=L_{i}$. Writing $\bigwedge^{\dim V_{i}}\kg$ as a sum of irreducible subrepresentations $\bigwedge^{\dim V_{i}}\kg=\oplus_{k=1}^q E_{k}$ and letting act $\exp(\kappa_{\mu})$, one sees that $L_{i}$ has to be included in some $E_{k_{0}}$ where $k_{0}\in \{1, \dots, q\}$. Write  $E:=E_{k_{0}}$ for short. Note the irreducible subrepresentation $E$ is also proximal. 

Let $\nu_{i}$ be the unique $\mu$-stationary measure on the projective space $P(E)$. Note $\nu_{i}$ is supported on $\bigwedge^{\dim V_{i}}\Gr(\kg, \dim V_{i})$ because the latter is compact and $G$-invariant.
By exponential convergence of density points \cite[Corollary 4.18]{BQ2}, the distribution of $(\Ad(\theta_{g})L_{i})_{g\sim \mu^{n}}$ on $P(E)$ converges exponentially fast to $\nu_{i}$ outside of an event of exponentially small measure. More precisely, one may find a pair of two $P(E)$-valued random variables $(\xi_{n}, \xi_{\infty})$ defined on a common probability space, such that $\xi_{n}$ has the same law as $(\Ad(\theta_{g})L_{i})_{g\sim \mu^{n}}$ and $\xi_{\infty}$ has law $\nu_{i}$, and satisfying  
$$ \mathbb{P}\bigl[ \dist(\xi_{n}, \xi_{\infty}) > e^{-c n} \bigr]\ll_{\mu} e^{-c n}$$
where $c=c(\mu)>0$.
As the Pl\"ucker embedding $\Gr(\kg, \dim V_{i}) \rightarrow P(\Lambda^{\dim V_{i}} \kg)$ is bi-Lipschitz, this allows to replace $(\Ad(\theta_{g})V_{i})_{g\sim \mu^{n}}$ by $V\sim \nu_{i}$ in order to establish the proposition.

We may also assume $\dim V_{i}+\dim W=\dim \kg$, for otherwise we may find $W'$ with the right dimension, containing $W$, and such that $\dist_{\angle}(\Ad(g)V_{i}, \, W')>\eps$ for some $g\in G$; and it is sufficient to establish the lemma for $W'$.

Now, letting $\uv$, $\uw$ be wedge products of orthonormal basis of $V$, $W$, setting $\varphi_{\uw}: E\rightarrow \bigwedge^d\kg\simeq \R, u\mapsto u\wedge \uw$, we have
$$\dist_{\angle}(V, \, W)=\| \uv \wedge \uw\|=\|\varphi_{\uw}(\uv)\| = \|\varphi_{\uw}\| \dist_{\angle}(\R \uv, \Ker \varphi_{\uw}).$$
By assumption, we know that $\|\varphi_{\uw}\|>\eps$. The result then follows from the H\"older regularity of the measure $\nu_{i}$ with respect to neighborhoods of hyperplanes, see \cite[Theorem 14.1]{BQ_book}.
\end{proof}

The second property, \Cref{avoid-bad-pencil-proba}, considers an arbitrary subspace $W\subseteq \kg$ and guarantees a partial transversality with $\Ad(\theta_{g})V_{i}$ for most $g\sim \mu^n$. In view of \Cref{subcritical-projection}, this is enough to ensure that the random projector $\pi_{||\Ad(\theta_{g})V_{i}}$ satisfies a subcritical projection theorem at scales above $e^{-n}$.

\begin{proposition}[Strict probabilistic submodularity]
\label[proposition]{avoid-bad-pencil-proba}
 Let $i\in \{1, \dots, m\}$ and $W\subseteq \kg$ a non-zero proper subspace. There exists $c=c(\mu)>0$ such that for all $n\geq 1$, $\rho\geq e^{-n}$
\begin{equation}\label{eq-nc-mungabpp}
\mu^n \left\{g\,:\, \max_{W'\in B_{\rho}(W)}\frac{\dim \Ad(\theta_{g})V_{i}\cap  W'}{\dim W'} < \frac{\dim V_{i}}{\dim \kg} \right\} \geq 1- O_{\data, \mu}( \rho^{c}).
\end{equation}
\end{proposition}

\begin{remark}
In the remainder of the paper, \Cref{avoid-bad-pencil-proba} will only be used via its consequence \Cref{mu-Sub}, see below. 
For the sole purpose of establishing \Cref{mu-Sub}, we only need  a weaker form of \Cref{avoid-bad-pencil-proba}, where the inequality $<$ in \eqref{eq-nc-mungabpp} is replaced by a large inequality $\leq $. This weaker result  can be obtained more directly, by combination of the submodular inequality from \Cref{submod-G-F} and the short argument entitled  ``proof  of \Cref{avoid-bad-pencil-proba}'' presented below. The stronger form displayed in \Cref{avoid-bad-pencil-proba} is obtained using the equality case of the submodular inequality, i.e. \Cref{submod-G-F-2}.  We believe  this more precise estimate can be helfpul in other contexts, and justifies  a slight detour.
\end{remark}

\noindent\emph{Strategy of proof}. To prove \Cref{avoid-bad-pencil-proba}, we first establish a purely geometric version of the statement (with no random variables). It is presented below as \Cref{avoidW'}, and relies Theorems \ref{submod-G-F}, \ref{submod-G-F-2}.
From there, we upgrade the geometric statement to the desired probabilistic result using \Cref{non-conc-Ad-thetagVi}.

\bigskip
\begin{comment}
We start with preliminaries, which will allow us to exploit to the submodular inequality from \Cref{Sec-submod}.
We let $\kg_{\C}=\kg\otimes \C$ be the complexification of $\kg$. Note $\kg_{\C}$ is a semisimple complex Lie algebra. We choose a Cartan subalgebra $\kh_{\C}\subseteq \kg_{\C}$, and write $\Phi_{\C}\subseteq \kh_{\C}^*$ the associated root system, $\kg_{\C}=\oplus_{\Phi_{\C}\cup \{0\}}\kg_{\C, \beta}$ the root space decomposition. We choose a set of positive roots $\Phi^+_{\C}\subseteq \Phi_{\C}$.
We denote by $\kn_{\C}=\oplus_{\beta\in \Phi^+_{\C}}\kg_{\C, \beta}$ the sum of positive root spaces, and set $\kb_{\C}=\kh_{\C}\oplus \kn_{\C}$ the associated Borel Lie algebra. Using negative roots, we define similarly $\kn^-_{\C}$, $\kb^-_{\C}$.

\end{comment}

We start with a geometric version of \Cref{avoid-bad-pencil-proba}.
We set $\ku=\oplus_{\alpha\in \Phi^+}\kg_{\alpha} \subseteq \kg$ the sum of positive root spaces for $\ka \overset{\ad}{\acts} \kg$, and $\kp=\kz_{\kg}(\ka)\oplus \ku$ where $\kz_{\kg}(\ka)$ is the centralizer of $\ka$ in $\kg$. Using negative roots, we define similarly $\ku^-, \kp^-$. 
\begin{lemma}
\label[lemma]{avoidW'}
Let $V\subseteq \kg$ be a non-zero proper subspace that is $\ad(\kp)$-invariant and satisfies either $V\subseteq \ku$ or $ \kp \subseteq V$.
Let $W \subseteq \kg$ be a non-zero proper subspace. Then there exists $g\in G$ such that
\begin{equation} \label{eq-sous-mod-kp-inv}
\frac{ \dim( \Ad(g)V \cap W)}{\dim W}< \frac{\dim V}{\dim \kg}.
\end{equation}
\end{lemma}
We begin the proof with the case  where $V$ is contained  in $\ku$.

\begin{proof}[Proof of \Cref{avoidW'} in the {case $V\subseteq \ku$}]
We argue by contradiction, assuming that \eqref{eq-sous-mod-kp-inv} fails. This implies that for every $W'\in \overline{\Ad(G)W}^{\Zar}$,
\begin{equation*}
\frac{ \dim(V \cap W')}{\dim W}\geq \frac{\dim V}{\dim \kg}.
\end{equation*}
By Theorems \ref{submod-G-F}, \ref{submod-G-F-2}, and the $G$-invariance of $\overline{\Ad(G)W}^{\Zar}$, we have in fact equality:
\begin{equation} \label{eq-sous-mod-kp-inv-fail}
\frac{ \dim(V \cap W')}{\dim W}= \frac{\dim V}{\dim \kg}.
\end{equation}

Note that the connected algebraic group with Lie algebra $\ka + \ku^-$ is an $\R$-split solvable $\R$-group.
Hence by the Borel fixed point theorem (over $\R$) \cite[Proposition 15.2]{Borel91}, the variety $\overline{\Ad(G)W}^{\Zar}$ contains some $\ad(\ka+\ku^-)$-invariant element $W'_{0}$. 
Up to replacing $W$ by $W'_{0}$, we may assume $W$ is $\ad(\ka+\ku^-)$-invariant. On the other hand, $V$ is $\ad(\kp)$-invariant and included in $\ku$ on which $\Ad(e^{V_{1}})$ acts trivially, whence for all  $v_{1}\in V_{1}$, 
$$V \cap \Ad(e^{v_{1}}) W \supseteq \Ad(e^{v_1})(V \cap W) = V \cap  W.$$
This passes to the Zariski-closure: for all $W''\in \overline{\Ad(e^{V_{1}}) W}^{\Zar}$,
$$V \cap W'' \supseteq V \cap  W.$$
As dimensions match due to \eqref{eq-sous-mod-kp-inv-fail}, we have in fact
$$V \cap W'' = V \cap  W.$$
Therefore, setting $S=\bigcup \left\{\,W''\,:\, W''\in \overline{\Ad(e^{V_{1}}) W}^{\Zar} \,\right\} $, we have justified 
\begin{equation}\label{eq-VSVW}
V \cap S = V \cap  W.
\end{equation}

We use \eqref{eq-VSVW} to show that $W$ is too big and get a contradiction.  For that we rely on the next two standard facts. 
(*) For every restricted  root $\alpha\in \Phi$ of $\ka \acts\kg$, for every $v\in \kg_{\alpha}\smallsetminus \{0\}$, there exists $v' \in \kg_{-\alpha}\smallsetminus \{0\}$ such that $(v, [v,w], w)$ is an $\sl_{2}$-triple \cite[Proposition 6.52]{Knapp02}. (**) Every non-zero $\ad(\ka+\ku)$-invariant  subspace of $\kg$ intersects $V_{1}$ non-trivially, and  contains $V_{1}$ if it is $\ad(\kp)$-invariant.\footnote{To justify (**), consider $E\subseteq \kg$ a non-zero $\ad(\ka+\ku)$-invariant subspace. By $\ad(\ka)$-invariance, it is enough to show $E\not\subseteq \oplus_{\alpha\neq \alpha_{\max}}\kg_{\alpha}=:W_{m}$ where $\alpha_{\max}$ is the longest root of $\Phi^+$ (i.e. $\kg_{\alpha_{\max}}=V_{1}$). If $E\subseteq W_{m}$, then by $\ad(\ka+\ku)$-invariance, $\Ad(G)E\subseteq \overline{\exp(\ad(\ku^-\oplus \kz_{\kg}(\ka)))E}\subseteq W_{m}$, contradicting the irreducibility of $\Ad(G)$. Hence the first claim on $E$. The second claim follows because the action of $\ad(\kz_{\kg}(\ka))$ on $V_{1}$ is irreducible (as can be checked by a similar argument).}  
%\marginpar{ref pour irred de $\ad(\kz_{0}(\kg))$ sur $V_{1}$? On l'utilise dans la footnote 13 et aussi plus tard en Section 7}

Getting back to the proof, recalling $W\cap \ku\neq \{0\}$ and $W$ is $\ad(\ka)$-invariant, we have that $W$ must intersect a restricted root space of $\ka \acts \kg$. By $\ad(\ku^-)$-invariance and (*), it follows $W\cap \ka \neq \{0\}$.  Let $\alpha_{\max}\in \ka^*$ be the longest root of $\Phi^+$.

Assume first that $W\cap \ka \not\subseteq \Ker \alpha_{\max}$. Let $x_{0}\in W\cap \ka$ with $\alpha_{\max}(x_{0})=1$. Then for every $v_{1}\in V_{1}$, we have that  $\Ad(e^{v_{1}})W$ contains the vector $\Ad(e^{v_{1}})x_{0}=x_{0}-v_{1}$. % Je pesne que c'est x_0 - v_0, W.H.
Passing to the limit, we deduce  $S\supseteq V_{1}$. 
But $V\supseteq V_{1}$ by (**), so \eqref{eq-VSVW} implies $W\supseteq V_{1}$. As the smallest $\ad(\ka+\ku^-)$-invariant subspace containing $V_{1}$ is $\kg$, we conclude $W=\kg$. Contradiction.

Otherwise, we would have $W \cap \ka \subseteq \Ker \alpha_{\max}$. 
In this case, $W\cap V_{1}=\{0\}$ due to the $\ad(\ku^-)$-invariance of $W$ and (*). Fact (**) (applied dually to $\ka+\ku^-$) also gives $W\cap W_{1}\neq \{0\}$, say $W\cap W_{1}\supseteq \R w_{1}$ with $w_{1}\neq 0$.
Use (*) to complete $w_{1}$ into an  $\sl_{2}$-triple $(v_{1}, [v_{1}, w_{1}], w_{1})$ where $v_{1}\in V_{1}$.
We have $\Ad(e^{t v_{1}})w_{1}=w_{1}+ t\ad(v_{1})w_{1}+ \frac{t^2}{2}\ad(v_{1})^2w_{1}$ where $\ad(v_{1})^2w_{1}= -2v_{1}\neq 0$, and it follows $\Ad(e^{t v_{1}})\R w_{1}$ converges to $\R v_{1}$ as $t\to +\infty$.
Hence $S\cap V_{1}\neq \{0\}$, so $W\cap V_{1}\neq \{0\}$ due to \eqref{eq-VSVW} and $V\supseteq V_{1}$. Contradiction.
\end{proof}

We now reduce the case $\kp\subseteq V$ in \Cref{avoidW'} to the case $V \subseteq \ku$ established above. For that we make use of the \emph{Weyl group of $G$}, identified with $N_{K}(\ka)/Z_{K}(\ka)$ where $N_{K}(\ka), Z_{K}(\ka)$ are respectively the stabilizer and fixator of $\ka$ in $K$ for the adjoint action. We denote by
\begin{equation} \label{def-iota}
\text{$\iota$ the unique element in the Weyl group such that $\iota(\ka^+) = -\ka^{+}$.}
\end{equation}
Alternatively, $\iota$ is the longest element of the Weyl group (for our choice of positive roots $\Phi^+$).
Note that $\iota$ is an involution, and identifying abusively $\iota$ with any representative in $K$, we have $\Ad\iota(\ku)=\ku^-$ and $\Ad\iota(\kp)=\kp^-$.

\begin{proof}[Proof of \Cref{avoidW'} in the {case $\kp\subseteq V$}] We omit the notation $\Ad$ for conciseness.
Recall $\kg$ is endowed with a $K$-invariant scalar product for which elements of $\ad(\ka)$ are self-adjoint.
The second condition implies that the restricted root spaces $(\kg_{\alpha})_{\alpha\in \Phi\cup \{0\}}$ are mutually orthogonal, in particular $\kp^\perp=\ku^-$ and $V^\perp\subseteq \ku^-$.
Acting with the longest element $\iota$ of the Weyl group $N_{K}(\ka)/Z_{K}(\ka)$ (here identified with a representative in $K$), see \eqref{def-iota}, and using that $ K$ preserves the scalar product, we get $(\iota V)^\perp=\iota (V^\perp) \subseteq \ku$. On the other hand, note\footnote{
Indeed, given any $w\in \kg$, the $\ad(w)$ anti-invariance of $\text{Kill}$ implies that the transpose of $\ad(w)$ is $\ad(-\vartheta w)$. It remains to check  that  $\vartheta$ switches $\kp$ and $\kp^{-}$: this is because $\vartheta\in \Aut(\kg)$ is an involution which acts on $\ka$ via $-\Id_{\ka}$, whence sends any restricted root space $\kg_{\alpha}$ to its opposite $\kg_{-\alpha}$.
}
 that $\ad(\kp^{-})$ is the transpose of $\ad(\kp)$ for our choice of Euclidean structure on $\kg$. Therefore $V^\perp$ is $\ad(\kp^-)$-invariant, then $(\iota V)^\perp$ is $\ad(\kp)$-invariant.
Applying the first case, studied previously, to $(\iota V)^\perp$, we deduce that for any subspace $W\subseteq \kg$, there exists $g'\in G$ such that
$$\dim (\iota V)^\perp \dim (\iota W)^\perp > \dim \kg \dim( g'(\iota V)^\perp \cap (\iota W)^\perp). $$
Note from the Iwasawa decomposition and the $\ad(\kp)$-invariance of $(\iota V)^\perp$ that we may ensure that $g'$ is in $K$, in which case $g'(\iota V)^\perp=(g'\iota V)^\perp$.
Applying \Cref{orthog-red-sub}, we finally obtain \eqref{eq-sous-mod-kp-inv} with $g=\iota^{-1} g'\iota$.
\end{proof}

We now use \Cref{non-conc-Ad-thetagVi} to upgrade \Cref{avoidW'} into the desired probabilistic statement.

\begin{proof}[Proof of \Cref{avoid-bad-pencil-proba}]
We claim that there exists a subspace $W_{1}\subseteq W$ such that
\begin{equation}\label{dimW1-vs-W}
\dim W_{1}>  \left(1- \frac{\dim V_{i}}{\dim \kg}\right) \dim W,
\end{equation}
and  satisfying for all $n\geq 1$, $\rho\geq e^{-n}$, and some constant $c=c(\mu)>0$,
\begin{equation}\label{transW1}
\mu^n\setbig{g\,:\, \dist_{\angle}(\Ad(\theta_g)V_{i}, W_{1})\leq \rho } \ll_{\data, \mu} \rho^c. 
\end{equation}
Note that this property is indeed sufficient: if $W'\in B_{\rho}(W)$ and $g\in G$ satisfy that $\dim  W'\cap \Ad(\theta_{g})V_{i} \geq \frac{\dim V_{i}}{\dim \kg}\dim W$, then for dimensional reasons, $\Ad(\theta_{g})V_{i}$ must intersect any subspace $W_{1}'\subseteq W'$ with $\dim W'_{1}=\dim W_{1}$. Choosing such $W_{1}'\in B_{O_{\data}(\rho)}(W_{1})$, we get $ \dist_{\angle}(\Ad(\theta_g)V_{i}, W_{1})\ll_{\data} \rho$, to which point \eqref{transW1} applies and concludes the proof.

It remains to check the claim. By \Cref{avoidW'}, we know there exists  $W_{1}\subseteq W$ for which \eqref{dimW1-vs-W} holds
and such that for some $g\in G$, we have $\Ad(g)V_{i} \cap W_{1}= \{0\}$, say
\begin{equation}\label{transv-eq}
\dist_{\angle}(\Ad(g)V_{i}, W_{1})>c_{0}
\end{equation}
for some $c_{0}=c_{0}(G,W)>0$. For any $W'$ in a small neighborhood of $W$, the above inequality \eqref{transv-eq} still holds with same inputs $(g, V_{i}, c_{0})$ and $W_{1}$ replaced by an appropriate perturbation $W'_{1}\subseteq W'$. By compactness of $\Gr(\kg)$, we may thus assume the constant $c_{0}$ to be independent from $W$, i.e., $c_{0}=c_{0}(G)$. Applying \Cref{non-conc-Ad-thetagVi}, we derive \eqref{transW1}, which concludes the proof.
\end{proof}

\bigskip
\subsection{Proof of dimensional interpolation under the walk} \label{Sec-proof-dimstab-conclusion}

As a last preliminary for the proof of \Cref{RW-preserve-dim}, we combine \Cref{avoid-bad-pencil-proba} and \Cref{subcritical-projection} to show that the random subspaces $\Ad(\theta_{g}V_{i})_{g\sim \mu^n}$ satisfy a subcritical projection theorem.

\begin{lemma}
\label[lemma]{mu-Sub}
Let $D, \eps, \delta>0$, let $n\geq1$ and $i\in \{1, \dots, m\}$. If $D^{-1},\eps\lll_{\mu}1$; $\delta\lll_{\data,\mu, \eps}1$; and $n\geq \sqrt{\eps}|\log\delta|$, then the distribution of $\Ad(\theta_{g}V_{i})_{g\sim \mu^n}$ satisfies $\SubP$ with parameters $(\delta, \eps, D\sqrt{\eps})$.
\end{lemma}

\begin{proof}
Taking $\delta\lll_{\data, \mu, \eps}1$ and noting the assumption on $n$ means $\delta^{\sqrt\eps}\geq e^{-n}$, we may apply  \Cref{avoid-bad-pencil-proba} to get for every non-zero subspace $W\subseteq \Gr(\kg)$
\begin{equation*}
\mu^n \left\{g\,:\, \max_{W'\in B_{\delta^{\sqrt\eps}}(W)}\frac{\dim W'\cap \Ad(\theta_{g})V_{i}}{\dim W'} > \frac{\dim V_{i}}{\dim \kg} \right\} \leq \delta^{c \sqrt\eps}
\end{equation*}
where $c=c(\mu)>0$.
Provided $\eps\lll_{\mu}1$, this allows us to apply our subcritical projection theorem (\Cref{subcritical-projection} and the first remark that follows it) to the random variable $(\Ad(\theta_{g})V_{i})_{g\sim \mu^n}$.
The claim follows.
\end{proof}

We can now combine Lemmas \ref{net-charts}, \ref{descrip-AdgB1}, \ref{mu-Sub}, and \Cref{subc-mult}, to conclude the proof of \Cref{RW-preserve-dim}.

\begin{proof}[Proof of \Cref{RW-preserve-dim}] Consider $s$, $\delta$, $\nu$, $(t_{i})_{i}$, $n$ as in \Cref{RW-preserve-dim}. Let $\eps>0$ be a parameter to be specified below. Apply \Cref{net-charts} with $r=\delta^{2/3}$, let $\varphi : \{\inj\geq \delta^{2/3}\} \rightarrow B^\kg_{1}$ the associated map, set $\tnu=\varphi_{\star}\nu$. Using \Cref{net-charts} item 1), and assuming $s\leq \frac{1}{3\lambda_{1}}$ (so $\delta^{1-s\lambda_{i}}\leq \delta^{2/3}$) and $\delta\lll_{\data,\eps}1$, we have for $i=1, \dots, m+1$, $$\sup_{v\in \kg}\tnu(B^\kg_{\delta^{1-s\lambda_{i}}}+v)\leq \delta^{-\eps}t_{i}.$$

We aim to plug these estimates in \Cref{subc-mult}, applied with $\eta=1$, and the random box $B^{\uV_{\theta_{g}}}_{\delta^\ur}:=\sum_{i}B^{\Ad(\theta_{g})V_{i}}_{\delta^{1-s\lambda_{i}}}$ where $g\sim \mu^n$. For that, we first need to check the required subcritical property for the subspaces $\Ad(\theta_{g})V_{i}$ : 
by \Cref{mu-Sub}, taking $\delta\lll_{\data, \mu, \eps}1$ and $\eps \lll_{\mu} s^2$ (so that $\delta^{r_{i+1}\sqrt\eps}\geq e^{-n}$),
the distribution of $(\Ad(\theta_{g})V_{i})_{g\sim \mu^n}$ satisfies $\SubP$ with parameters $(\delta^{r_{i+1}}, \eps, D\sqrt{\eps})$ where $D=D(\mu)>1$ is a large enough constant.
Up to increasing $D$, 
\Cref{subc-mult} then yields a set $E_{1}\subseteq G$ with $\mu^n(E_{1})\leq \delta^{\eps/D}$ such that for each $g\in G\smallsetminus E_{1}$, there is $\tF_{g}\subseteq \kg$ with measure $\tnu(\tF_{g})\leq \delta^{\eps/D}$ and such that
\begin{equation}\label{eq-control-tnuAdgB}
\sup_{v\in \kg}\tnu_{|\kg \smallsetminus \tF_{g}}\left(B^{\uV_{\theta_{g}}}_{\delta^\ur} +v \right)\leq \delta^{-D\sqrt\eps}\prod_{i}t^{j_{i}/d}_{i}.
\end{equation}

Moreover, by \Cref{descrip-AdgB1}, there is a subset $E_{2}\subseteq G$ such that
 $\mu^n(E_{2})\leq \delta^{\gamma}$ where $\gamma=\gamma(\mu, s, \eps)\in (0,1)$ and for $g\notin E_{2}$, the box $\Ad(g)B^\kg_{\delta}$ satisfies
\begin{equation}\label{eq-control-AdgB}
\delta^{\eps}B^{\uV_{\theta_{g}}}_{\delta^\ur}\subseteq \Ad(g)B^\kg_{\delta} \subseteq \delta^{-\eps}B^{\uV_{\theta_{g}}}_{\delta^\ur}.
\end{equation}

Put together, \eqref{eq-control-tnuAdgB} and \eqref{eq-control-AdgB} yield that for $g\in G\smallsetminus (E_{1}\cup E_{2})$,
$$\sup_{v\in \kg}\tnu_{|\kg \smallsetminus \tF_{g}}\left(\Ad(g)B^\kg_{\delta} +v \right)\leq \delta^{-2D \sqrt\eps}\prod_{i}t^{j_{i}/d}_{i}$$
provided $\eps\lll_{d}1$.

We now get back to $X$. Assume $s\leq \frac{1}{4\lambda_{1}}$, so that $\delta^{1-s\lambda_{i}}\in [\delta^{3/4}, \delta^{5/4}]$ for every $i$. Provided $\eps \lll 1$, we deduce from \eqref{eq-control-AdgB} that $B^\kg_{\delta^{4/3}}\subseteq \Ad(g)B^\kg_{\delta} \subseteq B^\kg_{\delta^{2/3}}.$
This allows to apply \Cref{net-charts} item 2), which yields that the subset $F_{g}=\varphi^{-1}(\tF_{g})\subseteq X$ satisfies $\nu(F_{g})\leq \delta^{\eps/D}$ and
$$\sup_{x\in X}\nu_{|X \smallsetminus F_{g}}\left(gB^G_{\delta}x \right)\leq \delta^{-3D\sqrt\eps}\prod_{i}t^{j_{i}/d}_{i}. $$
This concludes the proof, by taking $\eps$ small enough so that $3D\sqrt\eps\leq \eps_{1}$ and imposing $\eps_{2} \leq \frac{1}{2}\min(\eps/D, \gamma)$. 
\end{proof}

\section{Random walks increase dimension at one scale or another} \label{Sec-RWincrease}

In this section, we establish a supercritical decomposition property for the action of a Zariski-dense random walk on a simple homogeneous space. The main result is \Cref{sup-mult-X}. It implies \Cref{pr:splitnu}, and thus validates the second of the two key steps on which the main results of the paper rely (see \Cref{reductions}).

\bigskip
We keep the notations $G$, $K$, $\ka$, $\Phi^+$, $\norm{\cdot}$, $\Lambda$, $X$, $\mu$, $(\lambda_{i})$, $\data$ from \Cref{Sec-dim-pres}. \Cref{sup-mult-X} below ensures that a measure $\nu$ on $X$ which has dimension $\alpha$ above a scale $\delta$ can be partitioned into two submeasures $\nu=\nu_{1}+\nu_{2}$ such that for some $n=n(\mu, \delta)\geq 1$, and most $g\sim \mu^n$, the convolution $\delta_{g^{-1}}*\nu_{p}$ has improved dimension $\alpha+\eps$ at some appropriate scales $\delta^{t_{p}}$ where $t_{1}, t_{2}\in (0, 1)$ are absolute constants.

\begin{thm} \label{sup-mult-X}
Let $\varkappa, \eps, \delta \in (0, 1/2)$. 
Let $1 > t_1 > t_2 > 0$ be parameters such that the constant 
$t:=(t_{1}-t_{2})/(\lambda_{1}+\lambda_{2})$
  satisfies $t \lambda_{1}<\min(1-t_{1}, t_{2})$ and $3t\lambda_{1}<2t_{2}-t_{1}$.

Let $\nu$ be a Borel measure on $X$ of mass at most $\delta^{-\eps}$, which is supported on the compact set $\bigl\{\inj \geq \delta^{\frac{t_{2}-t\lambda_{1}}{2}}\bigr\}$, and satisfies for some $\alpha\in [\varkappa, 1-\varkappa]$, for all $\rho\in \bigl\{\delta^{t_{1}- t\lambda_{i}}\bigr\}_{i=1}^{m+1}\cup \bigl\{\delta^{t_{2}- t\lambda_{i}}\bigr\}_{i=1}^{m+1} \cup \bigl[\delta^{t_{1}- t\lambda_{2}}, \delta^{t_{1} - t\lambda_{1}}\bigr]$,
 $$\sup_{x\in X}\nu\bigl(B^G_{\rho}x\bigr) \leq \delta^{-\eps}\rho^{d\alpha}.$$

If $\eps, \delta\lll_{\data, \mu, \varkappa, (t_{p})_{p}}1$, then we can write $\nu=\nu_{1}+\nu_{2}$ where $\nu_{1}, \nu_{2}$ are mutually singular Borel measures satisfying the following.
Set $n=\lfloor |t \log \delta|\rfloor$.
There exists $\cE \subseteq G$ such that $\mu^n(\cE)\leq \delta^\eps$, and for every $p\in \{1, 2\}$,  $g\in G\smallsetminus \cE$, there exists $F_{p,g}\subseteq X$ satisfying $\nu(F_{p,g})\leq \delta^\eps$ and
$$\sup_{x\in X}\nu_{p | X\smallsetminus F_{p,g}} \bigl(g B^G_{\delta^{t_{p}}}x\bigr)\leq \delta^{t_{p}d (\alpha+\eps)}. $$
\end{thm}

%\begin{thm} \label{sup-mult-X}
 %Let $\varkappa, \eps, \delta \in (0, 1/2)$. Let $\nu$ be a Borel measure on $X$ which is supported on $\{\inj \geq \delta^{\eps}\}$ and satisfies for some $\alpha\in [\varkappa, 1-\varkappa]$, for all $\rho\geq \delta$,
 %$$\sup_{x\in X}\nu(B^G_{\rho}x) \leq \delta^{-\eps}\rho^{d\alpha}.$$
 %Let $t_{1},t_{2}\in (0,1)$ such that $t_{1}>t_{2}$ and setting $t=(t_{1}-t_{2})/(\lambda_{1}+\lambda_{2})$, we have $t \lambda_{1}<\min(1-t_{1}, t_{2})$ and $3t\lambda_{1}<2t_{2}-t_{1}$.

 %If $\eps, \delta\lll_{\data, \mu, \varkappa, (t_{p})_{p}}1$, then we can write $\nu=\nu_{1}+\nu_{2}$ where $\nu_{1}, \nu_{2}$ are mutually singular Borel measures satisfying the following. Set $n=\lfloor |t \log \delta|\rfloor$. There exists $E\subseteq G$ such that $\mu^n(E)\leq \delta^\eps$, and for all $g\in G\smallsetminus E$, $p\in \{1, 2\}$, there exists $F_{g,p}\subseteq X$ satisfying $\nu(F_{g,p})\leq \delta^\eps$ and
 %$$\sup_{x\in X}\nu_{p |X\smallsetminus F_{g,p}} (gB^G_{\delta^{t_{p}}}x)\leq \delta^{t_{p}d (\alpha+\eps)}. $$
%\end{thm}

\begin{remark} We explain the condition on $t_{1},t_{2}$.
Recall that for $g\sim \mu^n$, the box $\Ad(g)B^\kg_{\delta^{t_{p}}}$ can essentially be written $\Ad(\theta_{g})\sum_{i=1}^{m+1}B^{V_{i}}_{\delta^{t_{p} -t\lambda_{i}}}$ where $(V_{i})_{i=1}^{m+1}$ is a partial flag of $\kg$ determined by $\mu$, and $\theta_{g}$ denotes the first Cartan component of $g$ (see \Cref{descrip-AdgB1}).
The parameter $t$ is chosen so that the largest two side lengths of the box associated to $t_{1}$ correspond to the smallest two side lengths of the box associated to $t_{2}$. The condition $t \lambda_{1}<\min(1-t_{1}, t_{2})$ guarantees the exponents $t_{1} -t\lambda_{i}$ and $t_{2} -t\lambda_{i}$ are in $(0,1)$ for every $i$. Finally, the requirement $3t\lambda_{1}<2t_{2}-t_{1}$ further guarantees the boxes are not too distorted, meaning the exponents in fact all belong to some interval of the form $(\zeta,2\zeta)$ where $\zeta\in (0, 1)$. This last requirement is important to justify that additive translates $(\Ad(g)B^\kg_{\delta^{t_{p}}}+v)_{v\in \kg}$ represent the sets $(g B_{\delta_{t_{p}}}x)_{x\in X}$ in suitable charts (see \S\ref{Sec-covering charts}).

An example of suitable exponents $t_{1}, t_{2}$ is given by
$$t_{1}=1/2 \qquad t_{2}=7/16.$$
Note this choice is valid regardless of $\mu$.
\end{remark}

Let us sketch the strategy to prove \Cref{sup-mult-X}. Using \Cref{net-charts}, we first linearize $X$ at an appropriate scale (depending on $\mu$, $(t_{p})_{p}$, $\delta$). For $p\in \{1,2\}$, $n=\lfloor |t \log \delta|\rfloor$, and $\mu^n$-most $g$, this allows us to see translates of balls $gB^G_{\delta^{t_{p}}}x$ as Euclidean boxes $\Ad(g)B^\kg_{\delta^{t_{p}}}+v$ in the Lie algebra. Then we apply the multislicing supercritical decomposition \Cref{mult-sup-dec} to those boxes. To apply the latter, it is crucial to check that the corresponding random partial flag as $g\sim \mu^n$ satisfies a suitable supercritical alternative. Establishing this estimate is the essence of the present section. To do so, we first highlight a tailored version of the supercritical projection theorem,  \Cref{supercrit-weak-assump}, which takes place in a specific geometric setting, but has the advantage of involving less restrictive non-concentration assumptions on the random projector. Validating such non-concentration in the context of random walks is nonetheless quite subtle, this is the content of  \S\ref{Sec-nc-highestweight}.  
We then show  in \S\ref{Sec-supercrit-alt} that the supercritical alternative can be forced within  the framework of \Cref{supercrit-weak-assump} by replacing the involved set $A$ with a suitable slice. This part is also innovative. It is precisely to perform this reduction to a slice that we are constrained to prove a  dimensional increase \emph{alternative} as opposed to  a  dimensional increase at a \emph{definite} single scale.

\subsection{Projection theoretic preliminaries} \label{Sec-back-proj}
We record some useful results on discretized projection theorems.

\bigskip

For a subset $A\subseteq \R^d$, and $\alpha, \eps, \delta>0$, we set
\begin{equation}\label{defcO}
\begin{split}
\cO^{(\alpha, \eps)}_{\delta}(A) \defeq \bigl\{\, V\in \Gr(\R^d) : \exists A' \subseteq A \,\,&\text{ with }\,\,\cN_{\delta}(A')\geq \delta^\eps \cN_{\delta}(A)\\
& \text{ and }\, \,\cN_{\delta}(\pi_{V}A') < \delta^{-\alpha \dim V -\eps}\bigr\}.
\end{split}
\end{equation}
Note $\cO^{(\alpha, \eps)}_{\delta}(A)$ is dual to the exceptional set $\cE^{(\alpha, \eps)}_{\delta}(A)$ from \eqref{notation-cEaed}, as here we consider projections \emph{onto} rather than \emph{parallel to}.

The next result is a supercritical estimate under relatively weak non-concentration assumptions but in a specific geometric setting.

\begin{proposition}
\label[proposition]{supercrit-weak-assump}
Let $d,k\geq 1$ be integers with $d\geq 2k$.
Let $\varkappa, c, \eps, \delta\in (0, 1/2)$.
Consider $E_{1}, E_{2}\in \Gr(\R^d, k)$ such that $\dang(E_{1}, E_{2})\geq \delta^\eps$ and write $E=E_{1}\oplus E_{2}$.
If $k=1$, set $\sS=\Gr(E,1)$. 
%If $k\geq 2$, set 
%$$\sS =\{W\in \Gr(E,k)\,:\, W\in \{E_{1}, E_{2}\} \text{ or } W=\R v \oplus H \text{with }  \R v\in \Gr(E_{1},1), H\in \Gr(E_{2},k-1) \}.$$
If $k\geq 2$, set $\sS \subseteq \Gr(E,k)$ the collection of subspaces $W$ satisfying either $W\in \{E_{1}, E_{2}\}$ or $W=\R v+H$ where $\R v$  is a line  in $E_{2}$ and $H$ is a $(k-1)$-plane in $E_{1}$. 

Let $\sigma$ be a probability measure on $\Gr(\R^d, k)$ satisfying 
\begin{equation}
\label{eq:NCsS}
\forall \rho \geq \delta,\,\forall  W\in \sS,\quad
\sigma\set{V\,:\, \dang(V, W^\perp) <\rho}\leq \delta^{-\eps}\rho^c.
\end{equation}

Let $A \subseteq B^{E}_{1}$ such that $\cN_{\delta}(A)\geq \delta^{-2k\alpha+\eps}$ with $\alpha\in [\varkappa, 1-\varkappa]$, and satisfying
$$\sup_{i=1,2}\cN_{\delta}(\pi_{E_{i}}A)\leq \delta^{-k\alpha-\eps}, $$
while for all $\rho>\delta$,
$$\cN_{\rho}(\pi_{E_{1}}A)\geq \delta^{\eps} \rho^{-c}. $$

If $\eps, \delta \lll_{d, \varkappa, c} 1$, then
 $$\sigma\bigl( \cO_{\delta}^{(\alpha, \eps)}(A) \bigr)\leq \delta^{\eps}.$$
\end{proposition}

\begin{proof}
First, consider the special case $d = 2k$, i.e. $E = \R^d$.
In this case, the proof can be abstracted from \cite[Proof of Theorem 3]{Bourgain2010} for $k=1$, and more generally from \cite[Proof of Proposition 7]{He2020JFG} for abitrary $k$. It exploits Balog-Szemerédi-Gowers' theorem and the Pl\"unnecke-Ruzsa inequality to reduce the problem to a sum-product estimate for matrix algebras~\cite[Theorem 3]{He2019}.

Now we reduce the general case to this special case.
Let $\sigma$ and $A$ be as in the statement.
For $W \in \Gr(E)$, writing $W^\perp=(W^\perp \cap E) \oplus^\perp E^\perp$, we see directly from the definition of $\dang$ that
%\cite[Lemma 16]{He2020JFG}, 
\begin{equation}
\label{eq:lemma16}
 \dang(V, W^\perp) = \dang(V, E^\perp) \dang(\pi_E V, W^\perp \cap E).
\end{equation}
Hence, the non-concentration~\eqref{eq:NCsS} implies
\[
\forall \rho \geq \delta,\, \forall W\in \sS,\quad \sigma\set{V\,:\, \dang(\pi_E V, W^\perp) <\rho}\leq \delta^{-\eps}\rho^c.
\]
Thus viewing $E$ as $\R^{2k}$, we may apply the special case to the distribution $(\pi_E)_\star \sigma$.

For any subset $A' \subset A$, by~\cite[Lemma 18]{He2020JFG}, we have
\[
\forall V \in \Gr(\R^d),\quad \cN_\delta(\pi_{(\pi_E V)} A') \ll_d \dang(V, E^\perp)^{-d} \cN_\delta(\pi_V A').
\]
It follows that for every $V \in \cO_{\delta}^{(\alpha, \eps)}(A)$, either
\(
\pi_E V \in \cO_{\delta}^{(\alpha, (3dc^{-1} + 2)\eps)}(A),
\)
or
\(
\dang(V, E^\perp) < \delta^{3c^{-1}\eps}.
\)

To conclude, note that the probability of the first event can be bounded using the special case applied to $(\pi_E)_\star \sigma$: for $\eps, \delta \lll_{d,\varkappa, c} 1$,
\[
\sigma\setbig{ V : \pi_E V \in \cO_{\delta}^{(\alpha, (3dc^{-1} + 2)\eps)}(A)} \leq \delta^{2 \eps}.
\]
The probability of the second event can be bounded, noting that $E^\perp \subseteq E_1^\perp$ and using the non-concentration~\eqref{eq:NCsS}:
\[
\sigma\setbig{V : \dang(V, E^\perp) < \delta^{3c^{-1}\eps}} \leq \sigma\setbig{V : \dang(V, E_1^\perp) < \delta^{3c^{-1}\eps}} \leq \delta^{2 \eps}.
\]
Adding these bounds together finishes the proof of the proposition.
\end{proof}

%Let us pause here to briefly outline the proof of the supercritical alternative $\SAP$ we are aiming for.
%The main tool is \Cref{supercrit-weak-assump}.
%In order to use it we need, given a set $A \subset \kg$, to find a suitable slice $A \cap E$ in some suitable subspace $E = E_1 \oplus E_2 \in \Gr(\kg, 2 \dim V_1)$, establish appropriate non-concentration~\eqref{eq:NCsS} for the directions, and show a non-concentration property for the projection of the slice $\pi_{E_1}(A \cap E)$.
%The latter is dealt with right below using existing tools.
%The non-concentration of the direction of the projections, in view of \Cref{non-conc-Ad-thetagVi}, can be reduced to a Lie theoretic statement and this will be established in \Cref{Sec-nc-highestweight}. 
%The construction of the suitable slice, the most innovative part of this argument, is postponed to \Cref{Sec-supercrit-alt}.

\bigskip
\Cref{nc-proj} below claims that if a set $A$ satisfies some non-concentration in the sense that the uniform probability measure on $A$ is Frostman, then most of the projections of $A$ do as well, up to passing to a large subset of $A$ on which we have some control. 
This improves upon previous results of \cite{Bourgain2010, He2020JFG} which rely on a stronger form of non-concentration regarding projectors. \Cref{nc-proj} will be used later in \S\ref{Sec-supercrit-alt} to check the non-concentration property required for $A$ to apply \Cref{supercrit-weak-assump}. 
We use the shorthand $x^{(\rho)}:=B^{\R^d}_{\rho}(x)$, for $x\in \R^d$ and $\rho>0$.

\begin{lemma}
\label[lemma]{nc-proj}
Let $d>k\geq1$ be integers, let $\eps>0$, $ c,\alpha, \delta\in (0, 1)$. Let $\sigma$ be a probability measure on $\Gr(\R^d, k)$ satisfying for all non-zero subspaces $W\subseteq \R^d$, all $\rho \geq \delta$,
\begin{equation*}
\sigma \left\{V\,:\, \max_{W'\in B_{\rho}(W)}\frac{\dim V \cap W'}{\dim W'} > \frac{\dim V}{d} \right\} \leq \delta^{-\eps} \rho^{c}.
\end{equation*}

 Let $A \subseteq B^{\R^{d}}_{1}$ such that for all $x\in \R^d$, $\rho \geq \delta$,
 $$\cN_{\delta}\bigl(A\cap x^{(\rho)}\bigr)\leq \delta^{-\eps}\rho^{\alpha}\cN_{\delta}(A).$$

Let $D>1$, let $\cG$ be the set of $V\in \Gr(\R^d, k)$ satisfying the following: for every $A'\subseteq A$ with $\cN_{\delta}(A')\geq \delta^{\eps} \cN_{\delta}(A)$, there exists $A''\subseteq A'$ with $\cN_{\delta}(A'')\geq |\log \delta |^{-D} \cN_{\delta}(A')$ and satisfying $\cN_{\delta}(\pi_{V}A''\cap y^{(\rho)})\leq \delta^{-D\sqrt\eps}\rho^{\frac{k}{2d}\alpha}\cN_{\delta}(\pi_{V}A'')$ for all $\rho\geq \delta$, $y\in V$.

$$\text{If $D\ggg_{d,c}1 $ and $\delta\lll_{d, c,\eps}1$, then $\sigma(\cG)\geq 1- \delta^{\eps}.$}$$
\end{lemma}

The idea of proof is to use the subcritical projection Theorem~\ref{subcritical-projection} (in its single slicing form given by \Cref{subc-mult}) to see that the uniform measure on $A$ typically has  $\alpha k/d$-dimensional projections above the scale $\delta$, and deduce the announced result by a regularization argument.
\begin{proof}
We may argue under the extra condition $\eps\lll_{d, c}1$ otherwise the claim is trivial (by taking larger $D$).
We may assume $A$ to be $2\delta$-separated. 

We set $\nu$ the uniform probability measure on $A$. The non-concentration assumption on $A$ reads as:
$$\forall x\in \R^d,\, \forall \rho \geq \delta,\quad \nu(x^{(\rho)})\leq \delta^{-\eps}\rho^{\alpha}.$$
 %Let $\eta:=\sqrt{2\eps/c}$. Note that for $\rho \in [\delta, \delta^{\eta}]$, the upper bound $\delta^{-\eps}\rho^{\eta c}$ appearing in the non-concentration for $\sigma$ satisfies $\delta^{-\eps}\rho^{\eta c} \leq \rho^{\eta c/2}$.

We now apply our subcritical projection \Cref{subcritical-projection} to the random projector $(\pi_{V})_{V\sim \sigma}$.
More precisely, assume $\sqrt{\eps}<c/2$.
Given $\rho \in [\delta, \delta^{2\sqrt{\eps}/c}]$, set $\eps_{\rho}\in [\eps,1]$ such that $\rho^{\eps_{\rho}}=\delta^\eps$. Note that $\rho^{\sqrt{\eps_{\rho}}}\in [\delta, \delta^{2\eps/c}]$ due to $\rho^{\sqrt{\eps_{\rho}}}\leq \rho^{\sqrt{\eps}}$ and the prescribed upper bound on $\rho$. The non-concentration assumption on $\sigma$ yields
\begin{equation*}
\sigma \left\{V\,:\, \max_{W'\in B_{\rho^{\sqrt{\eps_{\rho}}}}(W)}\frac{\dim V \cap W'}{\dim W'} > \frac{\dim V}{d} \right\} \leq \rho^{c\sqrt{\eps_{\rho}}/2}.
\end{equation*}
Note this non-concentration estimate on $(V)_{V\sim \sigma}$ is also valid for $(V^\perp)_{V\sim \sigma}$ thanks to \eqref{grasmannian-distance} and \Cref{orthog-red-sub}.
Provided $D\ggg_{d,c}1$ and $\delta\lll_{d,c,\eps}1$, \Cref{subcritical-projection} (and the first remark following it) then guarantees that the distribution $(V^\perp)_{V\sim \sigma}$ satisfies $\SubP$ with parameters $(\rho, \eps_{\rho}, D \sqrt{\eps_{\rho}})$. We recall $\SubP$ was introduced in \Cref{Sub-crit-P}.

The two previous paragraphs allow to apply the slicing estimate from \Cref{subc-mult} with the random box $(B^{V^\perp}_{1}+B^{\R^d}_{\rho})_{V\sim \sigma}$ and the exponent $\eps_{\rho}$.
Up to taking larger $D$ and $\delta\lll_{d,\eps}1$, we obtain some $\cE_{\rho}\subseteq \Gr(\R^d)$ with $\sigma(\cE_{\rho})\leq \delta^{\eps/D}$ and such that for every $V\in \Gr(\R^d, k)\smallsetminus \cE_{\rho}$, there exists $F_{\rho, V}\subseteq \R^d$ such that $\nu(F_{\rho, V})\leq \delta^{\eps/D}$ and for all $y\in V$,
$$(\pi_{V}\nu_{|\R^d\smallsetminus F_{\rho, V}})(y^{(\rho)})\leq \delta^{-D\sqrt{\eps}}(\delta^{-\eps}\rho^{\alpha})^{k/d}\leq \delta^{-2D\sqrt\eps}\rho^{\alpha'} $$
where $\alpha':=\alpha k/d$.

Note that such an estimate automatically upgrades to a half-neighborhood of $\rho$, namely for all $y\in V$, $r\in [\rho^{2}, \rho]$,
$$(\pi_{V}\nu_{|\R^d\smallsetminus F_{\rho, V}})(y^{(r)})\leq \delta^{-2D\sqrt\eps}r^{\alpha'/2}.$$
Let $(\rho_{i})_{i\in I}$ be a collection of real numbers $\rho_{i}\in \bigl[\delta, \delta^{2\sqrt{\eps}/c}\bigr]$ such that \[
\bigl[\delta, \delta^{2\sqrt{\eps}/c}\bigr]\subseteq \bigcup\nolimits_{i \in I}\bigl[\rho^{2}_{i}, \rho_{i}\bigr]
\]
and $\abs{I} \leq O( \abs{\log \eps})$. Set $\cE = \bigcup\nolimits_{i \in I} \cE_{\rho_{i}}$, and for $V\notin \cE$, set $F_{V}=\bigcup\nolimits_{i \in I} F_{\rho_{i}, V}$.
Then $\sigma(\cE), \nu(F_{V})\leq \delta^{\eps/(2D)}$ and for all $y\in V$, $\rho\in \bigl[\delta, \delta^{2\sqrt{\eps}/c}\bigr]$,
$$(\pi_{V}\nu_{|\R^d\smallsetminus F_{V}})(y^{(\rho)})\leq \delta^{-2D\sqrt\eps}\rho^{\alpha'/2}.$$
As $\nu$ has mass $1$, this inequality also holds in the range $\rho\geq \delta^{2\sqrt{\eps}/c}$.
Getting back to $A$, % and noting $ |A|\ll |A\smallsetminus F_{V}|$ for $\delta\lll_{D,\eps}1$, we obtain 
this means for all $y\in V$, $\rho\geq \delta$,
\begin{equation}\label{eq-non-conc-proj-meas}
\absbig{(A\smallsetminus F_{V}) \,\cap \pi_{V}^{-1}y^{(\rho)}} \ll \delta^{-2D\sqrt\eps}\rho^{\alpha'/2} \abs{A}.
\end{equation}

Let $A'\subseteq A$ be a subset such that $|A'|\geq \delta^{\eps/(4D)} |A|$. Using \Cref{reg-measures}, we can extract $A''\subseteq A'\smallsetminus F_{V}$ such that $|A''|\gg_{d} |\log \delta|^{-O(1)} |A'|$ and which
 is regular\footnote{This phrasing is slightly abusive because $\cD_{\delta}$ is a priori not finer than $\pi_{V}^{-1}\cD_{\delta}$ in the sense given in \Cref{Sec-regularization}.
 However, we can consider a partition $\cP$ which is finer than $\pi_{V}^{-1}\cD_{\delta}$ and equivalent to $\cD_{\delta}$ in the sense that every $\cP$-cell is covered by $O_{d}(1)$ $\cD_{\delta}$-cells and conversely. In the argument above, we really mean $\cP$ instead of $\cD_{\delta}$.}
for $\pi_{V}^{-1}\cD_{\delta} \prec \cD_{\delta}$.
Then \eqref{eq-non-conc-proj-meas} becomes % still holds with $(A\smallsetminus F_{V}, \eps)$ replaced by $(A'', 2\eps)$ , 
\[
\forall y \in V,\, \forall \rho \geq \delta,\quad \absbig{A'' \,\cap \pi_{V}^{-1}y^{(\rho)}} \leq \delta^{-4D\sqrt\eps}\rho^{\alpha'/2} \abs{A''}.
\]
Dividing each side by the $\delta$-covering number of the intersection of $A''$ with a $\delta$-tube of axis parallel to $V^\perp$, we obtain
\begin{equation*}
\absbig{\pi_{V}A'' \cap y^{(\rho)}} \leq \delta^{-4D\sqrt\eps}\rho^{\alpha'/2}|\pi_{V}A''|.
\end{equation*}

In conclusion, we have seen that for some $D=D(d,c)>1$, for all $\eps\lll_{d,c}1$ and $\delta\lll_{d,c,\eps}1$, if we write $\cG'$ the set defined as $\cG$ but with $\eps$ replaced by $\eps':=\eps/(4D)$ and $D$ replaced by $D'=8D^{3/2}$, then we have $\sigma(\cG')\geq 1-\delta^{\eps'}$.
Then, arguing with $4D\eps$ from the start (noting here the assumptions on $\sigma$ and $A$ are still valid for this exponent), we obtain the desired estimate.
%replaced by $\eps/4D$ and $D$ replaced by $8D^3$. Noting the assumption of the lemma is stable by increasing $\eps$ to $4D\eps$, we obtain $\sigma(\cE_{reg})\geq 1-\delta^\eps$ provided $D\ggg_{d, c}1$, $\eps\lll_{D}1$, $\delta \lll_{d,c,D, \eps}1$. Noting the conclusion is stable by increasing $D$, we may reduce those constraints to $\eps\lll_{d,c}1$, $\delta \lll_{d,c, \eps}1$, thus finishing the proof of the lemma.
\end{proof}

\subsection{Non-concentration properties of highest weight subspaces} \label{Sec-nc-highestweight}

Recall from \Cref{descrip-AdgB1} that for $g\sim \mu^n$, the set $\Ad(g)B^\kg_{1}$ is a random Euclidean box in $\kg$, whose partial flag is given by $(\Ad(\theta_{g})V_{i})_{i=1, \dots, m+1}$ where $\theta_{g}$ refers to the first component of $g$ in the Cartan decomposition and $(V_{i})_{i=1, \dots, m+1}$ denotes the partial flag associated to the Lyapunov exponent $\kappa_{\mu}\in \ka^{++}$ of $\mu$, see \eqref{Cartan-filtration}.
We established in Sections \ref{Sec-submod} and \ref{Sec-dim-pres} some weak non-concentration estimates regarding the distribution of $\Ad(\theta_{g})V_{i}$ as $g\sim \mu^n$. Those were sufficient to apply the subcritical multislicing, and ultimately show the random walk on $X$ almost preserves the dimension of a given measure. In order to obtain a dimensional gain, we need a stronger estimate for at least one of the $V_{i}$'s. We focus on $V_{1}$ which is the highest weight subspace for $\Ad(G)\acts \kg$.

In general, the distribution of the random subspace $(\Ad(\theta_{g})V_{1})_{g\sim \mu^n}$ does not satisfy the non-concentration property required in Bourgain's projection theorem \cite{Bourgain2010} and its successive upgrades in \cite{He2020JFG, Shmerkin, BH24}. Recall this condition asks that for any $W\in \Gr(\kg)$ of complementary dimension, most realizations of $\Ad(\theta_{g})V_{1}$ are in direct sum with $W$. Although this property holds when $\Ad(G)$ is proximal (e.g. $G=\SL_{N}(\R)$), it fails drastically for an arbitrary simple Lie group (e.g. $G=\SO(N,1)$ with $N\geq 5$, see remark below and also \Cref{lack}).
This section still provides non-concentration estimates with respect to a smaller class of subspaces $W$, and which we will be able to exploit later through \Cref{supercrit-weak-assump}. Keeping in mind \Cref{non-conc-Ad-thetagVi}, which allows to reduce probabilistic statements to geometric ones, we focus on describing a collection of subspaces of $\kg$ which are in direct sum with some subspace from the orbit $\Ad(G)V_{1}$.

We call $(W_{i})_{i=1, \dots, m+1}$ the partial flag associated to $-\kappa_{\mu}$, in other terms $W_{m+1}=\kg$, and for $i=1, \dots,m$,
\begin{equation} \label{def-Wi}
W_{i}:=\bigoplus_{\alpha \in \Phi\cup \{0\} \,:\, \alpha(\kappa_{\mu})\leq \lambda_{m+2-i}} \kg_{\alpha}
\end{equation}
or equivalently, thanks to our choice of norm $\norm{\cdot}$,
\[
W_i = V_{m+1-i}^{\perp}.
\]
In particular, $W_{1}$ is the lowest weight subspace for $\Ad(G)\acts \kg$. For this subsection, we set $F_{0}:=V_{m}\cap W_{m}$. Note that $$\kg=V_{1}\oplus^\perp W_{1}\oplus^\perp F_{0}.$$

Recall  the subspaces $\ku, \kp, \ku^-, \kp^-\subseteq \kg$ from the previous section, namely $\ku=\oplus_{\alpha\in \Phi^+}\kg_{\alpha}$, $\kp=\kz_{\kg}(\ka)\oplus \ku$, and  $\ku^-, \kp^-$ are defined similarly using negative roots.  
We set $P,U, P^-,U^-\subseteq G$ the connected Lie subgroups of Lie algebras $\kp, \ku, \kp^-, \ku^-\subseteq \kg$. Therefore each $V_{i}$ is $\Ad(P)$-invariant while each $W_{i}$ is $\Ad(P^-)$-invariant.

For conciseness, \emph{we will omit the notation $\Ad$ throughout \Cref{Sec-nc-highestweight}}, meaning that given $g\in G$ and $V \subseteq \kg$, we will write $gV:=\Ad(g)V$.

\begin{proposition}
\label[proposition]{non-conc-AdgV1}
For every line $\R v\subseteq \kg$, every hyperplane $H\subseteq W_{1}$, there exists $g\in G$ satisfying
$$gV_{1} \cap (\R v + H + F_{0})=\{0\}. $$
\end{proposition}

\begin{remark}
Here we need to be particularly cautious.
 \Cref{non-conc-AdgV1} relies on the particular form of the space we want to make transverse to $V_{1}$ using the action of $G$, it is not just a consideration of dimensions.
To see why, consider the case where $G=\SO(N,1)$ with $N\geq 2$. Given a line $\R v\subseteq V_{1}$, we claim there exists a subspace $F'\subseteq F_{0}$ of dimension $N-1$ such that for all $g\in G$,
$$gV_{1} \cap (\R v\oplus W_{1} \oplus F')\neq \{0\}.$$
Note that $W_1 = \ku^-$ in this situation and that $\R v\oplus W_{1} \oplus F'$ has dimension $2N-1$, which for $N$ large, is much smaller than the codimension of $V_{1}$ in $\kg$ (equivalent to $\frac{1}{2}N^2$).
To check the claim, set $F'=[W_{1}, v]$ where $[\cdot, \cdot]$ is the Lie algebra bracket. Note that $\R v\oplus W_{1} \oplus F'$ is then $U^-$-invariant. On the other hand $V_{1}$ is $P$-invariant. This justifies the claim for $g\in U^-P=P^-P$, and it automatically upgrades to all $g\in G$ because $P^-P$ is Zariski-dense in $G$.
\end{remark}

In the proof of \Cref{non-conc-AdgV1}, we will consider the complexification $\kg_\C = \C \otimes_\R \kg$ of $\kg$ and work with a root system $\Phi_\C$ of $\kg_\C$.  More precisely, we choose a Cartan subalgebra $\kh_{\C}\subseteq \kg_{\C}$, and write $\Phi_{\C}\subseteq \kh_{\C}^*$ the associated root system, $\kg_{\C}=\oplus_{\Phi_{\C}\cup \{0\}}\kg_{\C, \beta}$ the root space decomposition. We choose a set of positive roots $\Phi^+_{\C}\subseteq \Phi_{\C}$.
We denote by $\kn_{\C}=\oplus_{\beta\in \Phi^+_{\C}}\kg_{\C, \beta}$ the sum of positive root spaces, and set $\kb_{\C}=\kh_{\C}\oplus \kn_{\C}$ the associated Borel Lie algebra. Using negative roots, we define similarly $\kn^-_{\C}$, $\kb^-_{\C}$.
The choices of $\kh_{\C}$ and $\Phi^+_\C$  need to have a certain compatibility with respect to their real analogues $\ka$ and $\Phi^+$ for $\kg$.
This is the purpose of the next lemma. 

%Recall from the previous section the subspaces $\ku, \kp, \ku^-, \kp^-\subseteq \kg$ are the real analogues of $\kn_{\C}, \kb_{\C}, \kn^-_{\C}, \kb^-_{\C}$, i.e. $\ku=\oplus_{\alpha\in \Phi^+}\kg_{\alpha}$, $\kp=\kz_{\kg}(\ka)\oplus \ku$, etc.  We  write $\ku_{\C}, \kp_{\C}, \ku^-_{\C}, \kp^-_{\C}\subseteq \kg_{\C}$ for their complexifications. 

\begin{lemma}
\label[lemma]{choice-Weylchamber}
We may choose a Cartan subalgebra $\kh_{\C}$ of $\kg_\C$, a root system $\Phi_\C \subset \kh_\C^*$ with positive roots $\Phi^+_{\C}$ satisfying that $\ka\subseteq \kh_{\C}$ and for every $\beta\in \Phi^+_{\C}$, we have $\beta_{|\ka}\in \Phi^+\cup \{0\}$. In this case, we have
$$\ku_{\C}\subseteq \kn_{\C}\subseteq \kb_{\C}\subseteq \kp_{\C}.$$
\end{lemma}

\begin{proof}
The first part of the lemma follows from \cite[\S 21.8]{Borel}.

%The complexification $\ka_{\C} = \ka \otimes \C$ is a commutative $\ad$-diagonalizable subalgebra of $\kg_{\C}$.
%Hence, it must be included in some Cartan subalgebra, which we can name $\kh_{\C}$. Note  that for every $\beta\in \Phi_{\C}\cup \{0\}$, we have $\beta_{|\ka}\in \Phi \cup \{0\}$. 
%We now  choose a suitable family of positive roots for $\Phi_{\C}$.
%Set $E\subseteq \kh^*_{\C}$ to be the real vector space spanned by $\Phi_{\C}$, set $F=\ka^*$.
%Let $p: E\rightarrow F, \gamma \mapsto \gamma_{|\ka}$.
%Let $\psi : F\rightarrow \R$ be a linear form such that $\Phi\cap \{\psi>0\} = \Phi^+$.
%Let $\varphi=\psi\circ p$.
%Note $\varphi$ may vanish on some roots from $\Phi_{\C}$.
%However, considering a small perturbation, we may find a linear form $\varphi' \in E^*$ such that $ \Phi_{\C}\cap \ker \varphi'=\{0\}$ and $ \Phi_{\C}\cap \{\varphi>0\}\subseteq \Phi_{\C}\cap \{\varphi'>0\}$.
%The set $\Phi^+_{\C}= \Phi_{\C}\cap \{\varphi'>0\}$ yields the desired set of positive roots.

To show $\ku_{\C}\subseteq \kn_{\C}\subseteq \kb_{\C}\subseteq \kp_{\C}$, note that for every $\alpha\in \Phi\cup \{0\}$, we have
$$\kg_{\alpha} \otimes \C = \oplus \{\kg_{\C, \beta}\,:\, p(\beta)=\alpha\}.$$
In particular for $\alpha \in \Phi^{+}$, all the roots $\beta$ contributing in the decomposition must be in $\Phi^+_{\C}$. This justifies $\ku_{\C}\subseteq \kn_{\C}$. As $\kb_{\C}=\oplus \{ \kg_{\C, \beta}\,:\, \beta\in \Phi^+_{\C}\cup \{0\}\}$, $\kp_{\C}=\oplus \{ (\kg_{\alpha})_{\C}\,:\, \alpha\in \Phi^+\cup \{0\}\}$, and $\beta\in \Phi^+_{\C}\cup \{0\}$ implies $p(\beta)\in \Phi^+\cup \{0\}$, we also have $\kb_{\C}\subseteq \kp_{\C}$.
\end{proof}

\begin{proof}[Proof of \Cref{non-conc-AdgV1}]
Recall that a simple Lie algebra over $\R$ is either absolutely simple or admits a complex structure.
In the second case, we have $\dim V_1 = 2$ and the stronger transversality that  for any $W \in \Gr(\kg, \dim \kg -2)$, there is some $g\in G$ such that $g V_1 \cap W = \{0\}$.
This is shown in \cite[Proposition 2.5]{He2020IJM}.
Thus, we may assume that $\kg$ is absolutely simple.
We will use that the root system of its complexification $\kg_{\C}$ is irreducible.

We may also assume $\R v \perp H \oplus F_{0}$.
Moreover, if $\R v \not\subset V_{1}$, then the lemma follows by taking $g=\Id$. We thus focus on the case $\R v\subseteq V_{1}$.

Assume the claim fails, then for every $u \in W_1$ there is a non-zero vector $v_1 \in V_1$ such that $\exp(u) v_1 \in \R v + F_0 + H$.
Expanding the exponential, we have $\exp(u) v_1 = v_1 + \ad(u) v_1 + \frac{1}{2} \ad(u)^2 v_1$ with $\ad(u)v_1 \in F_0$ and $\ad(u)^2 v_1 \in W_1$.
It follows that $v_1 \in \R v$ and we deduce that
\begin{equation*}
%\label{eq:adW2inH}
\forall u \in W_1,\quad \ad(u)^2 v \in H.
\end{equation*}

Consider the complexifications $V_{1, \C}$ and $W_{1, \C}$ of $V_1$ and $W_1$
and let $\cT$ denote the set of pairs $(L,S) \in \Gr(V_{1, \C}, 1) \times \Gr(W_{1, \C}, k - 1)$ satisfying
\begin{equation*}
%\label{eq:adW2inHoverC}
\forall u \in W_{1, \C},\quad \ad(u)^2 L \subset S.
\end{equation*}
By the above, it contains an element $(\C v, \C \otimes_{\R} H)$ defined over $\R$.
%We first show that $(\kg_{\C,\alpha_{\max}}, \oplus_{\alpha \in \cA(W_{1,\C}) \setminus \{-\alpha_{\max}\}} \kg_{\C,\alpha})$ is in the Zariski-closure of the $M$-orbit of $(v,H)$ and then show that this in not possible.

Fix a Cartan subalgebra $\kh_{\C}\subseteq \kg_{\C}$ containing $\ka$, and a choice of positive roots $\Phi^+_{\C}\subseteq \Phi_{\C}$ compatible with $\Phi^+$, as in \Cref{choice-Weylchamber}.
Set $\cB:=\Phi_{\C}\cup \{0\}$ and write $\kg_{\C}= \oplus_{\beta\in \cB }\kg_{\C, \beta}$ the root space decomposition of $\kg_{\C}$ (Cartan subalgebra included).
Note that $V_{1, \C}$ and $W_{1, \C}$ can be decomposed as subsums of root spaces, and write accordingly
\begin{equation*}
%\label{eq:decompV1C}
V_{1,\C}= \oplus_{\beta\in \cB(V_{1}) }\kg_{\C, \beta},\quad 
W_{1,\C}= \oplus_{\beta\in \cB(W_{1}) }\kg_{\C, \beta}.
\end{equation*}
Necessarily $\cB(V_{1})=- \cB(W_{1})$, and $\cB(V_{1})$ contains the highest root $\beta_{\max}$ (which is well defined because $\kg_{\C}$ is simple).

We claim that the pair 
\[
(L', S') = (\kg_{\C, \beta_{\max}}, \oplus_{\beta \in \cB(W_{1}) \setminus \{-\beta_{\max}\}} \kg_{\C,\beta})
\]
belongs to $\cT$.
This would lead to a contradiction because for non-zero $u \in \kg_{\C, - \beta_{\max}} \subset W_{1, \C}$, we have
\[
\ad(u)^2 \kg_{\C, \beta_{\max}} = \kg_{\C,  - \beta_{\max}}.
\]

To prove the claim, observe that the set $\cT$ is closed.
Moreover, recall that $\kg_0=\kz_{\kg}(\ka)$ denotes the centralizer of $\ka$ in $\kg$.
The closed connected group $G_0$ of Lie algebra $\kg_0$ and its complexification $G_{0,\C}$, acting on $\kg_\C$ via $\Ad$, preserve $V_{1, \C}, W_{1, \C}$.
Acting on $\Gr(V_{1, \C}, 1) \times \Gr(W_{1, \C}, k - 1)$ diagonally, they preserve the set $\cT$.

Note furthermore that $G_0$ acts irreducibly on $V_1$ and $W_1$.
Indeed, if a non-zero subspace $V'\subseteq V_1$ is $G_0$-invariant, then it is invariant under $G_0 U$, so $\Span \Ad(U^-) V'=\Span \Ad(U^- G_0 U) V'=\Span \Ad(G) V'=\kg$ where the last two equalities use respectively the Zariski-density of $U^- G_0 U$ in $G$, and the irreducibility of the action of $G$.
However, $\Ad(U^-) V'\subseteq V' + F_0 + W_1$, so necessarily $V'=V_{1}$.
The irreducibility of $G_0 \acts W_1$ is similar.

Let $(L,S)$ be a $\R$-point of $\cT$ (whose existence has been established  above).
The irreducibility of $G_0 \acts V_1$ implies the existence of $g_0 \in G_0$ such that
\begin{equation} \label{vhwr1}
g_0 L \not\subset \oplus_{\beta \in \cB(V_1)\setminus\{\beta_{\max}\}} \kg_{\C,\beta},
%\text{the coefficient of $v$ on $\kg_{\C, \alpha_{\max}}$ in the above decomposition is non-zero,}
\end{equation}
and that of $G_0 \acts W_1$ implies the existence of $g_0 \in G_0$ such that
\begin{equation} \label{vhwr2}
\kg_{\C, -\beta_{\max}} \not\subset g_0 S.
%\text{$H_{\C}$ does not contain $\kg_{\C, -\alpha_{\max}}$.}
\end{equation}
for otherwise, $\kg_{\C, -\beta_{\max}}$ would be contained in $\bigcap_{g_0 \in G_0} g_0 S$, a proper $G_0$-invariant subspace of $W_{1,\C}$ defined over $\R$.
As these are Zariski-open conditions in $g_0$ and $G_0$ is connected, there is some $g_0 \in G_0$ for which \eqref{vhwr1} and \eqref{vhwr2} hold simultaneously.
%\marginpar{Leger flou ici: est-ce que $G_{0}$ est algébrique?}
%\comW{Je pense $G_0$ est d'indice fini dans les $\R$-points du centralisateur de $\ka$, qui est algébrique.}

Fix such $g_0$ and consider an element $x$ in the Cartan subalgebra $\kh_{\C}$ of $\kg_{\C}$ such that the eigenvalues $(\beta(x))_{\beta\in \cB}$ are real and $\beta_{\max}(x) > \max_{\beta \in \cB \setminus \{\beta_{\max}\}} \beta(x)$.
For  $t\to+\infty$, observe that $\exp(t x)g_{0}L\to L'$  by \eqref{vhwr1}, and $\exp(t x) g_0 S \to S'$ by \eqref{vhwr2}.
Using that $\cT$ is $G_0$-invariant and closed, we deduce that $(L',W') \in \cT$, as desired.
\end{proof}

\Cref{non-conc-AdgV1} gives a certain class of subspaces that can be put in direct sum with $V_{1}$ modulo the action of $G$. We now work to extend this class.
The next lemma is preparatory to replace the $F_{0}$-component by subspaces $F$ of the form $F=g_{1}V_{m}\cap g_{2}V_{m}$. Recall $\iota$ denotes the longest element in the Weyl group $N_{K}(\ka)/Z_{K}(\ka)$, see \eqref{def-iota}. Below we abusively identify $\iota$ with any representative in $K$. Observe $\iota V_{1}=W_{1}$ and $\iota V_{m}=W_{m}$ (although for $i\neq 1, m$, we may have $\iota V_{i}\neq W_{i}$ depending on $\kappa_{\mu}$).
We also endow $U^-$ with the right invariant Riemmanian metric induced by $\norm{\cdot}_{|\ku^-}$.
(Note the associated distance on $U^-$ is not the distance on $G$ restricted to $U^-$, for $U^-$ is not totally geodesic in $G$).

\begin{lemma}
\label[lemma]{positioning}
Let $g_{1},g_{2}\in G$ and $r\in (0, 1/2]$ such that $$\dang(g_{1} V_{1}, g_{2}V_{m})\geq r.$$ Then there exists $g\in KB^{U^-}_{r^{-C}}$  with $C=C(G)>0$ such that
$$g_{1}\in gP \quad\text{and}\quad g_{2}\in g\iota P.$$
\end{lemma}

\begin{proof} 

We first check
\begin{equation}
\label{eqBPkappaP1}
\setbig{g \,:\, \dang(g V_{1}, W_{m})>0} \subseteq U^- P.
\end{equation}
Indeed, by Bruhat's decomposition: $G=\sqcup_{\omega} P^{-}\omega P$ where $\omega$ varies in the Weyl group of $G$. For $\omega$ different from the identity, we have $ \omega V_{1}\subseteq W_{m}$, whence $gV_{1}\subseteq W_{m}$ for any $g \in P^{-} \omega P$.
It follows that $g\in P^{-} P=U^{-}P$.

We deduce
\begin{equation} \label{eqBPkappaP2}
\setbig{ g \,:\, \dang(g V_{1}, W_{m}) \geq 1/2 } \subseteq B^{U^-}_{R} P.
\end{equation}
for some constant $R>0$ that is large enough depending on $G$ only. 
Indeed, in $G/P$, the left-hand side is compact while the family $(B^{U^-}_{n}P)_{n\geq 1}$ is an increasing sequence of open sets whose union is $U^- P/P$. Hence \eqref{eqBPkappaP2} follows from the previous paragraph.
%For that, note $GV_{1}$ is a manifold that is diffeomorphic to $G/S$ where $S=Stab_{G}(V_{1})$. (it is locally closed because it is compact, equal to $KV_{1}$). As $U\iota P$ is open in $G$, the projection $U\iota V_{1}$ is open in $G/S$. Moreover, the $U$-orbit maps $U \mapsto G.V_{1}, u\mapsto u.x$ are submersions for $x\in U.\iota V_{1}$. Hence given a compact set $\Omega \subseteq U.\iota V_{1}$, there exist $R_{\Omega}>0$, such that $\Omega\subseteq B^U_{R_{\Omega}}\iota P$. Thanks to \eqref{eqBPkappaP1}, we may apply this the case $\Omega=\{g\,:\, \dang(g V_{1}, V_{s-1})\geq 1/2\}$ and obtain \eqref{eqBPkappaP2}.

We now deduce
\begin{equation} \label{eqBPkappaP}
\setbig{g\,:\, \dang(g V_{1}, W_{m})\geq r } \subseteq B^{U^-}_{r^{-C}} P.
\end{equation}
for some $C=C(G)>1$. Note that  \eqref{eqBPkappaP2} justifies \eqref{eqBPkappaP} in the case where $r=1/2$. We infer the general case.
Consider $g$ such that $\dang(g V_{1}, W_{m})\geq r$. Let $v\in \ka^{++}$ satisfying\footnote{The condition  $\alpha(v)=1$ is only used to define $v$ in a deterministic way, so that it does not appear as subscript in the Vinogradov symbols that occur in the proof.} $\alpha(v)=1$ for all simple restricted roots $\alpha\in \Pi$. Taking $t\ggg_{G} |\log r|$, we have
$\dang(\exp(tv) g V_{1}, W_{m})\geq 1/2$. It follows from \eqref{eqBPkappaP2} that $\exp(tv) g\in B^U_{R} P$, i.e.,
$$g\in \exp(-tv) B^{U^-}_{R} P=\exp(-tv) B^{U^-}_{R} \exp(tv) P \subseteq B^{U^-}_{e^{c t}R} P$$
where $c=c(G)>0$.
This justifies \eqref{eqBPkappaP}.

To conclude, write $g_{2}=k_{2}\iota p_{2}$ with $k_{2}\in K$, $p_{2}\in P$ (relying for instance on the Iwasawa decomposition).
Then $g_{2}V_{m}=k_{2}W_{m}$, so the assumption of the lemma means $\dang(g_1 V_{1}, k_{2}W_{m})\geq r $.
It follows from \eqref{eqBPkappaP} that $k_{2}^{-1}g_{1}\in u^-P$ where $u^-\in B^{U^-}_{r^{-C}}$. We set $g= k_{2} u^-$, so $g_{1}\in gP$. Moreover $g_{2}=g (u^{-})^{-1}\iota p_{2}=(g \iota) (\iota^{-1} (u^{-})^{-1}\iota) p_{2} \in g\iota P$.
This concludes the proof.
\end{proof}

Combining \Cref{non-conc-AdgV1} and \Cref{positioning}, we are finally able to show

\begin{proposition}
\label[proposition]{nc-highest-weight}
Consider any subspaces $\R v, H, F\subseteq \kg$ with $\R v$ a line, $\dim H=\dim V_{1}-1$, and $F=g_{1}V_{m}\cap g_{2}V_{m}$ for some $g_{1},g_{2}\in G$. Let $r\in (0, 1/2]$ such that
$$\min\left\{\dang(H, g_{1}V_{m}), \,\dang(g_{1}V_{1}, g_{2}V_{m}) \right\}\geq r.$$
Then there exists $g\in G$ satisfying
$$\dang(gV_{1}, \R v + H + F) \geq r^C$$
where $C=C(G)>1$ is a large enough constant.
\end{proposition}

\begin{proof}
We start with a preliminary observation: By \Cref{non-conc-AdgV1} and compactness of  $\Gr(\kg)$, there exists $\eps_{0}=\eps_{0}(G)>0$ such that for any subspace $S\subseteq \kg$ of the form $S=\R v' +H'+F_{0}$ where $v'\in \kg$,  $\dim H'=\dim W_{1}-1$ and $\dist(H' \too W_{1})<\eps_{0}$, we have
$$\sup_{g\in G}\dang(gV_{1}, S)>\eps_{0}.$$

We now use \Cref{positioning} to reduce to the above observation.
By \Cref{positioning}, there exists $h\in KB^{U^-}_{r^{-O_{G}(1)}}$ such that $g_{1}\in hP$, $g_{2}\in h\iota P$.
Note that $S'':=h^{-1}(\R v+ H+F)= \R v'' + H'' +F_{0}$
where $v''=h^{-1}v$, and $H''=h^{-1}H$ satisfies $\dang(H'', V_{m}) \gg r^{O_{G}(1)}\dang(H,hV_m) \geq r^{O_{G}(1)}$. 
Let $v\in \ka^{++}$ with $\alpha(v)=1$ for all simple restricted roots $\alpha\in \Pi$.
The angle condition on $H''$ implies that for $t\ggg_{G} |\log r|$, we have $\dist(\exp(-tv) H'' \too W_{1})< \eps_{0}$. 
By the first paragraph, we get
$$\sup_{g\in G}\dang(gV_{1}, \exp(-tv) S'')>\eps_{0}.$$
Acting by $h \exp(tv)$, we obtain the claim.
\end{proof}

\subsection{Proof of the supercritical alternative}\label{Sec-supercrit-alt}

We establish a supercritical alternative for the Grassmannian distributions $(\Ad(\theta_{g})V_{i})_{g\sim \mu^n}$ where $i=1,m$. It is presented below as \Cref{sup-dec-X}. The proof utilizes the tools expounded in \S\ref{Sec-nc-random-box}, \S\ref{Sec-back-proj}, \S\ref{Sec-nc-highestweight}. Later, in \S\ref{Sec-proof-mult-Xff}, we will combine this supercritical alternative with the multislicing estimate from \Cref{mult-sup-dec} in order to deduce \Cref{sup-mult-X}.

\begin{proposition}[Supercritical alternative for random walks]
\label[proposition]{sup-dec-X}
Let $\varkappa, \eps, \delta \in (0,1/2)$, let $A\subseteq B^{\kg}_{1}$ be a non-empty subset satisfying for some
 $\alpha\in[\varkappa, 1-\varkappa]$, for $\rho \geq \delta$,
 \begin{equation*}
\max_{v \in \kg} \cN_\delta(A \cap B^\kg_\rho(v)) \leq \delta^{-\eps} \rho^{d \alpha } \cN_\delta(A).
\end{equation*}
Let $t>0$ and $n \geq t |\log \delta|$.

If $\eps, \delta\lll_{\data,\mu, \varkappa, t}1$, then there exists $A'\subseteq A$ such that
 $$\mu^n \setbig{ g \,:\, \Ad(\theta_{g})V_{i} \in \cO^{(\alpha, \eps)}_{\delta}(A')} \leq \delta^\eps. $$
for either $i = 1$ or $i = m$.
\end{proposition}

Recall the exceptional set $\cO^{(\alpha, \eps)}_{\delta}(A)$ defined in Equation \eqref{defcO}.
In terms of the $\SAP$ terminology from \Cref{SAP}, we obtain

%This point of view is better suited for the proof to come. Nonetheless, the case of projectors \emph{parallel to} $\Ad(\theta_{g})V_{i}$ can be deduced by changing the measure $\mu$ as we see below\footnote{in the proof of \Cref{cor-sup-alt-parallelto}.}.

\begin{corollary}
\label[corollary]{cor-sup-alt-parallelto}
Given $\varkappa, t>0$, there exists $\tau'=\tau'(\data, \mu, \varkappa,t)>0$ such that for $\delta\in (0, \tau')$ and $n\geq t |\log \delta|$, the distributions of $(\Ad(\theta_{g})V_{i})_{g\sim \mu^n}$ where $i \in \{1,m\}$ satisfy $\SAP$ with parameters $(\delta, \varkappa, \tau')$.
\end{corollary}

\begin{proof}
We just need to check that \Cref{sup-dec-X} is still valid if we replace $\cO^{(\alpha, \eps)}_{\delta}(A)$ by its dual $\cE^{(\alpha, \eps)}_{\delta}(A)$, which was used to define $\SAP$.
In other terms, we check that \Cref{sup-dec-X} holds for $V_{i}$ replaced by $V_{i}^\perp$.
We may identify $G$ with $\Ad(G)$.
We let $\mu'$ be the image of $\mu$ by the map $g\mapsto \transp{g}^{-1}$ where the left superscript ``$t$'' refers to the adjoint endomorphism of $(\kg, \langle \cdot,\cdot \rangle)$.
Recalling $g=\theta_{g}\exp(\kappa(g))\theta'_{g}$ denotes a Cartan decomposition of $g$ (see \eqref{def-Cartan-decomp}), and using that $\langle \cdot, \cdot \rangle$ is $K$-invariant and every element of $\exp(\ka)$ is self-adjoint,
we have
$$\transp{g}^{-1} =\theta_{g}\exp(-\kappa(g))\theta'_{g}.$$
Note the highest weight subspace of $\exp(-\kappa_{\mu})$ is $W_{1}$, and the orthogonal of its lowest weight subspace is $W_{m}$.
Therefore, applying \Cref{sup-dec-X} to $\mu'$ shows the proposition for $\mu$ is still valid with $(\Ad(\theta_{g})W_{i})_{i=1,m}$ in the place of $(\Ad(\theta_{g})V_{i})_{i=1,m}$. Recalling $V_{i}=W_{m+1-i}^\perp$, this justifies the corollary.
\end{proof}

\begin{proof}[Proof of \Cref{sup-dec-X}] Without loss of generality, we may suppose that $A$ is $2\delta$-separated. We may\footnote{Indeed, if the statement holds for some $\eps_{0}>0$ and every $\delta\in (0, \delta_{0}]$, then it holds automatically for every $0<\eps, \delta\leq \min (\eps_{0},\delta_{0})$.} also allow the upper bound on $\delta$ to depend on $\eps$.
We argue by contradiction assuming that for every subset $B\subseteq A$, for both $i = 1$ and $i = m$,
\begin{equation} \label{eq-ViO}
\mu^n \setbig{ g \,:\, \Ad(\theta_{g})V_{i}\in \cO^{(\alpha, \eps)}_{\delta}(B) } > \delta^\eps.
\end{equation}

We set $d=\dim \kg$, $k=\dim V_{1}$. Given $g\in G$, we write
\[
R_{g}:=\Ad(\theta_{g})V_{1} \quad \text{and} \quad S_{g}:=\Ad(\theta_{g})V_{m}.
\]
We consider $(g_{i})_{i=1, \dots, 4}$ four independent random variables of law $\mu^n$.
The next lemma says that with high probability, there is a large subset $A'$ of $A$ whose projections to $R_{g_{1}}, R_{g_{2}}, S_{g_{3}}, S_{g_{4}}$ are all small.
We may also require that the projection of $A'$ to $R_{g_{1}}$ satisfies some non-concentration properties.

\begin{lemma}
\label[lemma]{gi-proj}
If $\delta\lll_{\data, \mu, t, \eps}1$, then with $(\mu^n)^{\otimes 4}$-probability at least $\delta^{5\eps}$, the variable $(g_{i})_{i=1, \dots, 4}$ satisfies the following. There exists $A'\subseteq A$ such that $|A'|\geq \delta^{5\eps} |A|$ and
$$\max_{i=1,2}\cN_{\delta}(\pi_{R_{g_{i}}}A')<\delta^{-k \alpha -\eps} \quad \text{and} \quad \max_{i=3,4}\cN_{\delta}(\pi_{S_{g_{i}}}A')<\delta^{-(d-k)\alpha- \eps}$$
while for all $\rho\geq \delta$, $y\in R_{g_{1}}$,
$$\cN_{\delta}(\pi_{R_{g_{1}}}A'\cap y^{(\rho)})\leq \delta^{-M \sqrt\eps}\rho^{\varkappa/2}\cN_{\delta}(\pi_{R_{g_{1}}}A')$$
where $M>1$ only depends on $\data$, $\mu$, $t$.
\end{lemma}

\begin{proof}
Applying \eqref{eq-ViO} to $i=m$ and $B=A$, we obtain, with $\mu^n$-probability at least $\delta^\eps$ in $g_{4}$, there exists $A_{4}\subseteq A$ such that $|A_{4}|\geq \delta^\eps |A|$ and $\cN_{\delta}(\pi_{S_{g_{4}}}A_{4})<\delta^{-(d-k)\alpha-\eps}.$
Repeating the argument with $(A, g_{4})$ replaced by $(A_{4}, g_{3})$ we obtain with $\mu^n$-probability at least $\delta^\eps$ in $g_{3}$ some $A_{3}\subseteq A_{4}$ such that $|A_{3}|\geq \delta^\eps |A_{4}|$ and $\cN_{\delta}(\pi_{S_{g_{3}}}A_{3})<\delta^{-(d-k)\alpha-\eps}.$ Similarly, using now \eqref{eq-ViO} in the case $i=1$, and with $(A_{3}, g_{2})$, we obtain with $\mu^n$-probability at least $\delta^\eps$ in $g_{2}$ some $A_{2}\subseteq A_{3}$ such that $|A_{2}|\geq \delta^\eps |A_{3}|$ and $\cN_{\delta}(\pi_{R_{g_{2}}}A_{2})<\delta^{-k\alpha-\eps}.$

For the final step, we need to guarantee both small image and a non-concentration property for the projection to $R_{g_{1}}$.
The previous argument, repeated one more time allows for the first requirement.
Combined with \Cref{nc-proj} (applied with $2\eps$) and \Cref{avoid-bad-pencil-proba} (used at scales $\rho^t \geq \delta^t \geq e^{-n}$),
we obtain (assuming $\delta\lll_{\data, \mu, t, \eps}1$): with $\mu^n$-probability at least $\delta^{2\eps}$ in $g_{1}$, there are subsets $A''_{1} \subseteq A'_{1}\subseteq A_{2}$ satisfying $|A''_{1}|\geq \delta^{\eps} |A'_{1}|\geq \delta^{2\eps} |A_{2}|$ and
$$\cN_{\delta}(\pi_{R_{g_{1}}}A'_{1})<\delta^{-k\alpha-\eps}$$
while for all $\rho \geq \delta$, $y\in R_{g_{1}}$,
$$\cN_{\delta}(\pi_{R_{g_{1}}}A''_{1}\cap y^{(\rho)})\leq \delta^{-M\sqrt \eps}\rho^{\frac{k}{2}\alpha}\cN_{\delta}(\pi_{R_{g_{1}}}A''_{1}).$$
Taking $A':=A''_{1}$ concludes the proof of the lemma.
\end{proof}

The next lemma allows us to choose the spaces $R_{g_{1}}, R_{g_{2}}, S_{g_{3}}, S_{g_{4}}$ with good angle conditions.
\begin{lemma}
\label[lemma]{gi-spatial}
If $\eps\lll_{\data, \mu, t}1$ and $\delta\lll_{\data, \mu, t, \eps}1$, then with $(\mu^n)^{\otimes 4}$-probability at least $1-\delta^{6\eps}$, the variable $(g_{i})_{i=1, \dots, 4}$ satisfies the following for some $M>1$ depending only on $\data, \mu$.
\begin{enumerate}
\item For $i \neq j\in \{1,2,3,4\}$, $\dang(R_{g_{i}}, S_{g_{j}}) \geq \delta^{M\eps}$, 
\item $i =1,2$, $\dang(R_{g_{i}}, \, S_{g_{3}}^\perp+ S_{g_{3}}\cap S_{g_{4}}) \geq \delta^{M\eps},$
\item $\dang(R_{g_{1}}, \, R_{g_{2}}+ S_{g_{3}}\cap S_{g_{4}}) \geq \delta^{M\eps}.$
\end{enumerate}
\end{lemma}

\begin{proof}
Given $h\in G$, we have $\sup_{g\in G}\dang(\Ad(g)V_{1}, S_{h})=1$. It follows from \Cref{non-conc-Ad-thetagVi} that for some $M'=M'(\mu)>1$, and every $\eps\lll_{\mu, t}1$, $\delta\lll_{\data, \mu, t,\eps}1 $,
\begin{equation}\label{gi-spatial1}
\mu^n\setbig{g\,:\, \dang(R_{g}, S_{h}) < \delta^{M'\eps} } <\delta^{7\eps}.
\end{equation}

Moreover, conditionally to the event
%assuming $(g_{3},g_{4})$ satisfy
$$\dang(R_{g_{3}}, S_{g_{4}}) \geq \delta^{M'\eps},$$
\Cref{nc-highest-weight} and the observation that $\dang(S_{g_{3}}^\perp, S_{g_{3}})=1$ together yield
\[
\sup_{g \in G} \dang(\Ad(g) V_{1}, S_{g_{3}}^\perp+ S_{g_{3}}\cap S_{g_{4}})\geq \delta^{CM'\eps}
\]
where $C=C(G)>1$.
Invoking \Cref{non-conc-Ad-thetagVi}, we get that for $M\ggg_{\data,\mu}1$, and $\eps\lll_{\mu, t, M}1$, $\delta\lll_{\data, \mu, t,\eps}1 $,
\begin{equation}\label{gi-spatial2}
\mu^n\setbig{g\, : \,\dang(R_{g}, \, S_{g_{3}}^\perp+ S_{g_{3}}\cap S_{g_{4}} )< \delta^{M\eps}} \leq \delta^{7\eps}.
\end{equation}
Similarly, conditionally to the extra condition 
$\dang(R_{g_{2}}, S_{g_{3}}) \geq \delta^{M'\eps},$
we also get
\begin{equation}\label{gi-spatial3}
\mu^n\setbig{g\, : \,\dang(R_{g}, \, R_{g_{2}}+ S_{g_{3}}\cap S_{g_{4}} )< \delta^{M\eps}} \leq \delta^{7\eps}.
\end{equation}

Equations \eqref{gi-spatial1}, \eqref{gi-spatial2}, \eqref{gi-spatial3} together justify the lemma.
\end{proof}

We now fix, once and for all, a realization of the variables $(g_{i})_{i=1, \dots, 4}$, a subset $A'\subseteq A$, and a constant $M=M(\data, \mu, t)>1$ that satisfy the properties listed in \Cref{gi-proj} and \Cref{gi-spatial}. We set $F:=S_{g_{3}}\cap S_{g_{4}}$ and $E:=F^\perp$ so that
$$\kg=E\oplus^\perp F.$$

\begin{lemma}
\label[lemma]{A''}
Provided $\delta\lll_{\eps}1$, we may further assume $A'\cap E$ contains a subset $A''$ such that
 $$ \cN_{\delta}(A'')\geq \delta^{-2k\alpha +C\eps}.$$
 where $C>1$ only depends on $\data, \mu$.
\end{lemma}

\begin{proof} The first step is to show that $A'$ (or a rather large subset) has small projection to $F$.
Note that by construction of $A'$, we have $$\max_{i=3,4}\cN_{\delta}(\pi_{S_{g_{i}}} A')\leq \delta^{-(d-k)\alpha -\eps}. $$
The non-concentration assumption on $A$, combined with $|A'|\geq \delta^{5\eps} |A|$, also implies
$$\cN_{\delta}(A')\geq \delta^{-d\alpha +6\eps}.$$
On the other hand, the angle condition $\dang(R_{g_{3}}, S_{g_{4}}) \geq \delta^{M\eps}$ implies that every cylinder intersection $(\pi_{S_{g_{3}}}^{-1}B^{S_{g_{3}}}_{\delta}+v)\cap (\pi_{S_{g_{4}}}^{-1}B^{S_{g_{4}}}_{\delta}+v')$ has diameter $O(\delta^{1-M\eps})$, and in particular is covered by at most $O(\delta^{-dM\eps})$ balls of radius $\delta$. Combined with
 the submodular inequality from \cite[Lemma 2.6]{BH24}, we obtain
$$\cN_{\delta}(A') \cN_{\delta}(\pi_{F} A'_{1}) \ll \delta^{-d M\eps} \cN_{\delta}(\pi_{S_{g_{3}}} A')\cN_{\delta}(\pi_{S_{g_{4}}} A')$$
for some $A'_{1}\subseteq A'$ that satisfies $\cN_{\delta}(A'_{1}) \gg \delta^{dM\eps} \cN_{\delta}(A')$. Together, these inequalities imply
\begin{equation}\label{boundpiFA1'}
\cN_{\delta}(\pi_{F} A'_{1})\ll \delta^{-(d-2k)\alpha - (dM+8)\eps}.
\end{equation}

We now extract $A''$ from $A'_{1}$. Equation \eqref{boundpiFA1'}, the general inequality
$$\cN_{\delta}(A'_{1})\leq \cN_{\delta}(\pi_{F} A'_{1})\, \sup_{v\in \kg}\cN_{\delta}(A'_{1}\cap (\pi_{F}^{-1}B_{\delta}+v)) $$
and the lower bound $\delta^{-d\alpha +(dM+6)\eps} \ll \cN_{\delta}(A'_{1})$, together imply the existence of $A''\subseteq A'_{1}$ such that $\pi_{F}(A'')$ is included in a $\delta$-ball and $\cN_{\delta}(A'')\geq \delta^{-2k \alpha + C\eps}$ where $C= (2dM+15)\eps$. Up to translating $A$ and perturbating at scale $\delta$, we can assume $A''\subseteq A\cap E$. This concludes the proof.
\end{proof}

We now aim to show that for most elements $g$ selected by $\mu^n$, the projection of $A''$ to $R_{g}$ has $\delta$-covering number bigger than $\delta^{-k\alpha-\eps}$, yielding a contradiction with our assumption \eqref{eq-ViO}, case $i=1$. To do so, we look at the situation within the subspace $E$, in which we aim to apply the supercritical estimate under mild non-concentration assumptions \Cref{supercrit-weak-assump}. We set
$$L_{g}:=\pi_{E}(R_{g}).$$
The next lemma tells us that the projections of a subset of $E$ to either $R_{g}$ or $L_{g}$ have roughly the same covering numbers, provided $R_{g}$ is not too close to $F$, the orthogonal complement of $E$.

\begin{lemma}
\label[lemma]{redprojE}
Let $g\in G$ with $\dang(R_{g}, F)>r$ for some $r>0$.
For every subset $Z\subseteq E$, we have
\begin{equation}\label{redEeq-cn}
 r^{d}\cN_{\delta}(\pi_{R_{g}}Z) \ll_{\norm{\cdot}} \cN_{\delta}(\pi_{L_{g}}Z) \ll_{\norm{\cdot}} r^{-d} \cN_{\delta}(\pi_{R_{g}}Z)
 \end{equation}
 while for all $y\in L_{g}$, $\rho>0$,
 \begin{equation}\label{redEeq2-nc}
\cN_{\delta}(\pi_{L_{g}}Z \cap y^{(\rho)}) \ll_{\norm{\cdot}} r^{-d} \cN_{\delta}(\pi_{R_{g}}Z \cap (\pi_{R_{g}}y)^{(\rho)}).
 \end{equation}
\end{lemma}

\begin{proof} For the upper bound in \eqref{redEeq-cn}, see for instance \cite[Lemma 18]{He2020JFG}. The proof of the lower bound is similar. To check \eqref{redEeq2-nc}, note first by direct computation\footnote{More precisely, if $v\in E\cap L_{g}^\perp$, then for every $w\in R_{g}$, letting $e_{\perp, w}\in E^\perp$ such that $w+e_{\perp, w}\in L_{g}$, we have $v\perp (w+e_{\perp, w})$, whence $v\perp w$. This justifies one inclusion, and inversing the roles of $L_{g}$ and $R_{g}$, we obtain the converse.} that
$$E\cap L_{g}^\perp= E\cap R_{g}^\perp.$$
In particular for $y\in L_{g}$, $\rho>0$, we have $\pi_{R_{g}}(E \,\cap\, \pi_{L_{g}}^{-1}(y^{(\rho)}))=\pi_{R_{g}}y^{(\rho)}.$
Combined with \eqref{redEeq-cn}, we deduce
\begin{align*}
\cN_{\delta}(\pi_{L_{g}}Z \cap y^{(\rho)}) &= \cN_{\delta}(\pi_{L_{g}}(Z \cap \pi_{L_{g}}^{-1}y^{(\rho)}))\\
& \ll_{\norm{\cdot}} r^{-d}\cN_{\delta}(\pi_{R_{g}}(Z \cap \pi_{L_{g}}^{-1}y^{(\rho)}))\\
& \leq r^{-d}\cN_{\delta}(\pi_{R_{g}}Z \cap \pi_{R_{g}}y^{(\rho)}),
\end{align*}
and \eqref{redEeq2-nc} follows.
\end{proof}

The following lemma allows us to control how close $L_{g}$ is from a subspace $W$ of $\kg$ in terms of the position of $R_{g}$ with respect to $W+F$.
\begin{lemma}
\label[lemma]{Lg-angle}
For $g\in G$ and $W\in \Gr(\kg)$, we have
$$\dang(L_{g}, W) \geq \dang(R_{g}, W+F).$$
\end{lemma}
\begin{proof} It is a consequence of \cite[Lemma 16]{He2020JFG}.
\end{proof}

%\begin{proof}
%We may assume $R_{g}\cap F=\{0\}$. Given a subspace $S\subseteq \kg$, we let $e_{S}\in \bigwedge^*\kg$ denote a unit vector spanning the line $\bigwedge^{\dim S}S$. Recall also every endomorphism of $\kg$ induces an endomorphism of $\bigwedge^*\kg$. We then have
%\begin{align*}
%\dang(L_{g}, W)
%&\geq \dang(L_{g}, W+F) =\frac{\|\pi_{E}(e_{R_{g}})\wedge e_{W+F}\|}{\|\pi_{E}(e_{R_{g}})\| } \\
%&\geq \|e_{R_{g}}\wedge e_{W+F}\| \\
%& =\dang(R_{g}, W+F).
%\end{align*}
%The second line relied on the observations that $\pi_{E}(e_{R_{g}})\wedge e_{W+F}=e_{R_{g}}\wedge e_{W+F}$ and $\|\pi_{E}(e_{R_{g}})\|\leq \|e_{R_{g}}\|=1$.
%\end{proof}

Combining Lemmas \ref{redprojE} and \ref{Lg-angle}, we obtain that the features of $R_{g_{1}}, R_{g_{2}}, A'$ carry over within $E$ to $L_{g_{1}}, L_{g_{2}}, A''$.
\begin{corollary}
\label[corollary]{cor-features-LgiA''}
Provided $\delta\lll_{\data, \eps}1$, we have
$$\dang(L_{g_{1}}, L_{g_{2}}) \geq \delta^{C\eps},$$
$$\max_{i=1,2}\, \cN_{\delta}(\pi_{L_{g_{i}}}A'') \leq \delta^{-k\alpha - C\eps}$$
and for $\rho\geq \delta$, $y\in L_{g_{1}}$,
$$\cN_{\delta}(\pi_{L_{g_{1}}}A''\cap y^{(\rho)})\leq \delta^{-C \sqrt\eps}\rho^{\varkappa/2}\cN_{\delta}(\pi_{L_{g_{1}}}A'')$$
where $C=C(\data, \mu,t)>1$.
\end{corollary}

\begin{proof}
The first inequality (with $C=M$) follows from the lower bound in \Cref{Lg-angle}, in which we plug in the condition that $\dang(R_{g_{1}}, R_{g_{2}}+F)\geq \delta^{M\eps}$ from \Cref{gi-spatial} c).

The combination of Lemmas \ref{redprojE}, \ref{gi-spatial} $a)$, \ref{gi-proj}, yields for $i=1,2$,
$$\cN_{\delta}(\pi_{L_{g_{i}}}A'') \ll_{\data} \delta^{- dM\eps} \cN_{\delta}(\pi_{R_{g_{i}}}A'') \leq \delta^{-k\alpha - (dM+1)\eps },$$
whence the second inequality. It also gives the non-concentration estimate
\begin{equation}\label{ncpiLg1A''1}
\cN_{\delta}(\pi_{L_{g_{1}}}A''\cap y^{(\rho)})\ll_{\data} \delta^{- dM\eps- M\sqrt{\eps}}\rho^{\varkappa/2}\cN_{\delta}(\pi_{R_{g_{1}}}A').
\end{equation}

It remains to bound above $\cN_{\delta}(\pi_{R_{g_{1}}}A')$ using $\cN_{\delta}(\pi_{L_{g_{1}}}A'')$.
On the one hand, by construction of $A'$, we have
\begin{equation}\label{ncpiLg1A''2}
\cN_{\delta}(\pi_{R_{g_{1}}}A')\leq \delta^{-k\alpha-\eps}.
\end{equation}
On the other hand, the angle condition $\dang(L_{g_{1}}, L_{g_{2}}) \geq \delta^{M\eps}$ that we established above yields 
$$\cN_{\delta}(A'')\ll_{\data} \delta^{-M'\eps} \cN_{\delta}(\pi_{L_{g_{1}}}A'')\cN_{\delta}(\pi_{L_{g_{2}}}A'')$$
where $M'=O_{d}(M)$.
As $\delta^{-2k\alpha +C\eps}\leq \cN_{\delta}(A'')$ by construction of $A''$, and $\cN_{\delta}(\pi_{L_{g_{2}}}A'') \ll_{\data} \delta^{-k\alpha - (dM+1)\eps}$ by the above examination, we deduce
\begin{equation}\label{ncpiLg1A''3}
\delta^{-k\alpha + (M'+dM+1+C)\eps}\ll_{\data} \cN_{\delta}(\pi_{L_{g_{1}}}A'').
\end{equation}
Up to increasing $C$, Equations \eqref{ncpiLg1A''1}, \eqref{ncpiLg1A''2}, \eqref{ncpiLg1A''3} together justify the last inequality in the corollary.
\end{proof}

Let $g_5$ be a new independent random variable of law $\mu^n$.
We check that the random subspace $R_{g_5} \in \Gr(\kg, k)$ satisfies the non-concentration assumptions required in \Cref{supercrit-weak-assump} (with respect to the decomposition $E=L_{g_{1}}\oplus L_{g_{2}}$).

\begin{lemma}
\label[lemma]{ncLg5}
Let $W\in \Gr(E,k)$ such that $\max_{i=1,2}\dim (W\cap  L_{g_{i}}) \geq k-1$.
Assume $\delta\lll_{\eps}1$.
Then for every $\rho>\delta$,
$$\mu^n\setbig{g_{5}\,:\,\dang(R_{g_{5}}, W^\perp) \leq \rho }\leq \delta^{-C\eps} \rho^{c} $$
where $C=C(\data, \mu, t)>1$ and $c= c(\mu, t)>0$.
\end{lemma}

The idea of proof is the following. In the case where $W$ is either $L_{g_{1}}$ or $L_{g_{2}}$,  Propositions \ref{non-conc-Ad-thetagVi} and \ref{nc-highest-weight} reduce the estimate to certain angle conditions involving $R_{g_{1}}, R_{g_{2}}, S_{g_{3}}, S_{g_{4}}$, and those  have already been taken care of in \Cref{gi-spatial}. If we only know $\max_{i=1,2}\dim (W\cap  L_{g_{i}}) \geq k-1$, we can still conclude by exploiting  that the condition on the complement $\R v \oplus H$ in \Cref{nc-highest-weight} only concerns the $(k-1)$-dimensional subspace $H$ therein, while the line $\R v$ is free. 

\begin{proof}
%By \Cref{Lg-angle}, we may replace $\dang(L_{g_{5}}, W^\perp\cap E) $ by the quantity $\dang(R_{g_{5}}, W^\perp)$.
By \Cref{non-conc-Ad-thetagVi}, we only need to check the geometric statement
\begin{equation}\label{ncLg5-geom}
\sup_{g\in G}\dang(\Ad(g)V_{1}, W^\perp)\geq \delta^{O_{G}(M\eps)}.
\end{equation}
 Note $W^\perp$ is of the form $W^\perp=(W^\perp\cap E)+ F$ where $(W^\perp\cap E)$ contains a hyperplane of $L_{g_{i}}^\perp \cap E$ for some $i\in \{1,2\}$. By \Cref{nc-highest-weight} and the condition $\dang(R_{g_{3}}, S_{g_{4}}) \geq \delta^{M\eps}$, Equation \eqref{ncLg5-geom}  reduces to showing
\begin{equation}\label{ncLg5-geom2*}
\dang(L_{g_{i}}^\perp \cap E, S_{g_{3}})\geq \delta^{M\eps}.
\end{equation}\label{ncLg5-geom2}
Passing to the orthogonal and applying \Cref{Lg-angle}, we observe that
\begin{align*}
\dang(L_{g_{i}}^\perp \cap E, S_{g_{3}})
&= \dang(L_{g_{i}}+F, \,S_{g_{3}}^\perp )\\
&= \dang(L_{g_{i}}, \,S_{g_{3}}^\perp)\\
&\geq  \dang(R_{g_{i}}, \,S_{g_{3}}^\perp +F).
\end{align*}
The angle condition in \Cref{gi-spatial} b) therefore implies  \eqref{ncLg5-geom2*}, so \eqref{ncLg5-geom} holds. This  concludes the proof.
\end{proof}

\bfparagraph{Conclusion.}
We apply \Cref{supercrit-weak-assump} with the decomposition $E=L_{g_{1}}\oplus L_{g_{2}}$, the random subspace $(R_{g_{5}})_{g_{5}\sim \mu^n}$, the set $A''$, the exponent $C\sqrt \eps$, and the scale $\delta$. The required angle condition on $L_{g_{1}}, L_{g_{2}}$, as well as the the covering numbers conditions for $A''$ are satisfied thanks to \Cref{cor-features-LgiA''}. The non-concentration condition on the random subspace $(R_{g_5})_{g_5 \sim \mu^n}$ is fulfilled thanks to \Cref{ncLg5}.

It follows that, provided $\eps, \delta \lll_{\data, \mu, t, \varkappa} 1$, we have
\[
\mu^n\setbig{g_5 \,:\, R_{g_5}\in \cO_{\delta}^{(\alpha, \eps_0)}(A'')}\leq \delta^{\eps_0},
\]
% $$\mu^n\setbig{g_{5} \,:\, L_{g_{5}}\in \cO_{\delta}^{(\alpha, \eps_{0})}(A'')}\leq \delta^{\eps_{0}}$$
where $\eps_{0}=\eps_{0}(\data, \mu, t,\varkappa)>0$.
%But, thanks to \Cref{redprojE} and the angle condition $\dang(R_{g_{5}}, F)\geq \delta^{M\eps}$, we know that a small projection to $L_{g_{5}}$ yields a small projection to $R_{g_{5}}$. More precisely, provided $\eps\lll_{\data, \mu, t} \eps_{0}$, and $\delta \lll_{\data, \eps}1$, we have
%$$ R_{g_{5}}\in \cO_{\delta}^{(\alpha,\, \eps_{0}-2dM\eps)}(A'') \implies L_{g_{5}}\in \cO_{\delta}^{(\alpha, \eps_{0})}(A'').$$
%We deduce
% $$\mu^n\{g_{5} \,:\, R_{g_{5}}\in \cO_{\delta}^{(\alpha, \,\eps_{0}/2)}(A'')\}\leq \delta^{\eps_{0}}.$$
Such an estimate contradicts \eqref{eq-ViO}. This concludes the proof of the supercritical alternative.
\end{proof}

\subsection{Proof of the supercritical decomposition for random walks} \label{Sec-proof-mult-Xff}

We are finally able to conclude the proof of \Cref{sup-mult-X}. The argument below mimics the final step in the proof of \Cref{RW-preserve-dim}, but plugging the multislicing supercritical decomposition Theorem \ref{mult-sup-dec} instead of the subcritical multislicing Theorem \ref{subc-mult}. The work done in this section until now has been dedicated to establishing \Cref{sup-dec-X}, which is vital to use \Cref{mult-sup-dec}.

\begin{proof}[Proof of \Cref{sup-mult-X}]
Note that it is enough to show the claim with $\delta$ depending on $\eps$ as well. Indeed, if the statement holds for some parameters $(\eps, \delta)$, then it holds for all $(\eps', \delta)$ with $\eps'<\eps$.

By the assumption on $t_1$, $t_2$, we have $0 <  t_{2}- t\lambda_{1}<t_{1}+t\lambda_{1}<2 (t_{2}- t\lambda_{1})$.
We can thus choose $\zeta=\zeta(t_1,t_2,\mu) > 0$ such that $\zeta< t_{2}- t\lambda_{1}<t_{1}+t\lambda_{1}<2\zeta$.
Then $\zeta > (t_{2}-t\lambda_{1})/2$, and therefore $\nu$ is supported on $\{\inj \geq \delta^{\zeta}\}$.
Provided $\delta \lll_{G,(t_{p})} 1$, we apply \Cref{net-charts} with linearizing scale $r=\delta^{\zeta}$.
Consider $\varphi : \{\inj \geq \delta^{\zeta}\}\rightarrow B^\kg_{1}$ the associated map, set $\tnu=\varphi_{\star}\nu$. By \Cref{net-charts} item 1) and the dimension assumption on $\nu$, we have for $\delta\lll_{\data,\eps}1$, for $\rho\in \{\delta^{t_{1}- t\lambda_{i}}\}_{i=1}^{m+1}\cup \{\delta^{t_{2}- t\lambda_{i}}\}_{i=1}^{m+1} \cup [\delta^{t_{1}- t\lambda_{2}}, \delta^{t_{1} - t\lambda_{1}}]$,
 $$\sup_{v\in \kg}\tnu(B^\kg_{\rho}+v)\leq \delta^{-2\eps}\rho^{d\alpha}.$$

We aim to apply \Cref{mult-sup-dec} to $\tnu$ with the localization scale $\eta=1$, the parameters $\ur = (t_1 - t \lambda_i)_{1 \leq i \leq m + 1}$, $\us = (t_2 - t \lambda_i)_{1 \leq i \leq m + 1}$,
the probability space being $\Theta = K$ endowed with the distribution of $\theta_g$, $g \sim \mu^n$,
and the flag $\cV_{\theta} = (\Ad(\theta) V_i)_{1 \leq i \leq m+1}$, so that random boxes are
$$B^{\cV_{\theta}}_{\delta^\ur}:=\Ad(\theta)\sum_{i=1}^{m+1}B^{V_{i}}_{\delta^{t_{1}-t \lambda_{i}}},\qquad B^{\cV_{\theta}}_{\delta^\us}:=\Ad(\theta)\sum_{i=1}^{m+1}B^{V_{i}}_{\delta^{t_{2}-t \lambda_{i}}}.$$
Note that the conditions of application of the theorem are met. Indeed, our choice for $t$ and the observation that the vector $(\lambda_{i})_{i=1}^{m+1}$ is symmetric (i.e. $\lambda_{i}=-\lambda_{m+2-i}$) guarantee that some pair of consecutive entries for $\ur$ and $\us$ coincide, namely $r_{1}=s_{m}$, $r_{2}=s_{m+1}$. Moreover, the hypothesis $t \lambda_{1}<\min(1-t_{1}, t_{2})$ implies that all the exponents $r_{i}, s_{i}$ belong to $(0,1)$.
By \Cref{mu-Sub}, given $C>1$, provided $C\eps \lll_{\mu, t} 1$ and $\delta \lll_{\data, \mu,  C\eps} 1$, we have for each $i=1,\dots, m$ that the distribution of $(\Ad(\theta_{g})V_{i})_{g\sim \mu^n}$ satisfies $\SubP$ with parameters $(\delta^{r_{i+1}}, C\eps, D\sqrt{C\eps})$ where $D=D(\mu)>1$.
Finally, \Cref{cor-sup-alt-parallelto} guarantees that the distributions of $(\Ad(\theta_{g})V_{j})_{g\sim \mu^n}$ where $j=1,m$ together satisfy $\SAP$ with parameters $(\delta^{r_{2}-r_{1}}, \varkappa, \tau')$ where $\tau'=\tau'(\data, \mu, \varkappa, t_1, t_2)>0$ and provided $\delta<\tau'$.

 We may now apply \Cref{mult-sup-dec}. Under the conditions $\eps \lll_{\data, \mu, \varkappa, t_1, t_2}1$ and $\delta \lll_{\data, \mu, \varkappa, t_1, t_2, \eps}1$, we obtain a decomposition
 $$\tnu=\tnu_{1}+\tnu_{2}$$
where $\tnu_{1}=\tnu_{|A_{1}}, \tnu_{2}=\tnu_{|A_{2}}$ for some partition $B^\kg_{1}=A_{1} \sqcup A_{2}$, and a set $E'\subseteq G$ with $\mu^n(E')\leq \delta^\eps$ such that for each $p=1,2$, for $g\in G\smallsetminus E'$, there exists $\tF_{p,g}\subseteq \kg$ with measure $\tnu(\tF_{p,g})\leq \delta^\eps$ and such that
\begin{equation}\label{eq-tnu-tFg}
\sup_{v\in \kg}\tnu_{p |\kg \smallsetminus \tF_{p,g}}\left(B^{\cV_{\theta_{g}}}_{\delta^{\ut_{p}}} +v \right)\leq \leb(B^{\cV_{\theta_{g}}}_{\delta^{\ut_{p}}})^{\alpha+\eps_{0}}
\end{equation}
where $\ut_{p}=\ur$ if $p=1$, $\ut_{p}=\us$ if $p=2$, and $\eps_{0}=\eps_{0}(\data, \mu, \varkappa, (t_{p})_{p})>0$ is fixed.

Moreover, invoking \Cref{descrip-AdgB1}, there exists a subset $E''\subseteq G$ of mass
 $\mu^n(E'')\leq \delta^{\gamma}$ where $\gamma=\gamma(\mu, t_1, t_2, \eps)\in (0, \eps)$ and such that for $\delta\lll_{\mu, t_1, t_2, \eps}1$, for $g\notin E''$, the boxes $\Ad(g)B^\kg_{\delta^{t_{p}}}$ satisfy
\begin{equation}\label{eq-control-AdgB'}
\delta^{\eps} B^{\cV_{\theta_{g}}}_{\delta^{\ut_{p}}}\subseteq \Ad(g)B^\kg_{\delta^{t_{p}}} \subseteq \delta^{-\eps}B^{\cV_{\theta_{g}}}_{\delta^{\ut_{p}}}.
\end{equation}

Put together, \eqref{eq-tnu-tFg} and \eqref{eq-control-AdgB'} yield that for $g\in G\smallsetminus (E'\cup E'')$,
\begin{equation}\label{eq-control-tnuAdgBB}
\sup_{v\in \kg}\tnu_{p |\kg \smallsetminus \tF_{p,g}}\left(\Ad(g)B^\kg_{\delta^{t_{p}}} +v \right)\leq \delta^{-2d \eps} \leb(B^{\cV_{\theta_{g}}}_{\delta^{\ut_{p}}})^{\alpha+\eps_{0}}.
\end{equation}

We now pull back this information to $X$. Choosing $\eps\lll_{\mu, t_1, t_2}1$, we may also suppose 
$$\zeta< t_{2}- t\lambda_{1} - \eps <t_{1}+t\lambda_{1}+\eps<2\zeta.$$
It then follows from \eqref{eq-control-AdgB'} that $B^\kg_{\delta^{2\zeta}}\subseteq \Ad(g)B^\kg_{\delta^{t_{p}}} \subseteq B^\kg_{\delta^{\zeta}}.$
This allows us to apply \Cref{net-charts} item 2), and we obtain from \eqref{eq-control-tnuAdgBB} that for $F_{p,g}:=\varphi^{-1}(\tF_{p,g})\subseteq X$ and $\nu_{p}:=\nu_{ |\varphi^{-1}A_{p}}$, we have $\nu(F_{p,g})\leq \delta^\eps$ and provided $\delta\lll_{\data, \eps}1$,
$$\sup_{x\in X}\nu_{p |X \smallsetminus F_{p,g}}\left(gB^G_{\delta^{t_{p}}}x \right)\leq \delta^{-3d \eps} \leb(B^{\cV_{\theta_{g}}}_{\delta^{\ut_{p}}})^{\alpha+\eps_{0}}. $$
Observing that $\leb(B^{\cV_{\theta_{g}}}_{\delta^{\ut_{p}}}) \simeq_{\data} \delta^{dt_{p}}$ and taking $\eps\lll_{\data, \mu, \varkappa, t_1, t_2} 1$, the upper bound is smaller than $\delta^{dt_{p}(\alpha+\eps_{0}/2)}$. This concludes the proof.
\end{proof}

\appendix
\section{Proof of the multislicing machinery}
\label[appendix]{ss:Appendix-Multislicing}

We establish \Cref{subc-mult} and \Cref{mult-sup-dec}.

\subsection{Measure versus covering number}

We start by observing how measure upper bounds for cells of a partition of $\R^d$ can be deduced from covering number estimates.
Given a partition $\cQ$ of $\R^d$, we write  $\cN_{\cQ}$ for the associated covering number by $\cQ$-cells.
%The next lemma considers a finite measure $\nu$ on $\R^d$, and  expresses that measure upper bounds for $\cR$-cells can be deduced from covering number lower bounds for $\nu$-large subsets of $\R^d$. It will be helpful to deduce our multislicing estimates from covering number lower bounds.

\begin{lemma}
\label[lemma]{slicing-to-measures}
Let $\cQ$ be a partition of $\R^d$. Let $\nu$ be a Borel measure on $B^{\R^d}_{1}$ of mass at most $1$.
Assume that for some constants $C,c>0$, for every subset $A$ such that $\nu(A)\geq c$,  we have
$$\cN_{\cQ}(A)\geq C .$$
Then, there exists $E\subseteq \R^d$ such that $\nu(E)\leq c$ and for every $Q\in \cQ$,
$$\nu_{|\R^d\smallsetminus E}(Q)\leq C^{-1}. $$
\end{lemma}

\begin{proof}
Write $E:=\cup\{Q\in \cQ\,:\, \nu(Q)> C^{-1}\}$. As $\nu$ has mass at most $1$, we have $\cN_{\cQ}(E)< C$. The covering number hypothesis in the lemma then yields $\nu(E)< c$.
\end{proof}

\subsection{Regularization}
\label[appendix]{Sec-regularization}

It will be useful to assume some additional regularity on the measures and sets we will consider. We recall the corresponding notion of regular set, as well as a standard regularization procedure. We also record some weak regularity property for subsets of regular sets.

\bigskip
Given two partitions $\cQ, \cR$ of $\R^d$, we say $\cQ$ is refined by $\cR$, and write $\cQ\prec \cR$, if every $\cQ$-cell is a union of $\cR$-cells. Given $A\subseteq B^{\R^d}_{1}$, we write $\cQ(A)$ the set of $\cQ$-cells meeting $A$. We say that $A$ is \emph{regular} for $\cQ\prec \cR$ if for every $Q\in \cQ(A)$,
$$\cN_{\cR}(A\cap Q)= \frac{\cN_{\cR}(A)}{\cN_{\cQ}(A)}.$$
This notion generalizes to any finite filtration $\cQ_{1}\prec \dots \prec \cQ_{b}$ $(b\geq2)$ by asking regularity for each transition $\cQ_{i} \prec \cQ_{i+1}$.

%Note that for $0\leq r<s$ and $\delta \in (0, 1)$, the partition $\cD_{\delta^s}$ refines $\cD_{\delta^r}$ in the sense that  every $\cD_{\delta^r}$-cell is a union of $\cD_{\delta^s}$-cells. This relation is denoted by $\cD_{\delta^r}\prec \cD_{\delta^s}$. Given $A\subseteq B^{\R^d}_{1}$, we write $\cD_{\delta^s}(A)$ the set $\cD_{\delta^s}$-cells meeting $A$. We say that $A$ is \emph{regular} for $\cD_{\delta^{r}}\prec \cD_{\delta^s}$ if for every $Q\in \cD_{\delta^{s}}(A)$,
%$$\cN_{\cD_{\delta^r}}(A\cap Q)= \frac{\cN_{\cD_{\delta^r}}(A)}{\cN_{\cD_{\delta^s}}(A)}.$$
%This notion generalizes to any finite filtration $\cD_{\delta^{r_{1}}}\prec \dots \prec \cD_{\delta^{r_{m+1}}}$ by asking regularity for each transition $\cD_{\delta^{r_{i}}} \prec \cD_{\delta^{r_{i+1}}}$.

The next lemma allows to decompose any probability measure on $B^{\R^d}_{1}$ as the sum of mutually singular measures which are almost equidistributed among some $\cD_{\delta}$-cells  and whose supports satisfy a prescribed regularity.

%\begin{lemma}[Regularization procedure]\label{reg-measures}
%Let $m\geq1$,  $\ur\in \increasing_{m}$,  $\delta, \eps\in (0, 1/2)$ with $\delta\lll_{d,\eps}1$, and  $\nu$ a Borel probability measure on $B^{\R^d}_{1}$. There is a partition
%$$\R^d = (\bigsqcup_{k\in \cK} A_{k})\sqcup A_{\mathrm{bad}}  $$
%where $\nu(A_{\mathrm{bad}})\leq \delta^{\eps}$,  the index  set $\cK$ is finite of cardinality $|\cK|\ll_{d} |\log \delta|^{O(m)}$ and for each $k\in \cK$,
%\begin{itemize}
%\item  $A_{k}$ is a union of $\cD_{\delta}$-cells and is regular for $\cD_{\delta^{r_{1}}}\prec \dots \prec \cD_{\delta^{r_{m+1}}},$
%\item $\nu(A_{k})\geq \delta^{4\eps}$ and for $Q\in \cD_{\delta}(A_{k})$, we have
%$$2^{-1}\frac{\nu(A_{k})}{\cN_{\cD_{\delta}}(A_{k})}\leq \nu(Q)\leq 2\frac{\nu(A_{k})}{\cN_{\cD_{\delta}}(A_{k})}.$$
%\end{itemize}
%\end{lemma}

\begin{lemma}[Regularization procedure]
\label[lemma]{reg-measures}
Let $b\geq2$,  let $\delta, \eps\in (0, 1/2)$. Let $(\cQ_{i})_{i=1}^b $ be partitions of $\R^d$ such that $\cQ_{1}\prec \dots \prec \cQ_{b} \prec\cD_{\delta}$. Let $\nu$ be a Borel probability measure on $B^{\R^d}_{1}$
If $\delta\lll_{d,\eps}1$, then there is a partition
$$\R^d = (\bigsqcup_{k\in \cK} A_{k})\sqcup A_{\mathrm{bad}}  $$
where $\nu(A_{\mathrm{bad}})\leq \delta^{\eps}$,  the index  set $\cK$ is finite of cardinality $|\cK|\ll_{d} |\log \delta|^{O(b)}$ and for each $k\in \cK$,
\begin{itemize}
\item  $A_{k}$ is a union of $\cD_{\delta}$-cells and is regular for $\cQ_{1}\prec \dots \prec \cQ_{b},$
\item $\nu(A_{k})\geq \delta^{4\eps}$ and for $R\in \cD_{\delta}(A_{k})$, we have
$$2^{-1}\frac{\nu(A_{k})}{\cN_{\cD_{\delta}}(A_{k})}\leq \nu(R)\leq 2\frac{\nu(A_{k})}{\cN_{\cD_{\delta}}(A_{k})}.$$
\end{itemize}
\end{lemma}

\begin{proof}
Let $S=\cup_{R\in \cD_{\delta}(\supp \nu)}R$. Given $j\geq0$, write $S_{j}$  the union of cubes $R\in \cD_{\delta}$ such that
\[
 2^{-j - 1} < \nu(R) \leq 2^{-j}.
\]
In particular, $\nu(S_{j})\ll_{d} 2^{-j}\delta^{-d}$. We can write
$$\nu=\sum_{j\geq 0}\nu_{|S_{j}}=\sum_{j\in \cJ}\nu_{|S_{j}} + \nu_{|S_{\mathrm{bad}}}$$
  where  $\cJ$ is the set of $j\geq 0 $ such that $\nu(S_{j})\geq \delta^{\frac{3}{2}\eps}$ and $S_{\mathrm{bad}}=\cup_{j\notin \cJ} S_{j}$.  We note that  $|\cJ| \leq 10 d |\log \delta|$ and $\nu(S_{\mathrm{bad}})\leq \delta^{\frac{5}{4}\eps}$.

We now decompose each $S_{j}$ where $j\in \cJ$ into regular subsets. More precisely, an  iterated application of \cite[Lemma 2.5]{BH24} allows to write $S_{j}=(\bigsqcup_{\ell \in \cL_{j}} S_{j,\ell})\sqcup S_{j, \mathrm{bad}}$ where $\cN_{\cD_{\delta}}(S_{j, \ell})\geq \delta^{\frac{3}{2}\eps}\cN_{\cD_{\delta}}(S_{j})$, $\cN_{\cD_{\delta}}(S_{j, \mathrm{bad}})\leq \delta^{\frac{5}{4}\eps}\cN_{\cD_{\delta}}(S_{j})$, $|\cL_{j}|\ll_{d} |\log \delta|^{O(b)}$, and each  $S_{j,\ell}$ is a union of $\cD_{\delta}$-cells which is regular from $\cQ_{i}$ to $\cQ_{i+1}$ for all $1\leq i\leq b-1$. Note the fact that $\nu_{|S_{j}}$ is almost equidistributed among $\cD_{\delta}$-cells in $S_{j}$ implies $\nu(S_{j,\ell})\geq \delta^{4\eps}$ and $\nu(S_{j, \mathrm{bad}})\leq 2\delta^{\frac{5}{4}\eps}$.

The proof is concluded by taking $(A_{k})_{k\in \cK}=(S_{j,\ell})_{i\in \cJ, \ell\in \cL}$.
\end{proof}

\subsection{Intrinsic multislicing}

 We need the following statement which generalizes both \cite[Propositions 2.8, 2.9]{BH24} (linear case).

 Given $A\subseteq B^{\R^d}_{1}$, $\eps, \delta>0$, $\tau\in \R$, we set
\begin{equation}
\label{def-excep-set-intrinsic}
\begin{split}
\cI^{\eps, \tau}_{\delta}(A) \defeq \Bigl\{\, V\in \Gr(\R^d) : \exists A' \subseteq A \,\, &\text{ with }\,\,\cN_{\delta}(A')\geq \delta^{\eps} \cN_{\delta}(A) \\
 \text{ and }& \,\cN_{\delta}(\pi_{||V}A) < \delta^{\tau}\cN_{\delta}(A)^{\frac{\dim V^\perp}{d}} \,\Bigr\}.
\end{split}
\end{equation}

Given $i\in \{1, \dots, m\}$, $K\in \cD_{\delta^{r_{i}}}$, we write $V_{K, \theta, i}:=V_{Q(K), \theta, i}$ where $Q(K)\in \cD_{\eta}$ is the unique $\eta$-block containing $K$. Given a box $B^{\uV}_{\delta^\ur}$ and a set $A\subseteq \R^d$,  we write $\cN^{\uV}_{\delta^\ur}(A)$ the covering number of $A$ by translates of $B^{\uV}_{\delta^\ur}$ in $\R^d$.

\begin{proposition}[Intrinsic multislicing for covering numbers]
\label[proposition]{pr:subcritical}
Let $d>m\geq 1$, $\uj\in \cP_{m}(d)$, $\ur \in \increasing_{m}$,  $\delta\in (0, 1)$, $\eta\in [\delta^{r_{1}}, 1]$, $\eps, \eps'>0$, $(\tau_{i})_{i=1,\dots, m}\in \R^m$

Let $(\Theta, \sigma)$ be a probability space. For each $Q\in \cD_{\eta}$, consider a measurable map $\Theta \rightarrow  \cF_{\uj}, \theta \mapsto \uV_{Q,\theta}=(V_{Q,\theta, i})_{i}$.   Let   $A\subseteq B_1^{\R^d}$.

Assume that
\begin{enumerate}
\item  for all $i \in \{1, \dots, m\}$ and $K\in \cD_{\delta^{r_{i}}}$,
$$\sigma\setbig{\theta \,:\, V_{K,\theta, i} \in \cI^{\eps, \tau_{i}}_{\delta^{r_{i+1}}}(A\cap K)} \leq \delta^{r_{i+1}\eps}.$$

\item
 $\cN_{\delta^{r_{m+1}}}(A)\geq \delta^{r_{2}\eps'} \cN_{\delta^{r_{m+1}}}(\widetilde{A})$ for some $\widetilde{A}\subseteq B^{\R^d}_{1}$ containing $A$ and regular with respect to the filtration $\cD_{\delta^{r_1}} \prec \dots \prec \cD_{\delta^{r_{m+1}}}$.
%$A$ is regular with respect to the filtration $\cD_{\delta^{r_1}} \prec \dots \prec \cD_{\delta^{r_{m+1}}}$.
\end{enumerate}

If $\eps'\lll\eps$ and $\delta^{r_{2}} \lll_{d, \eps} 1$, then the exceptional set
\begin{equation*}
\begin{split}
\cE \defeq \Bigl\{\, \theta  \in \Theta : \exists A' \subseteq A & \,\text{ with }\, \cN_{\delta^{r_{m+1}}}(A')\geq \delta^{r_{2}\eps'} \cN_{\delta^{r_{m+1}}}(A)\\
\text{and } & \sum_{Q\in \cD_{\eta}} \cN^{\uV_{Q, \theta}}_{\delta^\ur}(A'\cap Q) < \delta^{\sum_{i=2}^{m+1} (\tau_{i}+\eps)r_{i}} \prod_{i=1}^{m+1} \cN_{\delta^{r_i}}(A)^{j_i/d} \,\Bigr\}
\end{split}
\end{equation*}
has measure $\sigma(\cE) \leq \delta^{r_{2}\eps'}$.
\end{proposition}

\begin{remark}
1) It is necessary here to impose some regularity on $A$. See \cite[Section 2.3]{BH24} for a counterexample when  condition $b)$  is removed.

2) The parameters $\tau_{i}$ may be positive or negative, in which case assumption $a)$ expresses either a subcritical or supercritical estimate.

%3) the term $\delta^{\sum_{i=2}^{m+1} (\tau_{i}+\eps)r_{i}}$ in the conclusion can be equivalently replaced by $\delta^{\sum_{i=2}^{m+1} (\tau_{i}+O(\eps'))r_{i}}$.
\end{remark}

\begin{proof} This is essentially the ouput of the proof of \cite[Proposition 2.8]{BH24} (linear case). We summarize it for completeness, and to help connect with \cite{BH24}. Up to replacing $\delta$ by $\delta^{r_{m+1}}$, and $\ur$ by $r_{m+1}^{-1}\ur$, we may assume $r_{m+1}=1$.
Noting that if the statement is true for some $\eps'$, then it is automatically true for smaller values of $\eps'$, we may assume throughout $\delta^{r_{2}\eps'}\lll_{d}1$. We may also suppose that $\tA$  is $2\delta$-separated, so that $\cN_{\delta}(\tA)=|\tA|$  and similarly for $A$ and other subsets. 
We then distinguish several cases.

 \begin{itemize}
 \item
\noindent \emph{Case $m=1$, $r_{1}=0$}.  Here we  have $\eta=1=\delta^{r_{1}}$, so there are only $O_{d}(1)$-many blocks $Q$ involved in the sum. For each of them, the assumption  $a)$ gives
\begin{equation*}
\begin{split}
\cE_{Q} \defeq \Bigl\{\, \theta  \in \Theta : \exists A'_{Q} \subseteq A\cap Q \,\,&\text{ with }\,\,|A'_{Q}|\geq \delta^{\eps} |A\cap Q|\\
& \text{ and } \,\,\cN^{\uV_{Q, \theta}}_{\delta^\ur}(A'_{Q}) < \delta^{\tau}  |A\cap Q|^{j_{2}/d} \,\Bigr\}
\end{split}
\end{equation*}
satisfies $\sigma(\cE_{Q})\leq \delta^\eps$. Let $\theta \notin \cup_{Q}\cE_{Q}$. Let  $A'\subseteq A$ such that  $|A'|\geq \delta^{\eps'} |A|$. Provided $\delta^{\eps'}\lll_{d}1$, there exists $Q_{0}\in \cD_{1}(A)$ such that $|A'\cap Q_{0}|\geq \delta^{2\eps'} |A|$. Using $\theta \notin \cE_{Q_{0}}$ and taking $\eps'\leq \eps/2$, we get
$$\cN^{\uV_{Q_{0}, \theta}}_{\delta^\ur}(A' \cap Q_{0}) \geq  \delta^{\tau}  |A\cap Q_{0}|^{j_{2}/d}  \geq  \delta^{\tau+\frac{2j_{2}}{d}\eps'}  |A|^{j_{2}/d} $$
and the proof is complete in this case.

\item
\noindent \emph{Case $m=1$, $r_{1}>0$}. Set $\rho:=\delta^{r_{1}}$. We partition $A$ into regular subsets. More precisely, provided $\delta\lll_{d, \eps'}1$, applying \Cref{reg-measures} to the uniform measure on $A$, we may write $A=\sqcup_{i\in I} A_{i} \sqcup A_{\mathrm{bad}}$ where $|I|\ll_{d} |\log \delta|^{O(1)}$, each $A_{i}$ is regular for $\cD_{\rho}\prec\cD_{\delta}$ and satisfies $|A_{i}|\geq \delta^{8\eps'}|A|$, and  $|A_{\mathrm{bad}}|\leq \delta^{2\eps'}|A|$. 
We then subdivide the ball $B^{\R^d}_{1}$ into $\cD_{\rho}$-blocks. Given $i\in I$ and a block $K\in \cD_{\rho}$, we use the shorthand $\cE_{i, K}=\cE_{\delta}^{\eps/2, \tau_{2}}(A_{i}\cap K)$. Note that  $\cE_{i, K}=\emptyset$ if $K\notin \cD_{\rho}(A_{i})$.  On the other hand, if $K\in \cD_{\rho}(A_{i})$, we have  by regularity of $A_{i}$, $\tA$,
$$|A_{i}\cap K|=\frac{|A_{i}|}{\cN_{\rho}(A_{i})}\geq \delta^{9\eps'} \frac{|\tA|}{\cN_{\rho}(\tA)}= \delta^{9\eps'} |\tA\cap K|\geq  \delta^{9\eps'} |A\cap K|.$$
Therefore, in any case, we have $\cE_{i,K}\leq \cE_{\delta}^{\eps, \tau_{2}}(A\cap K)$ (for $\eps'<\eps/18$), which in particular yields $\sigma(\cE_{i,K})\leq \delta^{\eps}$ by assumption a). 
For $\theta\in \Theta$, set $\cK_{\mathrm{bad}}(\theta)=\{K\in \cD_{\rho}(A)\,:\, \theta \in \cup_{i\in I}\cE_{i,K}\}$. 
For $A'\subseteq A$ and $i'\in I$, set
 $\cK_{\mathrm{large}}(A',i')=\{K\in \cD_{\rho}(A_{i'})\,:\, |A'\cap A_{i'}\cap  K| \geq \delta^{\eps/2} |A_{i'}\cap K| \}$. 
 It follows from these definitions and the regularity of $A_{i'}$ that
\begin{align*}
\sum_{Q\in \cD_{\eta}} \cN^{\uV_{Q, \theta}}_{\delta^\ur}(A'\cap Q)
&\gg_{d} \sum_{K\in \cD_{\rho}} \cN^{\uV_{K, \theta}}_{\delta^\ur}(A'\cap A_{i'}\cap K)\\
& \gg_{d} |\cK_{\mathrm{large}}(A',i') \smallsetminus \cK_{\mathrm{bad}}(\theta)| \delta^{\tau_{2}} \cN_{\rho}(A_{i'})^{-j_{2}/d} |A_{i'}|^{j_{2}/d}.
\end{align*}
Via Fubini's theorem and Markov's inequality, the set  $\cE_{1}=\{\theta\,:|\cK_{\mathrm{bad}}(\theta)| \geq \delta^{\eps/4}\cN_{\rho}(A)\}$ satisfies $\sigma(\cE_{1})\leq \delta^{\eps/4}$. Moreover, assuming $|A'|\geq \delta^{\eps'}|A|$ and choosing $i'\in I$ such that $|A'\cap A_{i'}|\geq \delta^{2\eps'}|A_{i'}|$, we have $ |\cK_{\mathrm{large}}(A',i')| \geq \delta^{12\eps'}\cN_{\rho}(A)$ provided $\delta^{\eps'}\lll1$, $\eps'\leq \eps/10$. For $\theta\notin \cE_{1}$,  this leads to 
\begin{align*}
\sum_{Q\in \cD_{\eta}} \cN^{\uV_{Q, \theta}}_{\delta^\ur}(A'\cap Q)
& \gg_{d} \delta^{\tau_{2}+13\eps'} \cN_{\rho}(A_{i'})^{j_{1}/d} |A_{i'}|^{j_{2}/d}.
\end{align*}
Noting $|A_{i'}|\geq \delta^{9\eps'} \abs{\tA}$ by construction, and $\cN_{\rho}(A_{i'})\geq \delta^{9\eps'} \cN_{\rho}(\tA)$  by regularity of $\tA$, the claim follows.

\item
\noindent \emph{Case $m\geq2$}.  Write $\cD^{\uV_{\cD_{\eta}, \theta}}_{\delta^\ur}$ the partition of $\R^d$ obtained from the cell $B^{\uV_{Q, \theta}}_{\delta^\ur}$ in each  block $Q\in \cD_{\eta}$. Write  $\cN^{\uV_{\cD_{\eta}, \theta}}_{\delta^\ur}$ the associated covering number. In particular, one has
$$ \sum_{Q\in \cD_{\eta}} \cN^{\uV_{Q, \theta}}_{\delta^\ur}(A\cap Q) \simeq_{d} \cN^{\uV_{\cD_{\eta}, \theta}}_{\delta^\ur}(A).$$
 Let $A'\subseteq A$ with $|A'|\geq \delta^{r_{2}\eps'} |A|$. Assuming $\delta\lll_{d, r_{2}, \eps'} 1$, we may extract  subsets $A'_{m}\subseteq \dots \subseteq A'_{1}\subseteq A'$ satisfying $|A'_{j}|\geq \delta^{2r_{2}\eps'}|A|$ and with each $A'_{j}$ regular for $\cD_{\delta^{r_{j+1}}} \prec \cD^{\uV_{\cD_{\eta}, \theta}}_{\delta^\ur}\vee \cD_{\delta^{r_{j+1}}}\prec \cD_{\delta}$.  %\marginpar{Pk supposer $\cD^{\uV_{\cD_{\eta}, \theta}}_{\delta^\ur}\vee \cD_{\delta^{r_{j+1}}}\prec \cD_{\delta}$? Ici et dans  note premier papier, ca n'a pas l'air nécessaire, si?}
 A repeated application of the submodular inequality for covering numbers  \cite[Lemma 2.6]{BH24} yields
$$\cN^{\uV_{\cD_{\eta}, \theta}}_{\delta^\ur}(A') \prod_{j=2}^{m}\cN_{\delta^{r_{j}}}(A') \gg_{d}  \prod_{j=1}^{m} \cN^{\uV^{(j)}_{\cD_{\eta}, \theta}}_{\delta^{(r_{j}, r_{j+1})}}(A''_{j})$$
where $\uV^{(j)}_{Q,\theta}=(V_{Q, \theta, j}, \R^d)$ and $A''_{j}\subseteq A'_{j}$ satisfies $\cN_{\delta^{r_{j+1}}}(A''_{j})\gg \cN_{\delta^{r_{j+1}}}(A'_{j})\gg \delta^{3r_{2}\eps'}\cN_{\delta^{r_{j+1}}}(A)$. We may then apply the previous  two cases to bound below the right-hand side, and the proposition follows. See the proof of  \cite[Proposition 2.8]{BH24} for more details.
\end{itemize}
\end{proof}

%The proof proceeds in $3$ steps. First we reduce to the case where $\nu$ is the uniform measure on a $\delta$-separated $A$ which is regular with respect to all transitions of scale $\cD_{\delta^{r_{i}}} \prec\cD_{\delta^{r_{i+1}}}$.  Second, we establish a multislicing supercritical decomposition for $A$ in terms of covering numbers. Third, we use the slicing argument to deduce $\nu$ has the announced property.  Steps $1$ and $3$ are standard. Step $2$ is inspired from [BH24].

\subsection{Proof of the subcritical multislicing}

We establish \Cref{subc-mult}.
For the rest of the section, \emph{we place ourselves in the setting of  \Cref{subc-mult}}.  We will also assume without loss of generality that $r_{m+1}=1$.

\begin{lemma}
\label[lemma]{mult-sup-dec-step1}
In order to prove  \Cref{subc-mult}, we may assume additionally that $\nu$ is the uniform probability measure on a set which is regular with respect to $\cD_{\delta^{r_1}} \prec \dots \prec \cD_{\delta^{r_{m+1}}}$ and intersects each $\cD_{\delta}$-cell in at most one point.
\end{lemma}

\begin{proof} Up to replacing $\nu$ by $\nu/\nu(\R^d)$, we may assume $\nu$ is a probability measure.
Consider the decomposition $\R^d=(\bigsqcup_{k\in \cK} A_{k})\sqcup A_{\mathrm{bad}}$ given by \Cref{reg-measures} applied with $r_{2}\eps'$ in the place of $\eps$. It is enough to establish  \Cref{subc-mult} for each $\nu_{|A_{k}}/\nu(A_{k})$. This in turn reduces to establishing  \Cref{subc-mult} for the probability measure $\cN_{\cD_{\delta}}(A_{k})^{-1}\sum_{Q\in \cD_{\delta}(A_{k})}\delta_{x_{Q}}$ where $x_{Q}$ denotes the center of the cell $Q$. Hence the lemma.
\end{proof}

We hereafter work under the extra assumption of \Cref{mult-sup-dec-step1} and set $A=\supp \nu$.

\begin{lemma}
\label[lemma]{mult-sub-dec-step2}
If $\delta^{r_{2}} \lll_{d, \eps}  1$ and $0<\eps'\lll\eps$, then the exceptional set
\begin{equation*}
\begin{split}
\cE \defeq \Bigl\{\, \theta  \in \Theta : \exists A' \subseteq A \,&\text{ with }\,\cN_{\delta}(A')\geq \delta^{r_{2}\eps'} \cN_{\delta}(A)\\
 \text{and } & \sum_{Q\in \cD_{\eta}} \cN^{\uV_{Q, \theta}}_{\delta^\ur}(A'\cap Q))< \delta^{(\tau+\eps) \sum_{i=2}^{m+1}r_{i}}\prod_{i}t^{-j_{i}/d}_{i} \,\Bigr\},
\end{split}
\end{equation*}
satisfies $\sigma(\cE)\leq \delta^{r_{2}\eps'}$.

\end{lemma}

\begin{proof} Given $i \in \{ 1,\dotsc,m+1 \}$, note that the condition on $t_{i}$ amounts to
$$\sup_{Q \in \cD_{\delta^{r_i}}} |A\cap Q|\ll_d t_{i} |A|.$$
Using the conditions of separation and regularity on $A$, we deduce
\begin{equation*}%\label{eq-dadte}
\cN_{\cD_{\delta^{r_i}}}(A) \gg_d  t^{-1}_{i}.
\end{equation*}
Then \Cref{pr:subcritical} yields the claim.
\end{proof}

\begin{proof}[Proof of \Cref{subc-mult}]
It follows from the combination of \Cref{mult-sub-dec-step2} and \Cref{slicing-to-measures}.
\end{proof}

\subsection{Proof of the supercritical multislicing  decomposition}
We establish \Cref{mult-sup-dec}.
Recall the notation $\cE_\delta^{(\alpha,\tau)}$ defined in \eqref{notation-cEaed}.

\bigskip

We  need the following upgrade on the supercritical alternative property \ref{SAP}.

\begin{lemma}
\label[lemma]{alternative->decomp}
 Let $\sigma_{1},\sigma_{2}$ be  probability measures on $\Gr(\R^d)$, let $\varkappa, \tau, \delta>0$.  Assume $(\sigma_{1},\sigma_{2})$ has the supercritical alternative property $\SAP$  with parameters $(\delta, \varkappa,\tau)$.   Then $(\sigma_{1},\sigma_{2})$ satisfies the following \emph{decomposition property}.

Let $A\subseteq B^{\R^d}_{1}$ be any non-empty subset satisfying for some $\alpha\in[\varkappa, 1-\varkappa]$, for  $\rho \geq \delta$,
\begin{equation}\label{nc-A-tau/4}
\max_{v \in \R^d} \cN_\delta(A \cap B_\rho(v)) \leq \delta^{-\tau/4} \rho^{d \alpha } \cN_\delta(A).
\end{equation}
If $\delta\lll_{d,\tau}1 $, then there exists a decomposition  $A=A_{1}\sqcup A_{2}$ such that
 $$\max_{p=1, 2; A_{p}\neq \emptyset}\sigma_{p}\bigl(\cE^{(\alpha, \tau/4)}_{\delta}(A_{p})\bigr) \leq \delta^{\tau/4}. $$
\end{lemma}

\begin{remark}
This is in fact an equivalence: the above property clearly implies $\SAP$ for $(\sigma_{1}, \sigma_{2})$ with parameters $(\delta,\varkappa, \tau')$ where $\tau'$ only depends on $\tau$.
\end{remark}

\begin{proof}
By the assumption $\delta\lll_{d,\tau}1$, it  is sufficient to prove the above decomposition property when $A$ is $2\delta$-separated and occurences of $\tau/4$ are replaced by $\tau':=\tau/3$.

Applying $\SAP$, we get some subset $S_{1}\subseteq A$ and $p_{1}\in \{1, 2\}$ such that $\sigma_{p_{1}}(\cE^{(\alpha, 3\tau')}_{\delta}(S_{1}))\leq \delta^{3\tau'}$.
If $\abs{A\smallsetminus S_{1}}\geq \delta^{2\tau'}\abs{A}$, then observe that $A\smallsetminus S_{1}$ also satisfies the non-concentration property \eqref{nc-dim-dalpha}.
This allows to  apply  $\SAP$ to $A\smallsetminus S_{1}$, yielding a subset $S_{2}\subseteq A\smallsetminus S_{1}$ and $p_{2}\in \{1, 2\}$ such that $\sigma_{p_{2}}(\cE^{(\alpha, 3\tau')}_{\delta}(S_{2}))\leq \delta^{3\tau'}$.
We can iterate the procedure,  stopping  at the first step $n$ for which  $\abs{A\smallsetminus \cup_{k\leq n} S_{k}} < \delta^{2\tau'}\abs{A}$.
If $\abs{A_{1}}\leq \abs{A_{2}}$, we set $A_{1}=\cup_{k\leq n\,:\,  p_{k}=1}S_{k}$ and $A_{2}=A\smallsetminus A_{1}$. Else we set $A_{2}=\cup_{k\leq n\,:\,  p_{k}=2}S_{k}$ and $A_{1}=A\smallsetminus A_{2}$. Note that in each case the union of the $S_{k}$'s in a given $A_{p}$ occupies a large proportion of $A_{p}$:
$$ |A_{p}\smallsetminus \sqcup_{k}S_{k} |\leq \delta^{\frac{3}{2}\tau'} |A_{p}|. $$
We may now apply \Cref{lm:exhaustion} below with $\Theta = \Gr(\R^d)$, $\lambda$  the counting measure on $A$, $\cP(\theta, A')$  the predicate $\cN_\delta\bigl(\pi_{||\theta} A'\bigr) < \delta^{- \alpha \dim \theta^\perp - 3\tau'}$, the parameter $\rho = \delta^{\tau'}$, and alternatively $(\sigma, A')= (\sigma_p, A_p)$ for $p\in \{1,2\}$.
The claim follows.
\end{proof}

The following lemma, used above, is an abstraction of the exhaustion argument used in the proof of \cite[Proposition 25]{He2020JFG} or \cite[Theorem 2.1]{BH24}.

\begin{lemma}
\label[lemma]{lm:exhaustion}
Let $(\Theta,\sigma)$ be a probability space and $(A,\lambda)$ a finite measure space.
Let $\cP(\theta, A')$ be a predicate with variables $\theta \in \Theta$ and $A' \subset A$.
Assume it is decreasingly monotone in $A'$, in the sense that $\cP(\theta, A')$ implies $\cP(\theta, A'')$ whenever $A'' \subset A'$.
Consider for a measurable subset $A'$ and a parameter $\rho \in (0, 1/16)$, the set
\[
\cE^\cP(A', \rho) = \setbig{ \theta \in \Theta \,:\, \exists A'' \subset A',\, \lambda(A'') \geq \rho \lambda(A') \text{ and } \cP(\theta, A'') }.
\]
If $(S_i)_{i \in I}$ is a finite family of disjoint measurable subsets of $A$, whose union $S :=  \bigsqcup_{i \in I} S_i$ is contained in some $A' \subset A$ with $\lambda(A' \setminus S) \leq \rho^{3/2} \lambda(A')$, then we have
\[
\sigma\bigl( \cE^\cP( A',\rho ) \bigr) \leq 2 \rho^{-1} \sup_{i \in I} \sigma\bigl(\cE^\cP(S_i, \rho^{3/2})\bigr).
\]
\end{lemma}
\begin{proof}
For $i \in I$, let $a_i = \lambda(S_i)/\lambda(S)$ be a weight.
For $J \subset I$, write $a_J = \sum_{j \in J} a_j$.
On account of \cite[Lemma 20]{He2020JFG}, it suffices to show
\[
\cE^\cP(A', \rho) \subset \bigcup_{J : a_J \geq \rho/2} \bigcap_{j \in J} \cE^\cP(S_i, \rho^{3/2}).
\]

Let $\theta \in \cE^\cP(A', \rho)$.
This means there is some $A'' \subset A'$ with $\lambda(A'') \geq \rho \lambda(A')$ such that $\cP(\theta, A'')$ holds.
Consider $J_\theta$ the set of indices $i \in I$ satisfying $\lambda(A'' \cap S_i) \geq \rho^{3/2} \lambda(S_i)$.
By definition and the monotonicity of $\cP$, we have $\theta \in \bigcap_{j \in J_{\theta}} \cE^\cP(S_j, \rho^{3/2})$.
Covering $A''$ by $(S_j)_{j \in J_\theta}$, $(A''\cap S_j)_{j \notin J_\theta}$ and $A' \setminus S$, we obtain
\begin{align*}
\rho \lambda(A') \leq \lambda(A'') &\leq \sum_{j \in J_\theta} \lambda(S_j) + \sum_{j \notin J_\theta} \rho^{3/2} \lambda(S_j) + \lambda(A' \setminus S)\\
& \leq a_{J_\theta} \lambda(S) + \rho^{3/2} \lambda(S) + \rho^{3/2}\lambda(A')
\end{align*}
resulting in $a_{J_\theta} \geq \rho - 2 \rho^{3/2} \geq \rho / 2$.
This shows the desired inclusion.
\end{proof}

We now engage in  the proof of \Cref{mult-sup-dec}.
For the rest of the section, \emph{we place ourselves in the setting of  \Cref{mult-sup-dec}}.
It will be convenient to merge $\ur$ and $\us$ into a tuple $\ut=(t_{k})_{k}\in \increasing_{q}$ defined by
$$\{ r_{1}<\dots < r_{m+1}\}\cup \{ s_{1}<\dots < s_{m+1}\} =\{ t_{1}<\dots < t_{q+1}\}.$$
We will  assume without loss of generality that $t_{q+1}=1$. Recall that by assumption on $\nu$, we have $\nu(\R^d)\leq \delta^{-c}$. Moreover, if $\nu(\R^d)<\delta^\eps$, then the statement is trivially true (taking $A_{p, \theta}=\emptyset$). This allows to assume $\nu(\R^d)\in [\delta^\eps, \delta^{-c}]$. Up to renormalizing, it is enough to establish the statement in the case where $\nu(\R^d)=1$, and with slightly better dimensional gain $\delta^{u \tau/(99d) }$ instead of $\delta^{u \tau/(100d)}$.
From here, similarly to \Cref{mult-sup-dec-step1} we may further assume that $\nu$ is the uniform probability measure on a set $A$ which is regular with respect to $\cD_{\delta^{t_{1}}} \prec\dots \prec \cD_{\delta^{t_{q+1}}}$ and intersects each $\cD_{\delta}$-cell in at most one point. This reduction is valid up to aiming for a slightly better dimensional gain, say $\delta^{u \tau/(98d)}$ instead of $\delta^{u \tau/(99d) }$.

In this context, we establish the following set-theoretic version of \Cref{mult-sup-dec}. We recall the notation $u=r_{i_{1}+1}-r_{i_{1}}$ comes from the statement of \Cref{mult-sup-dec}

\begin{lemma}
\label[lemma]{red-mult-sup-dec-step2}
If $\eps'\lll \eps$, and $\eps, \tau, c \lll_{d, t_{2}, u, \tau'}1$, and $\delta\lll_{d, t_{2}, u,  \tau', \eps} 1$, then there exists a decomposition 
%\marginpar{On demande $\cN_{\delta}(A_{1})>0$ dans $\cE_{1}$ pour avoir $\cE_{1}=\emptyset$ si $A_{1}=\emptyset$}
$$A=A_{1}\sqcup A_{2}$$
such that, writing
\begin{equation*}
\begin{split}
\cE_{1} \defeq \Bigl\{\, \theta  \in \Theta : \exists A' \subseteq A_{1} \,\,&\text{ with }\,\,\cN_{\delta}(A')\geq \delta^{t_{2}\eps'} \cN_{\delta}(A_{1})>0\\
 \text{and } & \sum_{Q\in \cD_{\eta}}\cN^{\uV_{Q, \theta}}_{\delta^\ur}(A'\cap Q) < \delta^{- u \tau'/(90d)} \leb(B^{\uV_{Q, \theta}}_{\delta^\ur})^{-\alpha} \,\Bigr\},
\end{split}
\end{equation*}
we have $\sigma(\cE_{1})\leq \delta^{t_{2}\eps'}$, and similarly with $\cE_{2}$ defined using $(A_{2}, \uW_{Q, \theta}, \us)$ in the place of $(A_{1}, \uV_{Q, \theta}, \ur)$.
\end{lemma}

\begin{proof} We may assume throughout the proof that $\delta$ is small enough depending on $c$, $\eps'$ as well (not only $d, t_{2}, u,  \tau', \eps$). This because if the conclusion is valid for some specific values $c,\eps'$, say depending on $d,t_{2},u,\tau', \eps$, then it also holds for any smaller of values of $c,\eps'$.

Note the non-concentration assumption on $\nu$ amounts to: for  $v\in \R^d$, and $\rho \in \{\delta^{t_{k}}\}_{k=1}^{q+1} \cup  [\delta^{r_{i_{1}+1}}, \delta^{r_{i_{1}}}]$,
$$|A\cap B_{\rho}(v)|\leq \delta^{-c}\rho^{\alpha d} |A|.$$
Taking $\rho = \delta^{t_k}$ and pigeonholing, we find
\begin{equation}\label{eq-dadte}
\forall k \in \{ 1,\dotsc, q+1 \},\quad
\cN_{\delta^{t_{k}}}(A) \gg_{d}  \delta^{-t_{k}\alpha d+c}.
\end{equation}
On the other hand, taking    $\rho \in [\delta^{r_{i_{1} +1}}, \delta^{r_{i_{1}}}]$, and recalling the conditions of separation and regularity on $A$, we find for all $v\in \R^d$,
\begin{equation}\label{reg-all-scales}
\cN_{\delta^{r_{i_{1} +1}}}(A\cap B_{\rho}(v)) \ll_{d} \delta^{-c}\rho^{\alpha d} \cN_{\delta^{r_{i_{1} +1}}}(A).
\end{equation}

Recall $u$ from the statement of \Cref{mult-sup-dec}, and assume  that for some $i\in \{1, \dots, m+1\}$, we have
$$\cN_{\delta^{r_i}}(A) \geq  \delta^{-r_i \alpha d - d^2(\tau+\eps+c)- u\tau'/5}.$$
By \Cref{pr:subcritical} and \eqref{eq-dadte}, if  $\eps'\lll\eps$ and $\delta^{r_{2}}\lll_{d, \eps}1$, then letting $A=A_{1}$, we get $\sigma(\cE_{1})\leq \delta^{r_{2}\eps'}\leq \delta^{t_{2}\eps'}$, thus completing the proof.

We now deal with the case where for every $i\in \{1, \dots, m+1\}$, we have
\begin{equation}\label{eq-smallNrm}
\cN_{\delta^{r_i}}(A) <  \delta^{-r_i \alpha d - d^2(\tau+\eps+c)- u\tau'/5}.
\end{equation}
Recall $i_{1}$ from the statement of \Cref{mult-sup-dec}. Let $K\in \cD_{\delta^{r_{i_{1}}}}(A)$ and let $\Delta^K$ be a similarity sending $K$ onto $[0, 1)^d$. Set   $A_{K}=A\cap K$ and $A^K=\Delta^K A_{K}$.
In order to apply our $\SAP$ hypothesis, we first check a suitable non-concentration property for $A^{K}$. Namely, provided  $\delta \lll_{d,c}1$ , we claim that for all $\rho\geq \delta^u$, $v\in \R^d$
\begin{equation}\label{eq-ncAK}
\cN_{\delta^{u}}(A^{K}\cap B_{\rho}(v))\leq \delta^{-d^2(\tau+\eps+2c)- u\tau'/5}\rho^{d\alpha} \cN_{\delta^{u}}(A^{K}).
\end{equation}
To see why, note that \eqref{eq-ncAK} reduces to: for all $\rho\geq \delta^{r_{i_{1}+1}}$, $v\in \R^d$
\begin{equation}\label{eq-ncA-K}
\cN_{\delta^{r_{i_{1}+1}}}(A_{K}\cap B_{\rho}(v))\ll_{d} \delta^{-c -d^2(\tau+\eps+c)- u\tau'/5}\left(\frac{\rho}{\delta^{r_{i_{1}}}}\right)^{d\alpha} \cN_{\delta^{r_{i_{1}+1}}}(A_{K}).
\end{equation}
To check \eqref{eq-ncA-K}, note we may assume $\rho \in [\delta^{r_{i_{1}+1}}, \delta^{r_{i_{1}}}]$, and replace $A_{K}$ by $A$ on the left handside. Note also from regularity  that $\cN_{\delta^{r_{i_{1}+1}}}(A_{K}) \simeq_{d}\frac{\cN_{\delta^{r_{i_{1}+1}}}(A)}{\cN_{\delta^{r_{i_{1}}}}(A)}$, and then apply  \eqref{reg-all-scales}, \eqref{eq-smallNrm} to conclude. This justifies \eqref{eq-ncAK}.

Set $\sigma^K_{1}$ (resp $\sigma^K_{2}$) the law of $V_{K, \theta, i_{1}}$ (resp. $W_{K, \theta, i_{2}}$) as $\theta\sim \sigma$. By hypothesis, $(\sigma^K_{1}, \sigma^K_{2})$ satisfies $\SAP$ with scale $\delta^{u}$ and parameter $(\varkappa, \tau')$.
Provided $\tau+\eps+2c \leq u\tau'/(20d^2)$ and $\delta^u\lll_{d, \tau'}1$, Equation \eqref{eq-ncAK} allows to use the $\SAP$ assumption (via its upgrade from \Cref{alternative->decomp}) to  obtain  a decomposition:
 $A^{K}=A^{K}_{1}\sqcup A^{K}_{2}$ such that, setting $\tau''=\tau'/4$,
 $$\max_{p=1, 2; A^{K}_{p}\neq \emptyset}\sigma^K_{p}(\cE^{(\alpha, \tau'')}_{\delta^{u}}(A^{K}_{p})) \leq \delta^{u\tau''}.$$
For $p=1$, this means that with $\sigma^K_{1}$-probability at least $1-\delta^{u\tau''}$,  for every  subset ${A_{1}^K}'\subseteq A_{1}^K$ with $\cN_{\delta^u}({A_{1}^K}')\geq \delta^{u \tau''} \cN_{\delta^u}({A_{1}^K})>0$, we have
$$\cN_{\delta^u}(\pi_{|| V_{K,\theta,i_{1}}}{A_{1}^K}')\geq \delta^{-u\alpha \dim V^\perp_{K, \theta, i_{1}}- u\tau''}.  $$
For $p=2$, a similar statement  with $(\sigma^K_{2}, A_{2}^K, W_{K, \theta, i_{2}})$ in the place of $(\sigma^K_{1}, A_{1}^K, V_{K, \theta, i_{1}})$.

We now normalize back from $[0, 1)^d$ to $K$. To this end, note that for any set $S\subseteq [0, 1)^d$, any subspace $V\subseteq \R^d$,
$$\cN_{\delta^u}(\pi_{|| V}S) \simeq_{d} \cN_{\delta^{r_{i_{1}+1}}}(\pi_{|| V}(\Delta^{K})^{-1}S),$$ and note also from \eqref{eq-dadte}, \eqref{eq-smallNrm}, and the regularity of $A$, that
$$\delta^{-u \alpha d} \geq   \delta^{d^2(\tau+\eps+2c)+ u\tau'/5}\frac{\cN_{\delta^{r_{i_{1}+1}}}(A)}{\cN_{\delta^{r_{i_{1}}}}(A)} \simeq_{d} \delta^{d^2(\tau+\eps+2c)+ u\tau'/5}\cN_{\delta^{r_{i_{1}+1}}}(A\cap K).$$
Setting $A_{K, p}=(\Delta^K)^{-1}A^K_{p}$, $u'=u/r_{i_{1}+1}$ , and taking $\tau, \eps, c\lll_{d, u}1$ so that  $d^2(\tau+\eps+2c)+\frac{1}{5}u\tau' \leq \frac{9}{10}u\tau''$, we have
\begin{equation}\label{supercrit-sigmaKut}
\max_{p=1, 2}\sigma^K_{p}(\cI^{u'\tau'', -u' \tau''/20}_{\delta^{r_{i_{1}+1}}}(A_{K, p})) \leq \delta^{r_{i_{1}+1}u'\tau''}.
\end{equation}
(Recall the notation $\cI_{\delta}^{s,t}(A)$ was defined in \eqref{def-excep-set-intrinsic}, and is the empty set  if $A$ is empty).

Set $A^*_{p}=\cup_{K\in \cD_{\delta^{r_{i_{1}}}}}A_{K, p}$ and observe $A=A^*_{1}\sqcup  A^*_{2}$. In particular, at least one  $A^*_{p}$ satisfies $|A^*_{p}|\geq \delta^{t_{2}\eps'} |A|$. For such $p$, \eqref{supercrit-sigmaKut} combined with the $\SubP$ assumption in  \Cref{mult-sup-dec} allows  to apply \Cref{pr:subcritical} to $A^*_{p}$, provided $\eps'\lll \eps\leq u'\tau''$ and $\delta^{t_{2}}\lll_{d, \eps}1$.  Invoking \eqref{eq-dadte} and the regularity of $A$, and assuming $\eps, \tau\lll_{t_{2},u, \tau'}1$, we then obtain  that
the exceptional set $\cE_{p}$ associated to $A^*_{p}$ as in the statement of \Cref{red-mult-sup-dec-step2}  satisfies $\sigma(\cE_{p})\leq \delta^{t_{2}\eps'}$. If both $p=1,2$ satisfy $|A^*_{p}|\geq \delta^{t_{2}\eps'}|A|$, then we set $A_{p}=A^*_{p}$  to finish the proof. If only one $A^*_{p}$ satisfies $|A^*_{p}|\geq \delta^{t_{2}\eps'}|A|$, say $A_{p_{0}}$, we set $A=A_{p_{0}}$ and still obtain the claim with $\eps'$ replaced by $\eps'/2$.
\end{proof}

\begin{proof}[Proof of  \Cref{mult-sup-dec}]
It  follows from the combination of \Cref{red-mult-sup-dec-step2} and  \Cref{slicing-to-measures}.
\end{proof}

\section{Lack of transversality in orthogonal groups}
\label[appendix]{lack}

Set $G=\SO(d,1)$ and $\kg=\so(d,1)$. Fix a Cartan subspace $\ka\subseteq \kg$, a Weyl chamber $\ka^{+}$, and write $V_{1}\subseteq \kg$ the subspace of highest weight. In other terms, $V_{1}$ is characterized by the property that for every $v\in \ka^+$, the adjoint action $\ad(v) \acts \kg$ is by homothety on $V_{1}$ with ratio given by the maximal eigenvalue of $\ad(v)$. 
Observe that $\dim \kg=\frac{d(d+1)}{2}$ and $\dim V_{1}=d-1$.  The goal of \Cref{lack} is to show the following.

\begin{proposition}
\label[proposition]{SO71-count}
Assume $d=2n-1$ with $n\geq2$. Then for any $g_{1},\dots, g_{4}\in G$, the family $(\Ad(g_{i})V_{1})_{i=1, \dots, 4}$ is not in direct sum.  
\end{proposition}

\begin{remark}
We can always put $V_{1}$ and $\Ad(g)V_{1}$ in direct sum, taking $g$ an element in the Weyl group switching $\ka^+$ and $-\ka^+$, whence sending the subspace of highest weight $V_{1}$ to the one of lowest weight.  
\end{remark}

Setting $G_{\C}$ the complexification of $G$, and $V_{1,\C}=V_{1}\otimes \C$, \Cref{SO71-count} is  equivalent to checking that the subspaces $(\Ad(h_{i})V_{1, \C})_{i=1, \dots, 4}$ are   never in direct sum for $(h_{i})_{i}\in G_{\C}^4$.

\bigskip

We first recall a description of the Lie algebra $\kg_{\C}$ of $G_{\C}$ (see    \cite[Section 18]{FH91} for details).  For that, it is convenient to consider the quadratic form on $\C^{2n}$ given by
$$q(x)=\sum_{k=1}^nx_{k}x_{n+k}.$$
It is  represented by the symmetric matrix 
$$
Q:=\frac{1}{2}\begin{pmatrix}
0 & I_{n}\\ I_{n} &0 
\end{pmatrix}
$$
Note $\SO(q, \C)\sim G_{\C}$. The complex Lie algebra $\so(q):=\text{Lie}(\SO(q, \C)) \sim \kg_{\C}$ is then given by
$$\so(q)=\left\{\begin{pmatrix}
A & B\\ C & D
\end{pmatrix} \,:\,  {^tA}=-D, \,{^tB}=-B,\, {^tC}=-C \right\}$$
 where  $A,B,C,D$ run within $M_{n}(\C)$, and the prescript ${^t}$ refers to the transposition. The diagonal matrices in $\so(q)$ constitute a Cartan subalgebra $\kh$ of $\so(q)$ (of rank $n$). Set $H_{k}=E_{k,k}-E_{n+k, n+k}$. The elements $(H_{k})_{1\leq k\leq n}$ form a basis of $\kh$, whose dual basis we denote by $(L_{k})_{1\leq k\leq n}$.  The non-zero roots of $\ad(\kh) \acts \so(q)$ are then given by $\{\pm L_{k}\pm L_{l}\}_{1\leq k\neq  l \leq n}$. More precisely, the root space corresponding to $L_{k}-L_{l}$ is $\C Y_{k,l}$ where $Y_{k,l}=E_{k,l}-E_{n+l,n+k}$, the root space corresponding to $L_{k}+L_{l}$ is $\C Z_{k,l}$ where  $Z_{k,l}=E_{k,n+l}-E_{l,n+k}$, and the rootspace corresponding to $-L_{k}-L_{l}$ is $\C \,{^tZ_{k,l}}$. The roots 
$$\{L_{k}-L_{k+1}\,:\, k=1, \dots, n-1\}\cup \{L_{n-1}+L_{n}\} $$
form a basis of the root system, with set of positive roots given by $\{L_{k} \pm L_{l}\}_{k<l}$.  The associated Borel subalgebra $\kb_{\C}$ then corresponds to the upper triangular matrices in $\so(q)$ (and is parametrized by the upper triangular parts of the blocks $A$ and $B$). 

\bigskip
Let us now describe how $V_{1, \C}$ fits in $\so(q)$. Consider the linear change of variables $\varphi :\C^{2n}\rightarrow \C^{2n}$ characterized by $^t\varphi(e_{1})=e_{1}$, and ${^t\varphi}(e_{n+1})=e_{n+1}$ while  for $l=2, \dots, n$, we set $^t\varphi(e_{l})=\frac{1}{2}(e_{l}+e_{n+l})$ and $^t\varphi(e_{n+l})=\frac{\sqrt{-1}}{2}(e_{l}-e_{n+l})$.
These requirements mean equivalently that, denoting by $(\cdot\,|\,\cdot)$ the standard scalar product on $\C^{2n}$ (namely $(x \,|\,y):= {^tx}\,y$),   we have
$$(\varphi(x)\,|\, e_{1})= (x\,|\, e_{1}), \,\,\,(\varphi(x) \,|\, e_{n+1})= (x \,|\, e_{n+1}),$$
$$ (\varphi(x) \,| \,e_{l})= (x \,|\, \frac{1}{2}(e_{l}+e_{n+l})), \,\,\,(\varphi(x) \,| \,e_{n+l})= (x \,|\,  \frac{\sqrt{-1}}{2}(e_{l}-e_{n+l}))$$
It follows from these relations that
$q\circ \varphi^{-1}(x) = x_{1}x_{n+1} + \sum_{j\neq 1, n+1} x_{j}^2.$ In particular, $\SO(q\circ \varphi^{-1}, \R)\sim \SO(d,1)$. The Cartan subspace of $\so(q\circ \varphi^{-1}, \R)$ is given by 
$$ \ka =\begin{pmatrix}
\diag(t,0\dots, 0) & 0\\ 0 & -\diag(t,0\dots, 0) 
\end{pmatrix}= \R H_{1}$$
The conjugation map $\sC_{\varphi^{-1}}: g\mapsto \varphi^{-1}g\varphi$  sends $\SO(q\circ \varphi^{-1}, \R)$ into $\SO(q, \C)$, and $\so(q\circ \varphi^{-1}, \R)$ into $\so(q)$, thus yielding a real form of $\so(q)$. As $\varphi$ stabilizes $\Vect\{e_{1}, e_{n+1}\}$ and $\Vect\{e_{l}, e_{n+l}\,:\, l=2, \dots, n\}$, it must commute with every element in $\ka$, whence $\sC_{\varphi^{-1}}$ stabilizes $\ka$. This means that $\ka$ is also the Cartan subspace of the real form of $\so(q)$ given by $\sC_{\varphi^{-1}}(\so(q\circ \varphi^{-1}, \R))$.
The corresponding space $V_{1, \C}\subseteq \so(q)$ is then given by the positive eigenspace of $H_{1}$:
$$V_{1, \C}= \Vect_{\C}\set{Y_{1,l}, Z_{1,l} \,:\, 2\leq l\leq n}$$

We now give a more handy description of $V_{1,\C}$.  We write $\langle \cdot,\cdot\rangle_{q}$  the symmetric bilinear form associated to $q$, i.e., $\langle x,y\rangle_{q}=\frac{1}{2}\sum_{l=1}^n x_{l}y_{n+l}+x_{n+l}y_{l}$. Below, the superscript $\perp$ refers to the orthogonal for this bilinear form $\langle \cdot,\cdot\rangle_{q}$.

\begin{lemma}
\label[lemma]{Mw}
The map $V_{1,\C} \rightarrow \{e_{1}, e_{n+1}\}^\perp, M\mapsto 2M e_{n+1}$ is a linear isomorphism.
Given $w\in \{e_{1}, e_{n+1}\}^\perp$, the corresponding $M_{w}\in V_{1,\C}$ is given by: $\forall v\in \C^{2n}$,
$$ M_{w}v= -\langle w, v\rangle_{q} e_{1}+\langle e_{1}, v\rangle_{q} w.$$
\end{lemma}

\begin{proof}
Direct computation justifies that $Me_{n+1}\in \{e_{1}, e_{n+1}\}^\perp=\Vect\set{e_{j}\,:\,  j\neq 1, n+1}$.
Write 
$2M e_{n+1}=\sum_{j\neq 1,n+1} \lambda_{j}e_{j}=:w_{M}$. Observe from the description of $Y_{k,l}$ and $Z_{k,l}$ that $2M e_{j}= -\lambda_{n+j} e_{1}$ where subscripts are considered modulo $2n$. Consider $v=\sum_{j=1}^{2n} c_{j}e_{j}$.  Observe 
$c_{j}=2\langle e_{n+j}, v \rangle_{q}$. 
It follows that
\begin{align*}
Mv &= c_{1}M e_{1}+ c_{n+1}Me_{n+1} + \sum_{j\neq 1, n+1} c_{j} Me_{j}\\
&=0+ \langle e_{1}, v \rangle_{q} w_{M} - \sum_{j\neq 1, n+1} \langle e_{n+j}, v \rangle_{q} \lambda_{n+j} e_{1}\\
&=0+ \langle e_{1}, v \rangle_{q}w_{M} - \langle w_{M}, v \rangle_{q} e_{1}.
\end{align*}
This justifies that the map $V_{1,\C} \rightarrow \{e_{1}, e_{n+1}\}^\perp, M\mapsto 2M e_{n+1}$ is  injective with the desired inverse map. Surjectivity follows because dimensions match.
\end{proof}

The following general fact will also play a role.

\begin{lemma}
\label[lemma]{sym-eq}
Let $(E, \langle \cdot, \cdot\rangle)$ be a $\C$-vector space endowed with a non-degenerate symmetric $\C$-bilinear form.

Let $s\in \N^*$, let $(\eps_{a})_{a=1, \dots, s}, (u_{a})_{a=1, \dots, s}$ be tuples of vectors in $E$ such that the family $(\eps_{a})_{a=1, \dots, s}$ is free. Then the  next two statements are equivalent:
1) for all $v\in E$, 
$$\sum_{a} \langle \eps_{a}, v\rangle u_{a}= \sum_{a} \langle u_{a}, v\rangle \eps_{a}.$$
2) There exists a symmetric matrix $(\lambda_{a,b})_{1\leq a,b\leq s}\in M_{s}(\C)$ such that $u_{a}=\sum_{b}\lambda_{a,b}\eps_{b}$ for every $a$. 
\end{lemma}

\begin{proof}
1) $\implies$ 2). Using freeness of $(\eps_{a})_{a=1, \dots, s}$ and non-degeneracy of $\langle \cdot,\cdot\rangle$,  condition 1) implies that $\Vect\{u_{a}\}_{a=1, \dots, s}\subseteq \Vect\{\eps_{a}\}_{a=1, \dots, s}$. We can now write $u_{a}=\sum_{b}\lambda_{a,b}\eps_{b}$ for some coefficients $\lambda_{a,b}\in \C$. Condition 1) gives the relation 
$$\sum_{a,b} \lambda_{a,b} \langle \eps_{a}, v\rangle  \eps_{b}= \sum_{a,b} \lambda_{a,b}  \langle \eps_{b}, v\rangle \eps_{a}.$$
Using freeness again, we deduce for every $b$ that $\sum_{a}\lambda_{a,b} \eps_{a}= \sum_{a} \lambda_{b,a}  \eps_{a}$, and finally $\lambda_{a,b}=\lambda_{b,a}$.

2) $\implies$ 1). Direct computation.
\end{proof}

We are now able to conclude  that any four translates of $V_{1,\C}$ under $\SO(q,\C)$ are never in direct sum.

\begin{lemma}\label{complex-case}
 For all   $h_{1},h_{2},h_{3},h_{4}\in \SO(q,\C)$, the family 
$(\Ad(h_{a})V_{1,\C})_{a=1, \dots 4}$ is \underline{not} in direct sum.
\end{lemma} 

\begin{proof} It is enough to check the result for a Zariski-dense subset of tuples $(h_{a})_{1\leq a\leq 4}$. In particular, by irreducibility  of the action of $\SO(q,\C)$ on $\C^{2n}$, one may assume 
that the family $(h_{a}e_{1})_{1\leq a\leq 4}$ is free. 

For $a=1, \dots, 4$, let $w_{a}\in \{e_{1}, e_{n+1}\}^\perp$. By \Cref{Mw}, we have the linear relation $\sum_{a}\Ad(h_{a})M_{w_{a}}=0$ is equivalent to
\begin{align}\label{equation-relation-count}
\forall v\in \C^{2n}, \quad \sum_{a=1}^4 \langle h_{a}e_{1}, v\rangle_{q} h_{a}w_{a} = \sum_{a=1}^4 \langle h_{a}w_{a}, v\rangle_{q} h_{a}e_{1}. 
\end{align}
By \Cref{sym-eq}, this amounts to
\begin{align}\label{equation-relation-count2}
\text{$h_{a}w_{a}=\sum_{b=1}^4\lambda_{a,b}h_{b}e_{1}$ for some symmetric matrix $(\lambda_{a,b})_{1\leq a,b\leq 4}$}.
\end{align}
Hence, we are reduced to check the existence of a \emph{non-zero} symmetric matrix $(\lambda_{a,b})_{1\leq a,b\leq 4}\in \Sym_{4}(\C)$ such that for each $a$, the vector 
$\sum_{b}\lambda_{a,b}h_{b}e_{1}$ is  orthogonal to both $h_{a}e_{1}$ and $h_{a}e_{n+1}$ for $\langle \cdot, \cdot\rangle_{q}$.
This last condition defines a subspace of $M_{4}(\C)$ of dimension at least $16-8=8$. On the other hand, $\Sym_{4}(\C)$ has dimension $10$. As $8+10>\dim M_{4}(\C)$, those two subspaces must intersect non trivially.  A non-zero $(\lambda_{a,b})_{1\leq a,b\leq n}$ in the intersection yields via \eqref{equation-relation-count2} an example of (non all zero) $w_{a}$'s  such that 
\[
\sum_{a}\Ad(h_{a})M_{w_{a}}=0.\qedhere
\]
\end{proof}

\begin{proof}[Proof of \Cref{SO71-count}]
It follows from \Cref{complex-case} and the observation that a collection of subspaces of $\R^m$ is  in direct sum if and only if their complexifications is in direct sum in $\C^m$.
\end{proof}

\bibliographystyle{abbrv}
\bibliography{Equidis-G-simple}

@book{Borel,
	Author = {Borel, Armand},
	Edition = {2nd enlarged ed.},
	Fseries = {Graduate Texts in Mathematics},
	Isbn = {0-387-97370-2},
	Issn = {0072-5285},
	Language = {English},
	Publisher = {New York etc.: Springer-Verlag},
	Series = {Grad. Texts Math.},
	Title = {Linear algebraic groups.},
	Volume = {126},
	Year = {1991},
	Zbl = {0726.20030},
	Zbmath = {50185}}

@article{Lin25b,
	Author = {Lin, Zuo},
	Journal = {Preprint arXiv:2511.15696},
	Title = {Polynomially effective equidistribution for certain unipotent subgroups in quotients of perfect {Lie} groups},
	Year = {2025}}

@book{Knapp02,
	Author = {Knapp, Anthony W.},
	Edition = {2nd ed.},
	Fseries = {Progress in Mathematics},
	Isbn = {0-8176-4259-5},
	Issn = {0743-1643},
	Language = {English},
	Publisher = {Boston, MA: Birkh{\"a}user},
	Series = {Prog. Math.},
	Title = {Lie groups beyond an introduction},
	Volume = {140},
	Year = {2002},
	Zbl = {1075.22501},
	Zbmath = {1849081}}

@article{Petridis,
	Author = {Petridis, Giorgis},
	Doi = {10.1007/s00493-012-2818-5},
	Fjournal = {Combinatorica},
	Issn = {0209-9683},
	Journal = {Combinatorica},
	Language = {English},
	Number = {6},
	Pages = {721--733},
	Title = {New proofs of {Pl{\"u}nnecke}-type estimates for product sets in groups},
	Volume = {32},
	Year = {2012},
	Zbl = {1291.11127},
	Zbmath = {6319208}}

@article{BS,
	Author = {Breuillard, Emmanuel and de Saxc{\'e}, Nicolas},
	Doi = {10.4171/JEMS/1346},
	Fjournal = {Journal of the European Mathematical Society (JEMS)},
	Issn = {1435-9855},
	Journal = {J. Eur. Math. Soc. (JEMS)},
	Language = {English},
	Number = {11},
	Pages = {4273--4313},
	Title = {A subspace theorem for manifolds},
	Volume = {26},
	Year = {2024},
	Zbmath = {7927733}}

@article{EMS97,
	Author = {Eskin, A. and Mozes, S. and Shah, N.},
	Doi = {10.1007/PL00001616},
	Fjournal = {Geometric and Functional Analysis. GAFA},
	Issn = {1016-443X},
	Journal = {Geom. Funct. Anal.},
	Language = {English},
	Number = {1},
	Pages = {48--80},
	Title = {Non-divergence of translates of certain algebraic measures},
	Volume = {7},
	Year = {1997},
	Zbl = {0872.22009},
	Zbmath = {1000115}}

@book{Borel91,
	Author = {Borel, Armand},
	Doi = {10.1007/978-1-4612-0941-6},
	Edition = {Second},
	Isbn = {0-387-97370-2},
	Mrclass = {20-01 (20Gxx)},
	Mrnumber = {1102012},
	Mrreviewer = {F. D. Veldkamp},
	Pages = {xii+288},
	Publisher = {Springer-Verlag, New York},
	Series = {Graduate Texts in Mathematics},
	Title = {Linear algebraic groups},
	Url = {https://doi-org-s.qh.yitlink.com:8444/10.1007/978-1-4612-0941-6},
	Volume = {126},
	Year = {1991}}

@article{Lin25,
	Author = {Zuo Lin},
	Journal = {Preprint arXiv:2508.06705},
	Title = {{Quadratic forms of signature $(2,2)$ or $(3,1)$ I: effective equidistribution in quotients of $SL_{4}(\mathbb{R})$}},
	Year = {2025}}

@article{BHZ25,
	Author = {T. B{\'e}nard AND W. He AND H. Zhang},
	Journal = {Preprint arXiv:2508.09076},
	Title = {{Khintchine dichotomy and Schmidt estimates for self-similar measures on $\mathbb{R}^d$}},
	Year = {2025}}

@book{FH91,
	Author = {Fulton, William and Harris, Joe},
	Fseries = {Graduate Texts in Mathematics},
	Isbn = {0-387-97495-4; 0-387-97527-6},
	Issn = {0072-5285},
	Language = {English},
	Publisher = {New York etc.: Springer-Verlag},
	Series = {Grad. Texts Math.},
	Title = {Representation theory. {A} first course},
	Volume = {129},
	Year = {1991},
	Zbl = {0744.22001},
	Zbmath = {51906}}

@article{ABRS,
	Author = {Aka, Menny and Breuillard, Emmanuel and Rosenzweig, Lior and de Saxc{\'e}, Nicolas},
	Doi = {10.1007/s00039-018-0436-0},
	Fjournal = {Geometric and Functional Analysis. GAFA},
	Issn = {1016-443X},
	Journal = {Geom. Funct. Anal.},
	Language = {English},
	Number = {1},
	Pages = {1--57},
	Title = {Diophantine approximation on matrices and {Lie} groups},
	Volume = {28},
	Year = {2018},
	Zbl = {1440.11134},
	Zbmath = {6864558}}

@article{BH25-BrascampLieb,
	Arxiv = {arXiv:2511.11091},
	Author = {B{\'e}nard, Timoth{\'e}e and He, Weikun},
	Howpublished = {Preprint, {arXiv}:2511.11091 [math.{CA}] (2025)},
	Journal = {Preprint arXiv:2511.11091},
	Title = {Effective {Brascamp}-{Lieb} inequalities},
	Url = {https://arxiv.org/abs/2511.11091},
	Year = {2025}}

@article{LMWY25,
	Author = {Lindenstrauss, Elon and Mohammadi, Amir and Wang, Zhiren and Yang, Lei},
	Journal = {arXiv preprint arXiv:2503.21064},
	Title = {Effective equidistribution in rank 2 homogeneous spaces and values of quadratic forms},
	Year = {2025}}

@article{Kim24,
	Author = {Kim, Wooyeon},
	Doi = {10.1215/00127094-2024-0001},
	Fjournal = {Duke Mathematical Journal},
	Issn = {0012-7094},
	Journal = {Duke Math. J.},
	Language = {English},
	Number = {17},
	Pages = {3317--3375},
	Title = {Effective equidistribution of expanding translates in the space of affine lattices},
	Url = {projecteuclid.org/journals/duke-mathematical-journal/volume-173/issue-17/Effective-equidistribution-of-expanding-translates-in-the-space-of-affine/10.1215/00127094-2024-0001.full},
	Volume = {173},
	Year = {2024},
	Zbmath = {7974315}}

@article{Khal-Lue-Wei25,
	Author = {O. Khalil AND M. Luethi AND B/ Weiss},
	Journal = {arXiv:2502.19552},
	Title = {Measure rigidity and equidistribution for fractal carpets},
	Year = {2025}}

@article{BG2008,
	Author = {Bourgain, Jean and Gamburd, Alex},
	Doi = {10.1007/s00222-007-0072-z},
	Fjournal = {Inventiones Mathematicae},
	Issn = {0020-9910},
	Journal = {Invent. Math.},
	Language = {English},
	Number = {1},
	Pages = {83--121},
	Title = {On the spectral gap for finitely-generated subgroups of {{\(\text{SU}(2)\)}}},
	Volume = {171},
	Year = {2008},
	Zbl = {1135.22010},
	Zbmath = {5236752}}

@article{BHZ24,
	Author = {B{\'e}nard, Timoth{\'e}e and He, Weikun and Zhang, Han},
	Journal = {J. Amer. Math. Soc.},
	Pages = {587-623},
	Title = {Khintchine dichotomy for self similar measures},
	Volume = {39},
	Year = {2026}}

@article{BQ-survey12,
	Author = {Benoist, Yves and Quint, Jean-Fran{\c{c}}ois},
	Doi = {10.1007/s11537-012-1220-9},
	Fjournal = {Japanese Journal of Mathematics. 3rd Series},
	Issn = {0289-2316},
	Journal = {Jpn. J. Math. (3)},
	Language = {English},
	Number = {2},
	Pages = {135--166},
	Title = {Introduction to random walks on homogeneous spaces},
	Volume = {7},
	Year = {2012},
	Zbl = {1268.22010},
	Zbmath = {6134767}}

@article{Clozel-Oh-Ullmo,
	Author = {Clozel, Laurent and Oh, Hee and Ullmo, Emmanuel},
	Doi = {10.1007/s002220100126},
	Fjournal = {Inventiones Mathematicae},
	Issn = {0020-9910},
	Journal = {Invent. Math.},
	Language = {English},
	Number = {2},
	Pages = {327--351},
	Title = {Hecke operators and equidistribution of {Hecke} points},
	Volume = {144},
	Year = {2001},
	Zbl = {1144.11301},
	Zbmath = {1655627}}

@article{Chvatal79,
	Author = {Chvatal, V.},
	Doi = {10.1016/0012-365X(79)90084-0},
	Fjournal = {Discrete Mathematics},
	Issn = {0012-365X},
	Journal = {Discrete Math.},
	Language = {English},
	Pages = {285--287},
	Title = {The tail of the hypergeometric distribution},
	Volume = {25},
	Year = {1979},
	Zbl = {0396.60016},
	Zbmath = {3615603}}

@article{BH24,
	Author = {B{\'e}nard, Timoth{\'e}e and He, Weikun},
	Howpublished = {Preprint, {arXiv}:2409.03300},
	Journal = {To appear in Annals of Mathematics},
	Title = {Multislicing and effective equidistribution for random walks on some homogeneous spaces},
	Url = {https://arxiv.org/abs/2409.03300}}

@article{BCCT08,
	Author = {Bennett, Jonathan and Carbery, Anthony and Christ, Michael and Tao, Terence},
	Doi = {10.1007/s00039-007-0619-6},
	Fjournal = {Geometric and Functional Analysis. GAFA},
	Issn = {1016-443X},
	Journal = {Geom. Funct. Anal.},
	Language = {English},
	Number = {5},
	Pages = {1343--1415},
	Title = {The {Brascamp}-{Lieb} inequalities: {Finiteness}, structure and extremals},
	Url = {hdl.handle.net/20.500.11820/b13abfca-453c-4aea-adf6-d7d421cec7a4},
	Volume = {17},
	Year = {2008},
	Zbl = {1132.26006},
	Zbmath = {5258720}}

@article{GMO08,
	Author = {A. Gorodnik AND F. Maucourant AND H. Oh},
	Journal = {Ann. Sci. {Ec. Norm. Super.}},
	Number = {3},
	Pages = {385--437},
	Title = {Manin's and {Peyre}'s conjectures on rational points and adelic mixing},
	Volume = {41},
	Year = {2008}}

@article{Benard23-equidmass,
	Author = {B{\'e}nard, Timoth{\'e}e},
	Doi = {10.1007/s11856-022-2422-3},
	Fjournal = {Israel Journal of Mathematics},
	Issn = {0021-2172},
	Journal = {Isr. J. Math.},
	Language = {English},
	Number = {1},
	Pages = {417--422},
	Title = {Equidistribution of mass for random processes on finite-volume spaces},
	Volume = {255},
	Year = {2023},
	Zbl = {1525.37002},
	Zbmath = {7720427}}

@article{LMW,
	Author = {E. Lindenstrauss AND A. Mohammadi AND Z. Wang},
	Journal = {To appear in Annals of Mathematics},
	Title = {Effective equidistribution for some one parameter unipotent flows}}

@article{He2020JFG,
	Author = {He, Weikun},
	Doi = {10.4171/JFG/92},
	Fjournal = {Journal of Fractal Geometry},
	Issn = {2308-1309},
	Journal = {J. Fractal Geom.},
	Language = {English},
	Number = {3},
	Pages = {271--317},
	Title = {Orthogonal projections of discretized sets},
	Volume = {7},
	Year = {2020},
	Zbl = {1455.28007},
	Zbmath = {7290135}}

@article{BQ1,
	Author = {Benoist, Yves and Quint, Jean-Fran{\c{c}}ois},
	Doi = {10.4007/annals.2011.174.2.8},
	Fjournal = {Annals of Mathematics. Second Series},
	Issn = {0003-486X},
	Journal = {Ann. Math. (2)},
	Language = {English},
	Number = {2},
	Pages = {1111--1162},
	Title = {Stationary measures and closed invariants on homogeneous spaces},
	Volume = {174},
	Year = {2011},
	Zbl = {1241.22007},
	Zbmath = {5960723}}

@article{BQ2,
	Author = {Benoist, Yves and Quint, Jean-Fran{\c{c}}ois},
	Doi = {10.1090/S0894-0347-2013-00760-2},
	Fjournal = {Journal of the American Mathematical Society},
	Issn = {0894-0347},
	Journal = {J. Am. Math. Soc.},
	Language = {English},
	Number = {3},
	Pages = {659--734},
	Title = {Stationary measures and invariant subsets of homogeneous spaces. {II}},
	Volume = {26},
	Year = {2013},
	Zbl = {1268.22011},
	Zbmath = {6168604}}

@article{BQ3,
	Author = {Benoist, Yves and Quint, Jean-Fran{\c{c}}ois},
	Doi = {10.4007/annals.2013.178.3.5},
	Fjournal = {Annals of Mathematics. Second Series},
	Issn = {0003-486X},
	Journal = {Ann. Math. (2)},
	Language = {English},
	Number = {3},
	Pages = {1017--1059},
	Title = {Stationary measures and invariant subsets of homogeneous spaces. {III}.},
	Volume = {178},
	Year = {2013},
	Zbl = {1279.22013},
	Zbmath = {6220728}}

@article{BFLM,
	Author = {Bourgain, Jean and Furman, Alex and Lindenstrauss, Elon and Mozes, Shahar},
	Doi = {10.1090/S0894-0347-2010-00674-1},
	Fjournal = {Journal of the American Mathematical Society},
	Issn = {0894-0347},
	Journal = {J. Am. Math. Soc.},
	Language = {English},
	Number = {1},
	Pages = {231--280},
	Title = {Stationary measures and equidistribution for orbits of nonabelian semigroups on the torus},
	Volume = {24},
	Year = {2011},
	Zbl = {1239.37005},
	Zbmath = {5849079}}

@article{Shmerkin,
	Author = {Shmerkin, Pablo},
	Doi = {10.4171/jems/1283},
	Fjournal = {Journal of the European Mathematical Society (JEMS)},
	Issn = {1435-9855,1435-9863},
	Journal = {J. Eur. Math. Soc. (JEMS)},
	Mrclass = {28A75 (05D99 26A16 28A80 49Q15)},
	Mrnumber = {4634691},
	Number = {10},
	Pages = {4155--4204},
	Title = {{A non-linear version of Bourgain's projection theorem}},
	Url = {https://doi.org/10.4171/jems/1283},
	Volume = {25},
	Year = {2023}}

@book{BQ_book,
	Author = {Benoist, Yves and Quint, Jean-Fran{\c c}ois},
	Doi = {10.1007/978-3-319-47721-3},
	Fseries = {Ergebnisse der Mathematik und ihrer Grenzgebiete. 3. Folge},
	Isbn = {978-3-319-47719-0; 978-3-319-47721-3},
	Issn = {0071-1136},
	Language = {English},
	Publisher = {Cham: Springer},
	Series = {Ergeb. Math. Grenzgeb., 3. Folge},
	Title = {Random walks on reductive groups},
	Volume = {62},
	Year = {2016},
	Zbl = {1366.60002},
	Zbmath = {6644704}}

@article{He2020IJM,
	Author = {He, Weikun},
	Doi = {10.1007/s11856-020-2032-x},
	Fjournal = {Israel Journal of Mathematics},
	Issn = {0021-2172},
	Journal = {Isr. J. Math.},
	Language = {English},
	Number = {2},
	Pages = {593--627},
	Title = {Random walks on linear groups satisfying a {Schubert} condition},
	Volume = {238},
	Year = {2020},
	Zbl = {1452.60023},
	Zbmath = {7247266}}

@article{Bourgain2010,
	Author = {Bourgain, Jean},
	Doi = {10.1007/s11854-010-0028-x},
	Fjournal = {Journal d'Analyse Math{\'e}matique},
	Issn = {0021-7670},
	Journal = {J. Anal. Math.},
	Language = {English},
	Pages = {193--236},
	Title = {The discretized sum-product and projection theorems},
	Volume = {112},
	Year = {2010},
	Zbl = {1234.11012},
	Zbmath = {5864352}}

@article{HS2022,
	Author = {He, Weikun and de Saxc{\'e}, Nicolas},
	Doi = {10.1215/00127094-2021-0045},
	Fjournal = {Duke Mathematical Journal},
	Issn = {0012-7094},
	Journal = {Duke Math. J.},
	Language = {English},
	Number = {5},
	Pages = {1061--1133},
	Title = {Linear random walks on the torus},
	Volume = {171},
	Year = {2022},
	Zbl = {1518.37005},
	Zbmath = {7513356}}

@article{HLL_Aff,
	Author = {He, Weikun and Lakrec, Tsviqa and Lindenstrauss, Elon},
	Doi = {10.1093/imrn/rnaa322},
	Fjournal = {IMRN. International Mathematics Research Notices},
	Issn = {1073-7928},
	Journal = {Int. Math. Res. Not.},
	Language = {English},
	Number = {11},
	Pages = {8003--8037},
	Title = {Affine random walks on the torus},
	Volume = {2022},
	Year = {2022},
	Zbl = {1503.37017},
	Zbmath = {7560362}}

@incollection{HLL_Nil,
	Author = {He, Weikun and Lakrec, Tsviqa and Lindenstrauss, Elon},
	Booktitle = {Analysis at large. Dedicated to the life and work of Jean Bourgain},
	Doi = {10.1007/978-3-031-05331-3_7},
	Isbn = {978-3-031-05330-6; 978-3-031-05333-7; 978-3-031-05331-3},
	Language = {English},
	Pages = {131--171},
	Publisher = {Cham: Springer},
	Title = {Equidistribution of affine random walks on some nilmanifolds},
	Year = {2022},
	Zbl = {1518.58018},
	Zbmath = {7672487}}

@article{BenoistSaxce,
	Author = {Benoist, Yves and de Saxc{\'e}, Nicolas},
	Doi = {10.1007/s00222-015-0636-2},
	Fjournal = {Inventiones Mathematicae},
	Issn = {0020-9910},
	Journal = {Invent. Math.},
	Language = {English},
	Number = {2},
	Pages = {337--361},
	Title = {A spectral gap theorem in simple {Lie} groups},
	Volume = {205},
	Year = {2016},
	Zbl = {1357.22003},
	Zbmath = {6641518}}

@article{BG2012,
	Author = {Bourgain, Jean and Gamburd, Alex},
	Doi = {10.4171/JEMS/337},
	Fjournal = {Journal of the European Mathematical Society (JEMS)},
	Issn = {1435-9855},
	Journal = {J. Eur. Math. Soc. (JEMS)},
	Language = {English},
	Number = {5},
	Pages = {1455--1511},
	Title = {A spectral gap theorem in {SU}{{\((d)\)}}},
	Volume = {14},
	Year = {2012},
	Zbl = {1254.43010},
	Zbmath = {6095887}}

@article{HS2023,
	Author = {He, Weikun and De Saxc{\'e}, Nicolas},
	Doi = {10.1017/etds.2025.3},
	Fjournal = {Ergodic Theory and Dynamical Systems},
	Issn = {0143-3857},
	Journal = {Ergodic Theory Dyn. Syst.},
	Language = {English},
	Number = {10},
	Pages = {3090--3147},
	Title = {Semisimple random walks on the torus},
	Volume = {45},
	Year = {2025},
	Zbmath = {8090326}}

@article{LM,
	Author = {Lindenstrauss, E. and Mohammadi, A.},
	Doi = {10.1007/s00222-022-01162-5},
	Fjournal = {Inventiones Mathematicae},
	Issn = {0020-9910},
	Journal = {Invent. Math.},
	Language = {English},
	Number = {3},
	Pages = {1141--1237},
	Title = {Polynomial effective density in quotients of {{\(\mathbb{H}^3\)}} and {{\(\mathbb{H}^2 \times \mathbb{H}^2\)}}},
	Volume = {231},
	Year = {2023},
	Zbl = {1514.37010},
	Zbmath = {7661776}}

@article{Yang,
	Author = {Yang, Lei},
	Doi = {10.4007/annals.2025.202.1.3},
	Fjournal = {Annals of Mathematics. Second Series},
	Issn = {0003-486X},
	Journal = {Ann. Math. (2)},
	Language = {English},
	Number = {1},
	Pages = {189--264},
	Title = {Effective version of {Ratner}'s equidistribution theorem for {{\(\mathrm{SL}(3,\mathbb{R})\)}}},
	Volume = {202},
	Year = {2025},
	Zbmath = {8063318}}

@article{KhalilLuethi,
	Author = {Khalil, Osama and Luethi, Manuel},
	Doi = {10.1007/s00222-022-01171-4},
	Fjournal = {Inventiones Mathematicae},
	Issn = {0020-9910},
	Journal = {Invent. Math.},
	Language = {English},
	Number = {2},
	Pages = {713--831},
	Title = {Random walks, spectral gaps, and {Khintchine}'s theorem on fractals},
	Volume = {232},
	Year = {2023},
	Zbl = {1526.11040},
	Zbmath = {7676259}}

@article{Gan24,
	Author = {Gan, Shengwen},
	Doi = {10.1093/imrn/rnae008},
	Eprint = {https://academic.oup.com/imrn/article-pdf/2024/9/7944/57442834/rnae008.pdf},
	Issn = {1073-7928},
	Journal = {International Mathematics Research Notices},
	Month = {02},
	Number = {9},
	Pages = {7944-7971},
	Title = {{Exceptional Set Estimate Through Brascamp--Lieb Inequality}},
	Url = {https://doi.org/10.1093/imrn/rnae008},
	Volume = {2024},
	Year = {2024}}

@article{He2019,
	Author = {He, Weikun},
	Doi = {10.1007/s11854-019-0071-1},
	Fjournal = {Journal d'Analyse Math{\'e}matique},
	Issn = {0021-7670},
	Journal = {J. Anal. Math.},
	Language = {English},
	Number = {2},
	Pages = {637--676},
	Title = {Discretized sum-product estimates in matrix algebras},
	Volume = {139},
	Year = {2019},
	Zbl = {1453.11019},
	Zbmath = {7172359}}

\end{document}